\documentclass{SCAEOL}
\numberwithin{equation}{section}
\usepackage{subfigure}
\usepackage{graphicx}
\usepackage{color}
\begin{document}

\Year{2013} %
\Month{January}
\Vol{56} %
\No{1} %
\BeginPage{1} %
\EndPage{XX} %
\AuthorMark{Yong-Xia Hao {\it et al.}}

\newtheorem{discussion}{Discussion}
\newtheorem{algorithm}{Algorithm}

\title{Sparse approximate solution of fitting surface to scattered points by MLASSO model}{}

%

\author[1]{Yong-Xia Hao}{}
\author[2]{Chong-Jun Li}{Corresponding author}
\author[2]{Ren-Hong Wang}{}

\address[{\rm1}]{Faculty of Science, Jiangsu University, Zhenjiang {\rm 212000}, China;}
\address[{\rm2}]{School of Mathematical Sciences, Dalian University of Technology, Dalian {\rm 116024}, China;}
\Emails {yongxiahaoujs@ujs.edu.cn; chongjun@dlut.edu.cn; renhong@dlut.edu.cn}\maketitle


 {\begin{center}
\parbox{14.5cm}{\begin{abstract}
The goal of this paper is to achieve a computational model and corresponding efficient
algorithm for obtaining a sparse representation of the fitting surface to the given scattered data.
The basic idea of the model is to utilize the principal shift invariant (PSI) space and
the $l_{1}$ norm minimization. In order to obtain different sparsity of the approximation solution,
the problem is represented as a multilevel LASSO (MLASSO) model with different regularization parameters.
The MLASSO model can be solved efficiently by the alternating direction method of multipliers.
Numerical experiments indicate that compared to the AGLASSO model and the basic MBA algorithm in \cite{Lee},
the MLASSO model can provide an acceptable compromise between the minimization of the data mismatch term and the sparsity of the solution.
Moreover, the solution by the MLASSO model can reflect the regions of the underlying surface where high gradients occur.\vspace{-3mm}
\end{abstract}}\end{center}}



\renewcommand{\baselinestretch}{1.2}
\begin{center} \renewcommand{\arraystretch}{1.5}
{\begin{tabular}{lp{0.8\textwidth}} \hline \scriptsize
{\bf Citation:}\!\!\!\!&\scriptsize Y-X Hao, C-J Li, R-H Wang. SCIENCE CHINA Mathematics  journal sample. Sci China Math, 2013, 56, doi: 10.1007/s11425-000-0000-0\vspace{1mm}
\\
\hline
\end{tabular}}\end{center}

\baselineskip 11pt\parindent=10.8pt  \wuhao
\section{Introduction}
Sparse representation of a function via a linear combination of a small number of
functions has recently received a lot of attention in several mathematical fields such as
approximation theory \cite{Cohen,Kunis,Rauhut1,Rauhut2}, compressed sensing, signal and image
processing \cite{Cand1,Cand2,Cand3,Cand4} etc.
The problem can be described as follows. Consider a linearly dependent set of $n$ functions
$\{\varphi_{i}\}_{i=1}^{n}$ and a function $f$ represented as $f = \sum_{i=1}^{n}X_{i}\varphi_{i}.$
Since the set of functions is not linearly independent, this representation is not unique.
The problem is then to find the sparsest solution, i.e., the coefficient vector $X=(X_{1},X_{2},\ldots,X_{n})$ has
as many zero components as possible (referred to minimizing the $l_{0}$ norm of the vector $X$). This
optimization problem is NP-hard, since the $l_{0}$ norm is nonconvex and discontinuous.
Hence, much attention has been paid to solutions minimizing $\|X\|_{1}= \sum_{i=1}^{n}|X_{i}|$ instead.

In this paper, we consider the problem of reconstructing a surface from scattered data using a sparse representation.
The scattered data fitting problem arises in many
applications, such as signal processing, computer graphics and neural networks \cite{Johnson}. In a typical
scattered data reconstruction problem, we are given a set of scattered points
$\Xi = \{x_{1},x_{2},\ldots,x_{N}\}\subseteq \mathbb{R}^{2}$ and associated noisy function values
$\tilde{f}=f|_{\Xi}+\textbf{n}=\{f_{1},f_{2},\ldots,f_{N}\}$, where $\textbf{n}$ is the error vector.
Then we seek a function $g$ which fits the given
data $\{(x_{i},f_{i})\}_{i=1}^{N}$ well. There are a lot of existing methods and algorithms in the literature.
Various methods can be found in a survey on scattered data interpolation \cite{Bozzini}. For approximation methods, B-splines
have a solid mathematical foundation and have been used in many literatures, such as \cite{Lee2,Wang3} etc.
Wavelet frames have also been used to reconstruct implicit surfaces from unorganized point sets in $\mathbf{R}^{3}$ \cite{Dong2}.
In order to control the local and global fitting error simultaneously, adaptive methods are presented in \cite{Casta1,Schmitt}.
The adaptivity is achieved by a portion of the data with a patch, testing the fit for satisfaction
within a given tolerance, and subdividing the patch if the tolerance is not met \cite{Forsey}.
In addition, several approximation methods employ a multilevel structure to approximate data efficiently.
In particular, a multilevel scheme based on B-splines is proposed in \cite{Lee}
to approximate scattered data.  These methods run on the approximation space $S=\bigcup_{j=1}^{J}S_{j}$,
where $S_{j}\subseteq S_{j+1}$ are  principal shift invariant (PSI) spaces generated by a
single compactly supported function $\varphi$. The multilevel approximation procedure is as follows: for each level $j$,
the point set $(\Xi, \triangle g_{j-1})$ is approximated by a function $g_{j}\in S_{j}$ obtained by the least square method, where
$\triangle g_{j-1}|_{\Xi} = \tilde{f}-(g_{1} + \cdots + g_{j-1})|_{\Xi}$. The procedure is terminated
until certain conditions are satisfied. Then the final approximation surface is
$$g = g_{1} + g_{2} + \cdots + g_{J}, \; g_{j}\in S_{j}, \; j=1,\ldots,J.$$
However, the methods above do not produce the sparse representation of the surface.

In this paper, we present an efficient method to obtain a sparse representation
of the fitting surface to the given scattered points. We still choose the space $S$ defined above as the approximation space.
But instead of using the multilevel scheme, we put all the basis functions of $S_{j}, 1\leq j\leq J$ together as a
frame of $S$. Denote the basis functions of $S_{j}$ as $\{\varphi_{i}^{j}\}_{i=1}^{n_{j}}, 1\leq j\leq J$,
then $S=span\{\varphi_{i}^{j}\}_{i,j=1}^{n_{j},J}$. We then try to find the fitting surface $g\in S$ as:
$$g=\sum_{j=1}^{J}g_{j}, \; g_{j}=\sum_{i=1}^{n_{j}}X_{i}^{j}\varphi_{i}^{j}.$$
Since the functions $\{\varphi_{i}^{j}\}_{i,j=1}^{n_{j},J}$ are linearly dependent, the representation
of $g$ as above is not unique and we will seek a relatively sparse one. The choice of the space $S$ makes a
sparse representation of $g$ exist and the function $\varphi_{i}^{j}$ can be constructed in
a multilevel way. We use the similar approach as those used in compressed sensing, i.e., to use
the $l_{1}$ norm of the coefficient vector as the regularization term.
Thus, the problem can be represented by the following minimization
\begin{equation}\label{eq:1.1}
   \min_{g\in S} \sum_{i=1}^{N}(g(x_{i})-f_{i})^{2}+\sum_{j=1}^{J}\lambda_{j}\|\textbf{X}_{j}\|_{1},
\end{equation}
where $\mathbf{X}_{j}=(X_{1}^{j},\ldots,X_{n_{j}}^{j}), 1\leq j\leq J$ are coefficient vectors.
The parameters $\lambda_{j}> 0, j=1,\ldots,J$ are called the regularization parameters,
which serve as a weight to adjust the balance between the two terms. Large values of $\lambda_{j}$
will lead to a sparse function $g$, at the cost of a potentially large fitting error, while
small values of $\lambda_{j}$ will lead to a small fitting error, but with a potentially not too sparse fitting
function $g$. In addition, different values of $\lambda_{j}$ can lead to different sparsity of $g$.
Let $\textbf{f}$ denote the column vector $\{f_{i}\}_{i=1}^{N}$, then the formulation (\ref{eq:1.1}) is equivalent to the following minimization:
\begin{equation}\label{eq:1.2}
\min_{\textbf{X}_{j},1\leq j\leq J} \|\sum_{j=1}^{J}\mathbf{A}_{j}\textbf{X}_{j}-\textbf{f}\|_{2}^{2}+\sum_{j=1}^{J}\lambda_{j}\|\textbf{X}_{j}\|_{1},
\end{equation}
where $\mathbf{A}_{j}$ is the observation matrix defined by
$$\mathbf{A}_{j}(i,k)=\varphi_{k}^{j}(x_{i}), \; i=1,2,...,N, \; k=1,2,...,n_{j},\; j=1,2,...,J.$$
For the case $\lambda_{1}=\cdots=\lambda_{J}$, the model reduces to the model presented in \cite{Shen}.
Moreover, the $l_{1}$ related minimization
resulting from our proposed model can be efficiently solved using the alternating direction method of multipliers (ADMM) \cite{Yuan2,Gabay}. 

This framework combines the ideas developed in compressed sensing with well-known
concepts arising in adaptive and multilevel finite element methods. The solution of the $l_{1}$ minimization problem
and the multilevel basis functions are used to control the grid refinement and adaptivity.
Since we aim at finding a sparse representation of the surface, we discard those
coefficients which are smaller than a certain threshold. Then only large coefficients of
the solution are left which indicate important contributions of the underlying surface. Moreover,
these large coefficients belong to the parts of the surface that have large fluctuations.
It seems a little similar for our method and the approach in \cite{DeVore}, both using the PSI space.
However, the method in \cite{DeVore} was used to approximate functions expressed as a infinite sum of wavelet decomposition by a finite
sum, while we deal with scattered data fitting problem by representing the fitting function as
a finite sum directly with certain accuracy and more sparse coefficients.
The behavior of our method is demonstrated via four examples: a discontinuous
function, a non-smooth function, a smooth function and the Franke test function. In addition,
we compare the numerical results with the AGLASSO model and the basic MBA algorithm in \cite{Lee} followed
by the same thresholding.

The rest of the paper is organized as follows. In the next section, we recall the main
ingredients of the PSI space which will be used here. Moreover, we will propose the sparsity
based regularization model for scattered data fitting. In Section 3, the ADMM
algorithm will be applied to solve the minimization problem resulted from the
proposed model. Numerical experiments are also performed to illustrate the algorithm.
Section 4 is the conclusion.

\section{Sparse solution of PSI approach to scattered data approximation}
For a given set of scattered points $\{x_{i}\}_{i=1}^{N}\subseteq \Omega \subseteq \mathbb{R}^{2}$
and the corresponding noisy data $\{f_{i}\}_{i=1}^{N}$, our task is to reconstruct a fitting surface
with a sparse representation.

\subsection{PSI space and $l_{1}$ regularization}
Let $\Omega \subseteq \mathbb{R}^{2}$ be a bounded domain of interest where all data lie in and
let $\varphi$ be a carefully chosen, compactly supported function (e.g. uniform B-spline,
box spline, radial basis functions). Denote
$$\Lambda=\{k\in \mathbb{Z}^{2}: supp(\varphi(\cdot/h-k))\cap\Omega\neq\emptyset\},\;S^{h}(\varphi, \Omega)=\{f|f=\mathop\sum\limits_{k\in \Lambda}c_{k}\varphi(\cdot/h-k)\},$$
where $h>0$ is a scaling parameter that controls the refinement of the space.
Denote $S_{j}=S^{h/2^{j-1}}(\varphi, \Omega), S_{1}\subseteq \cdots \subseteq S_{J}$,
we then look for a function $g\in S_{J}$ which fits closely the given
data. Then $g$ is composed of a sequence
of functions as $$g = g_{1} + g_{2} + \cdots + g_{J},$$ where $g_{i}\in S_{i},i =1, 2,...,J.$

Here we choose a proper PSI space generated by B-spline as the approximation space $S$ since it enjoys desirable
properties for data fitting. It has a simple structure and provides good approximation
to smooth functions, which leads to simple and accurate algorithms. Moreover, it can be
associated to a wavelet or frame system and hence one can solve the fitting problem by
making use of the advantages that a wavelet (frame) system can offer \cite{Shen}. These advantages include
sparse approximation of functions in the wavelet (frame) domain, multilevel structure of
basis  functions, adaptivity to the data, norm equivalence, etc.

Recall that a function $\varphi$ is said to satisfy the Strang-Fix conditions of order $m$ if
$$\hat{\varphi}(0)\neq 0, \; D^{\alpha}\hat{\varphi}(2\pi j)=0, \; \forall j\in \mathbb{Z}^{2}\setminus0,\; |\alpha|<m.$$
Denote $$W_{2}^{m}=\{f\in L_{2}(\mathbb{R}^{2}):\|f\|_{W_{2}^{m}}=\sqrt{2\pi}\|(1+|\cdot|^{m})\hat{f}\|_{L_{2}(\mathbb{R}^{2})}<+\infty\},$$
where $\hat{f}$ is the Fourier transform of the function $f$ and $|\cdot|=||\cdot||_{2}$ denotes the Euclidean norm.
Then if $\varphi$ satisfies the Strang-Fix conditions \cite{deBoor,Jia}, a PSI space provides good approximation to $W_{2}^{m}$ (see \cite{deBoor,Jia}),
i.e., $\varphi$ satisfies the Strang-Fix conditions of order $m$ if and only if for all $f\in W_{2}^{m}$,
$$\inf\limits_{s\in S^{h}(\varphi, \Omega)}\|f-s\|_{L^{2}(\mathbb{R}^{2})}=O(h^{m}).$$
Particularly, the B-spline $B_{m}$ of order $m$ satisfies the Strang-Fix conditions of
order $2$ for all $m\geq2$ \cite{Dong2}. For more detailed discussions on PSI space, see \cite{deBoor}.

Obviously, the union set of the basis functions of $S_{j}$ is not linearly independent. Thus the representation
of $g$ is not unique and we want to determine a relatively sparse one, i.e., a representation with
as many vanishing coefficients as possible. Every function
$g_{j}\in S_{j}$ can be written as $$g_{j} = \sum_{k\in I_{j}}X^{j}_{k}\varphi(2^{j-1}x/h-k),$$
where $$I_{j}=\{k\in \mathbb{Z}^{2},  supp(\varphi(2^{j-1}\cdot/h-k))\cap\Omega\neq\emptyset\}.$$
Let $\textbf{X}_{j}$ and $\textbf{f}$ denote the column vector $\{X^{j}_{k}\}_{k\in I_{j}}$ and
$\{f_{i}\}_{1\leq i\leq N}$ respectively, then the problem can be formulated as follows.
\begin{equation}\label{eq:2.1}
\min_{\textbf{X}_{j},1\leq j\leq J} \|\sum_{j=1}^{J}A_{j}\textbf{X}_{j}-\textbf{f}\|_{2}^{2}+\sum_{j=1}^{J}\lambda_{j}\|\textbf{X}_{j}\|_{1},
\end{equation}
where $A_{j}$ is the observation matrix defined by
$$A_{j}(i,k)=\varphi(2^{j-1}x_{i}/h-k), \; k\in I_{j}, \; i=1,2,...,N, \; j=1,2,...,J.$$

Obviously, the model (\ref{eq:2.1}) balances the fitting accuracy and the $l_{1}$ norm. In order to achieve
a sparse representation, small coefficients are neglected. That is, after obtaining the solution
$\{\textbf{X}_{j}\}_{j=1}^{J}$ of the model (\ref{eq:2.1}), we discard the small elements of $\textbf{X}_{j},1\leq j\leq J$.
Then the final solution only has large values left which indicate important contributions
(fluctuations) of the real surface. Furthermore, comparing with the multilevel approximation
approach given in \cite{Lee}, our method has the advantages of simplicity. Another important distinction is
that it can be interpretable as a sparse strategy for reconstructing scattered data.

\subsection{The MLASSO model}
The model (\ref{eq:2.1}) is related to the LASSO model in some extent.
Recall that the mathematical model of LASSO is:
$$\min_{\textbf{X}} \|A\textbf{X}-\textbf{f}\|_{2}^{2}+\mu\|\textbf{X}\|_{1}.$$
It was proposed originally in \cite{Tibshirani}, and plays a very influential role in variable selection
and dimensionality reduction. The Group LASSO (GLASSO) model proposed in \cite{Yuan1} solves the convex
optimization problem:
$$\min_{\textbf{X}_{j},1\leq j\leq J} \|\sum_{j=1}^{J}A_{j}\textbf{X}_{j}-\textbf{f}\|_{2}^{2}+\mu \sum_{j=1}^{J}\|\textbf{X}_{j}\|_{2},$$
where $\mu>0$ is a regularization parameter. The GLASSO model was proposed to perform
variable selection on groups of variables for linear regression models. It has many applications
in areas such as computer vision, data mining, etc. Meier et al. in \cite{Meier} extended the
GLASSO to logistic regression. The GLASSO does not, however, yield sparsity within
a group. Moreover, GLASSO suffers from estimation inefficiency and selection inconsistency.
To remedy these problems, the adaptive GLASSO method (AGLASSO) is proposed in \cite{Wang1} as:
$$\min_{\textbf{X}_{j},1\leq j\leq J} \|\sum_{j=1}^{J}A_{j}\textbf{X}_{j}-\textbf{f}\|_{2}^{2}+\sum_{j=1}^{J}\mu_{j}\|\textbf{X}_{j}\|_{2}.$$

Obviously, the model (\ref{eq:2.1}) reduces to the LASSO model when $\mu=\lambda_{1}=\cdots=\lambda_{J}$.
Moreover, the model (\ref{eq:2.1}) acts like the LASSO at the multilevel.
Therefore, we denote the model (\ref{eq:2.1}) as the multilevel LASSO model (MLASSO).
Compared with the AGLASSO model, the MLASSO model considers an additional penalty on the $l_{1}$ norm
instead of $l_{2}$ norm of the regression coefficient vector,
and it produces as election of variables with sparsity among different levels.
It is known that when $l_{2}$ norm regularization term is applied to the data set, the resulting
surface tends to be smooth without sharp discontinuities but have undesirable oscillations
near the discontinuities \cite{Shen}. Recently, several surface reconstruction approaches have been proposed
to preserve surface discontinuity by replacing $l_{2}$ regularization using more sophisticated
regularization, e.g., the Huber approximation of $l_{2}$ norm of function derivatives in \cite{Stevenson},
the local kernel regularization in \cite{Gijbels} and the non-local means regularization in \cite{Dong2}.
In our method, we use the $l_{1}$ regularization instead, in order to obtain the relatively sparse solution.
The regularization term $\|\cdot\|_{1}$ in (\ref{eq:2.1}) penalizes the roughness of the solution.

\subsection{The parameter selection}
The determination of the proper value of $\{\lambda_{j}\}$ in the MLASSO model is an important
problem and depends on the variance of the noise $\textbf{n}$, the properties of $\{A_{j}\}$ and $\|\cdot\|_{1}$.
An appropriate choice of the regularization parameters is of vital importance for the
quality of the resulting estimate and has been the subject of extensive research \cite{Belge}.

In recent years, there has been a growing interest in sophisticated regularization techniques
which use multiple constraints as a mean of improving the quality of the solution \cite{Belge}.
Among them, only a few papers discussed the choice of multiple regularization parameters.
However, most of them discuss the case that the regularization term is the $l_{2}$ norm.
For example, a multi-parameter generalization of heuristic L-curve has been proposed in \cite{Belge},
a knowledge of noise (covariance) structure is required for a choice of parameter in \cite{Bauer1,Bauer2,Chen},
some reduction to a single parameter choice is suggested in \cite{Brezinski}.
At the same time, the discrepancy principle, which is widely used and known
as the first parameter choice strategy proposed in the regularization theory \cite{Phillips}, has
been discussed in a multi-parameter context in \cite{Lu}. Of course, there might be many different
methods to choose the regularization parameters satisfying certain principle.
In fact, the choice of the regularization parameters in the regularization modes, such as the LASSO model
and the standard Tikhonov regularization model, have not yet solved. No choice is available for all the
models. Basically, different models need different methods to decide the parameters, and even the same one
may have different methods to choose, such as \cite{Zou,Wang1,Hanke,Wang4,Hao}.

For the MLASSO model, instead of discussing similar methods as listed above, we propose two
simple criteria for the choice of the parameters $\{\lambda_{j}\}_{j=1}^{J}$
according to the sparsity and the support of the basis functions of the different levels.
On one hand, since the length of $\mathbf{X}_{j}$ is less than that of $\mathbf{X}_{j+1}$, so if one
wants to obtain sparser solution, the parameters of the last several levels should be larger. In particular,
choosing the same value for all the parameters, i.e., $\lambda_{1}=\cdots=\lambda_{J}$,
obtains the global sparsity in the whole space $S$. On the other hand, the support of the basis functions of the $j$-th level
is larger than that of the  $(j+1)$-th level, so if one wants smaller support, the parameters of the
first several levels should be larger. Hence, the more appropriate choice of the parameters
is to make the parameters of the first and last several levels greater than the middle several levels.
We can only give such a qualitative guideline, since quantitative guidance for regularization parameters is rarely available.
In this way, the solution has smaller support and sparser representation with good approximation accuracy.
Numerical results also confirm this selection method as shown in the Tables 1-8.
Moreover, in contrast with the methods listed above, this choice is easier and cheaper in the sense of computational cost.

Based on the above, the proposed method offers some interesting advantages: \\
1) Related with the LASSO model, the MLASSO model has a strong statistical background;\\
2) Compared to the LASSO model, the parameters in the MLASSO model can be chosen differently, thus one can attain high flexibility for
   the approximation accuracy and the sparsity of the solution;  \\
3) Compared to the AGLASSO model, the MLASSO model can provide the sparsity of the solution by choosing different parameters and
   the distribution of the solution can reflect the large fluctuations of the underlying surface;\\
4) Utilizing the ADMM algorithm, the MLASSO model can be efficiently solved.   \\

\section{Numerical algorithm}
\subsection{Algorithm}
In this section, we use the ADMM algorithm to solve the minimization model (\ref{eq:2.1}) for experimental evaluation.
It turns out that ADMM is equivalent to or closely related ro many famous algorithms, such as the Douglas-Rachford splitting method
in PDE literature \cite{Douglas,Eckstein,Lions}, the Bregman iterative algorithms for $l_{1}$ problems in signal processing \cite{Goldstein}
and many others. In particular, we refer to \cite{Esser,Wu-Tai,Steidl,Setzer} for the relationship between ADMM and the split Bregman iteration scheme which is
very influential in the area of image processing. Convergence analysis of the ADMM was given in \cite{Eckstein2,He,Boyd}. 

In order to solve the MLASSO model (\ref{eq:2.1}) and guarantee the convergence, we denote
$$A=(A_{1},A_{2},\ldots,A_{J}),\quad \textbf{X}=(\textbf{X}^{T}_{1},\textbf{X}^{T}_{2},\ldots,\textbf{X}^{T}_{J})^{T},$$
$\lambda=diag(\lambda_{1},\lambda_{2},\ldots,\lambda_{J})$ is a diagonal matrix with $\lambda_{j},j=1,...,J$ as the main diagonal.
Then we can rewrite (\ref{eq:2.1}) as
\begin{equation}\label{eq:3.1}
  \begin{array}{l}
    \mathop{}\min\limits_{\textbf{X}} \|A\textbf{X}-\textbf{f}\|_{2}^{2}+\|\lambda\textbf{X}\|_{1}. \\
  \end{array}
\end{equation}
By introducing an auxiliary variable $\textbf{d}=\lambda\textbf{X}$, we convert the unconstrained minimization problem (\ref{eq:3.1}) into a constrained one:
\begin{equation} \label{eq:3.2}
  \begin{array}{l}
    \mathop{}\min\limits_{\textbf{X},\textbf{d}=\lambda\textbf{X}} \|A\textbf{X}-\textbf{f}\|_{2}^{2}+\|\textbf{d}\|_{1}. \\
  \end{array}
\end{equation}
In this way, the MLASSO model (\ref{eq:2.1}) is turned into a classical $l_{1}$ minimization problem. The
augmented Lagrangian function of problem (\ref{eq:3.2}) is
$$L(\textbf{X},\textbf{d},\textbf{b})=\|A\textbf{X}-\textbf{f}\|_{2}^{2}+\|\textbf{d}\|_{1}+\frac{\beta}{2}\|\lambda\textbf{X}-\textbf{d}\|_{2}^{2}+<\textbf{b},\lambda\textbf{X}-\textbf{d}>,$$
where $\beta>0$ is a parameter of the algorithm. Then by applying the ADMM method, given the initialization $\{\textbf{X}^{0},\textbf{d}^{0},\textbf{b}^{0}\}$ and the parameters
$\{\lambda_{j}\}_{j=1}^{J}$, it results in the following optimization algorithm:
\begin{equation}\label{eq:3.3}
\left\{
  \begin{array}{l}
    \textbf{X}^{k+1} =\arg\min\limits_{\textbf{X}}\|A\textbf{X}-\textbf{f}\|_{2}^{2}
    +\frac{\beta}{2}\|\lambda\textbf{X}-\textbf{d}^{k}+\frac{\textbf{b}^{k}}{\beta}\|_{2}^{2},\\
    \textbf{d}^{k+1} = \arg\min\limits_{\textbf{d}}\|\textbf{d}\|_{1}
    +\frac{\beta}{2}\|\lambda\textbf{X}^{k+1}-\textbf{d}+\frac{\textbf{b}^{k}}{\beta}\|_{2}^{2},\\
    \textbf{b}^{k+1} = \textbf{b}^{k}+(\lambda\textbf{X}^{k+1}-\textbf{d}^{k+1}).\\
  \end{array}
\right.
\end{equation}
First of all, the above algorithm is convergent, since it is just the classical ADMM for two block of variables. Secondly, the system can be simply solved. The first step of each iteration in (\ref{eq:3.3}) is
$$(2A^{T}A+\beta \lambda^{T}\lambda)\textbf{X}=2A^{T}f+\beta\lambda^{T} (\textbf{d}^{k}-\frac{\textbf{b}^{k}}{\beta}).$$
This linear system is positive definite and therefore it can be solved by the conjugate gradient method (CG).
When the matrix is ill-posed, i.e., its condition number is huge, the convergence rate of the CG will be very slow. Under this case,
the preconditioned CG \cite{Massimiliano,Sajo-Castelli,Benzi} can be used instead to reduce the condition number of the coefficient matrix
and improve the convergence speed.
The second subproblem has a simple analytical solution based on soft-thresholding operator \cite{Donoho}, that is
$$\textbf{d}^{k+1} = T_{\frac{1}{\beta}}(\lambda\textbf{X}^{k+1}+\frac{\textbf{b}^{k}}{\beta}),\\$$
where $T_{\theta}$ is the soft-thresholding operator defined by
$$T_{\theta}:x=[x_{1},x_{2},\ldots,x_{M}]\rightarrow T_{\theta}(x)=[t_{\theta}(x_{1}),t_{\theta}(x_{2}),\ldots,t_{\theta}(x_{M})],$$
where
$$t_{\theta}(\xi)=sgn(\xi)max\{0,|\xi|-\theta\}.$$
The complete description of the algorithm for solving the model (\ref{eq:2.1}) is provided as Algorithm 1 as follows:
\begin{algorithm}(Adapted ADMM for solving the MLASSO model (\ref{eq:2.1}))\\
Step 1) Set $J$ and the initial values $\{\textbf{X}^{0},\textbf{d}^{0},\textbf{b}^{0}\},$ choose appropriate sets of parameters
$\{\lambda_{j}\}_{j=1}^{J}$, $\beta$ and two thresholds $\sigma,\epsilon$; \\
Step 2) For $k=0,1,\ldots$, perform the iteration (\ref{eq:3.3}) until convergence;  \\
Step 3) Assume $\tilde{\textbf{X}}$ is the solution obtained from Step 2), if $|\tilde{\textbf{X}}(i)|\leq\sigma$, i.e., the absolute value of the $i$-th element of $\tilde{\textbf{X}}$ is less than $\sigma$, set $\tilde{\textbf{X}}(i)=0$;  \\
Step 4) The final solution $\textbf{X}$ are $\tilde{\textbf{X}}$ after the treatment of Step 3).
\end{algorithm}
In our numerical experiments, the initializations are $\textbf{X}^{0}=\textbf{d}^{0}=\textbf{b}^{0}=0, \beta=1$ and
the stopping criteria is $$\|\textbf{d}^{k}-\lambda\textbf{X}^{k}\|_{2}\leq\epsilon.$$

\begin{figure}[htbp]
\renewcommand{\figurename}{Fig.}
 \centering
  \subfigure[900 randomly scattered data points for $f_{1}$.]{
  \includegraphics[width=0.2\textwidth]{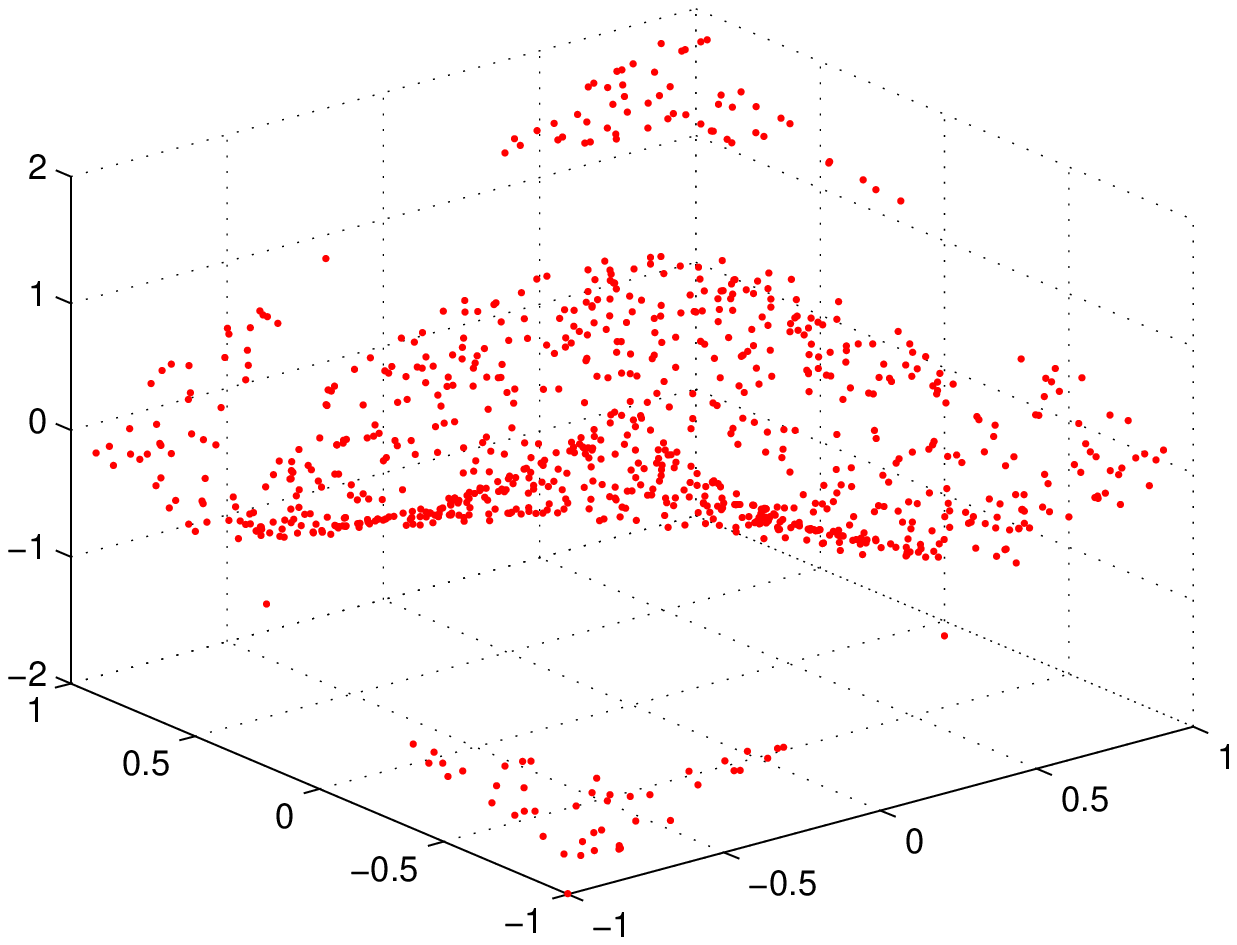}}
  \subfigure[The method in \cite{Lee}.]{
  \includegraphics[width=0.2\textwidth]{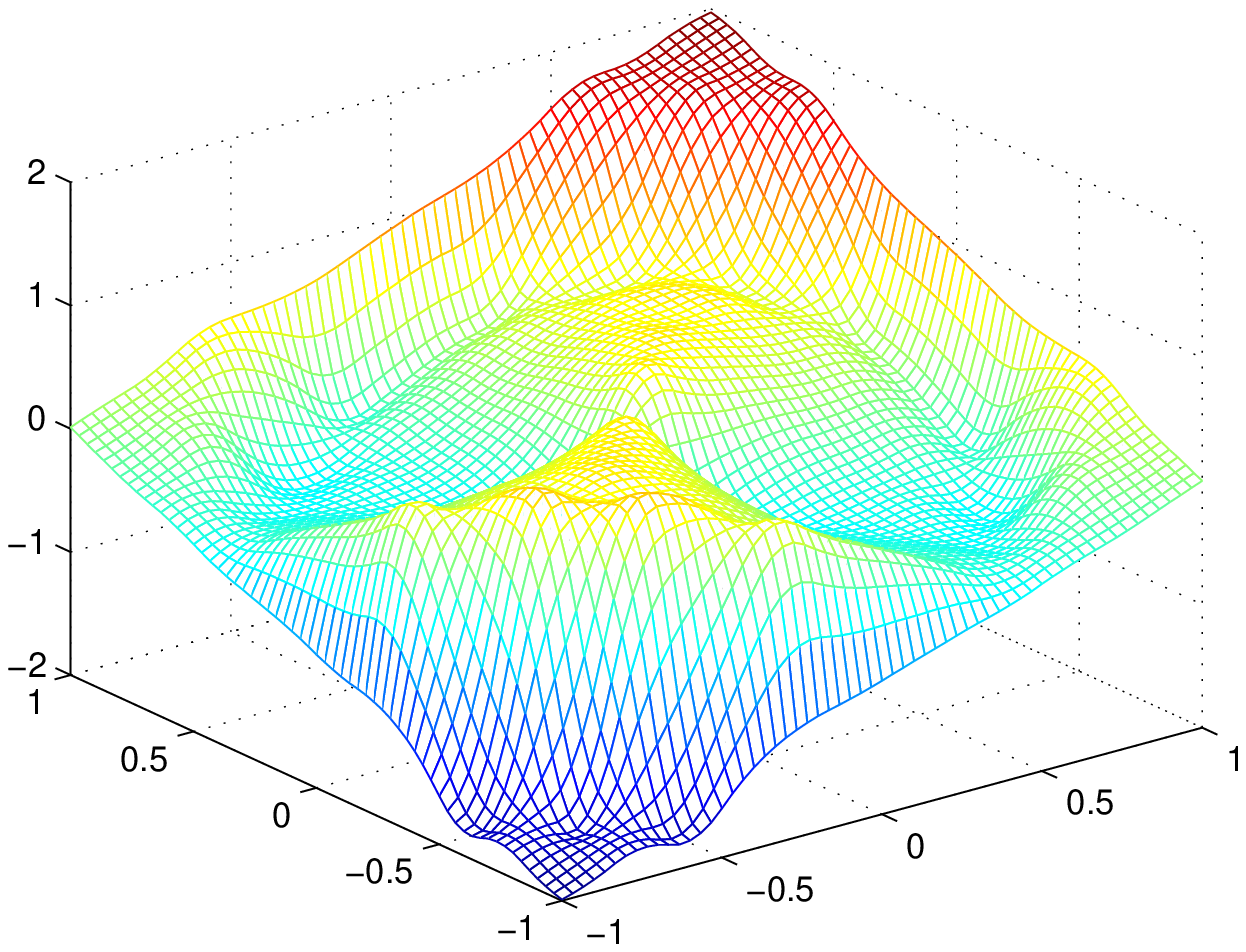}}
  \subfigure[$J=4$ for MLASSO: $\lambda_{1}=\lambda_{2}=\lambda_{3}=\lambda_{4}=0.01$.]{
  \includegraphics[width=0.2\textwidth]{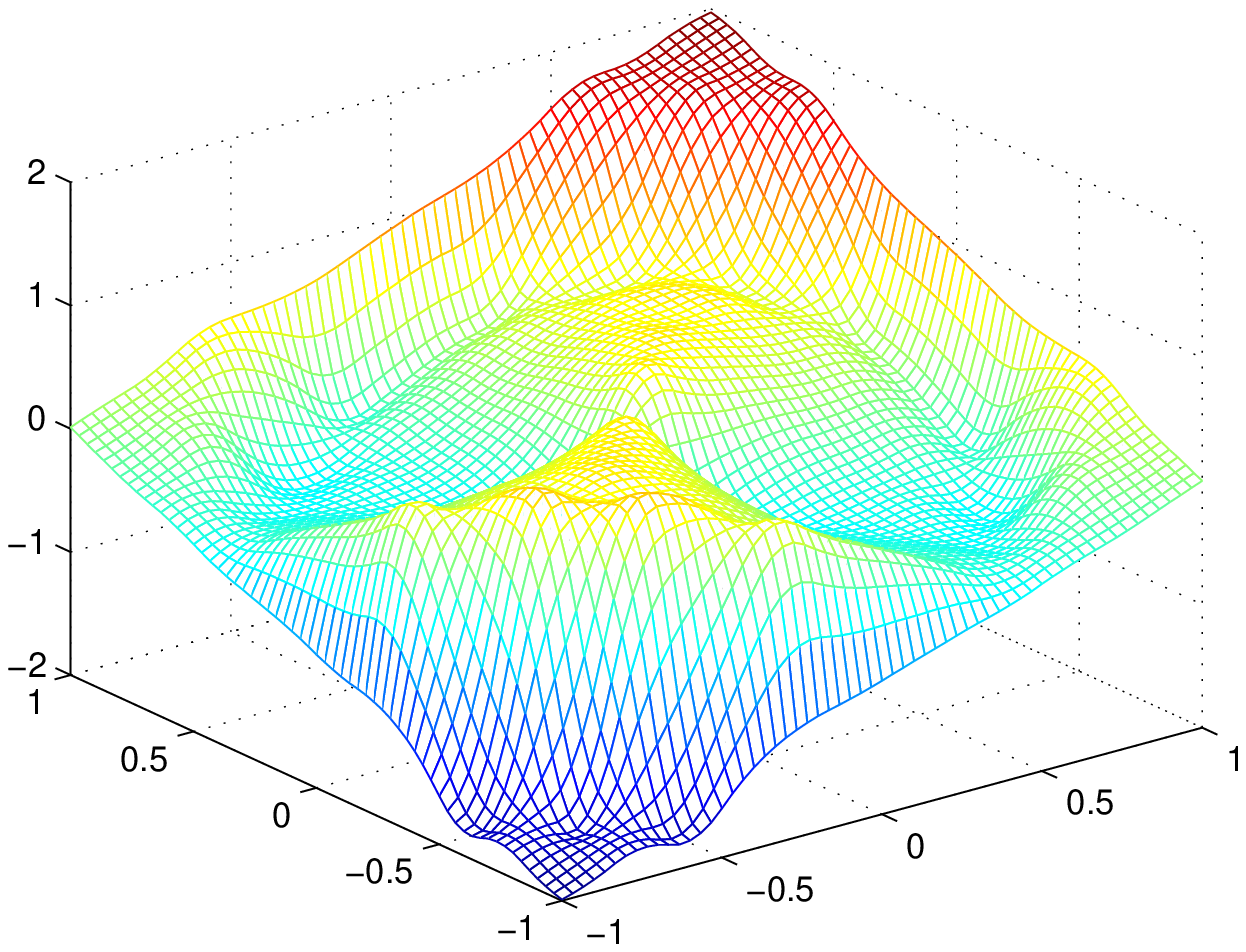}}
  \subfigure[$J=4$ for AGLASSO: $\lambda_{1}=\lambda_{2}=\lambda_{3}=\lambda_{4}=0.01$.]{
  \includegraphics[width=0.2\textwidth]{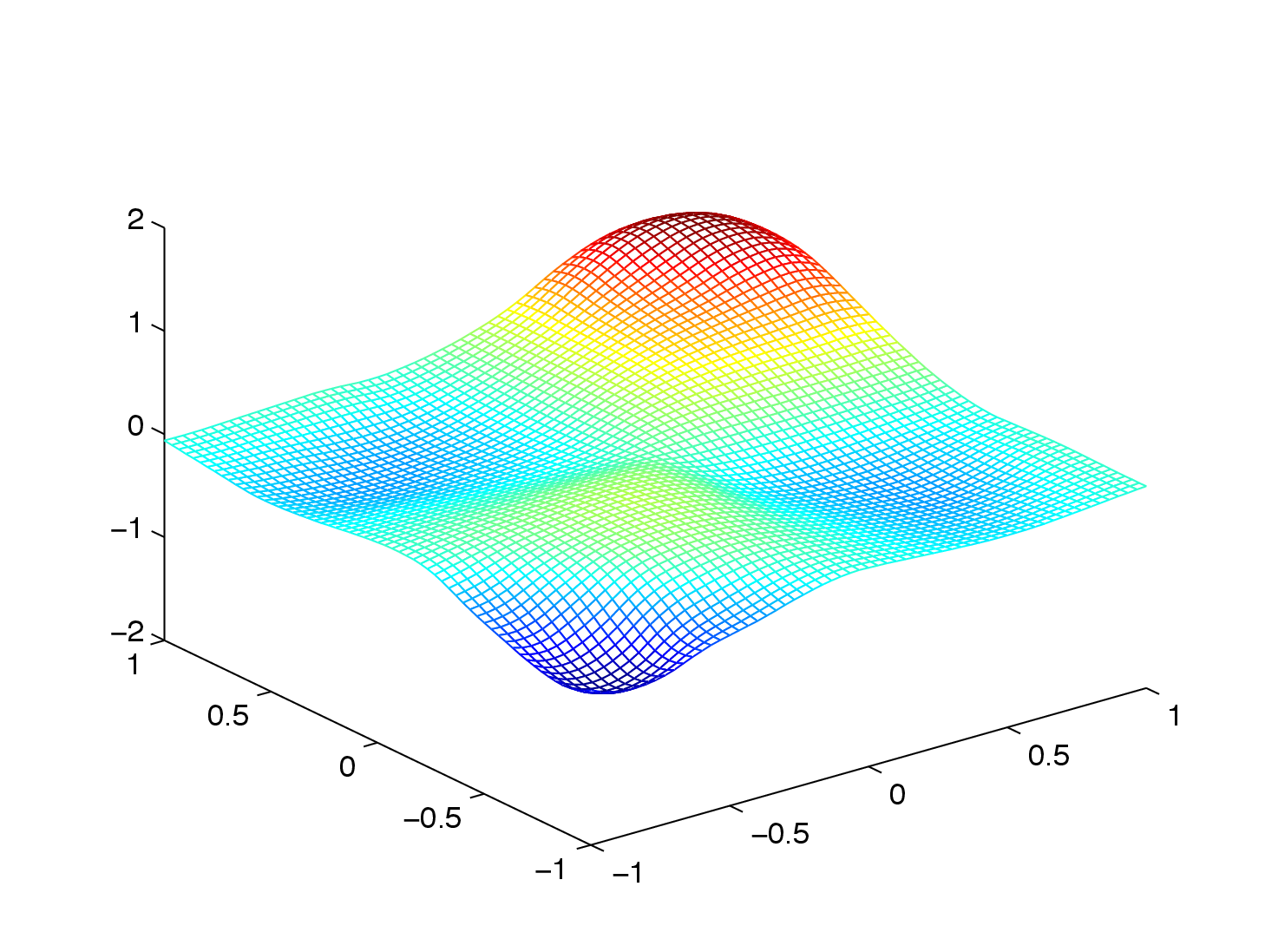}}
  \subfigure[900 randomly scattered data points for $f_{2}$.]{
  \includegraphics[width=0.2\textwidth]{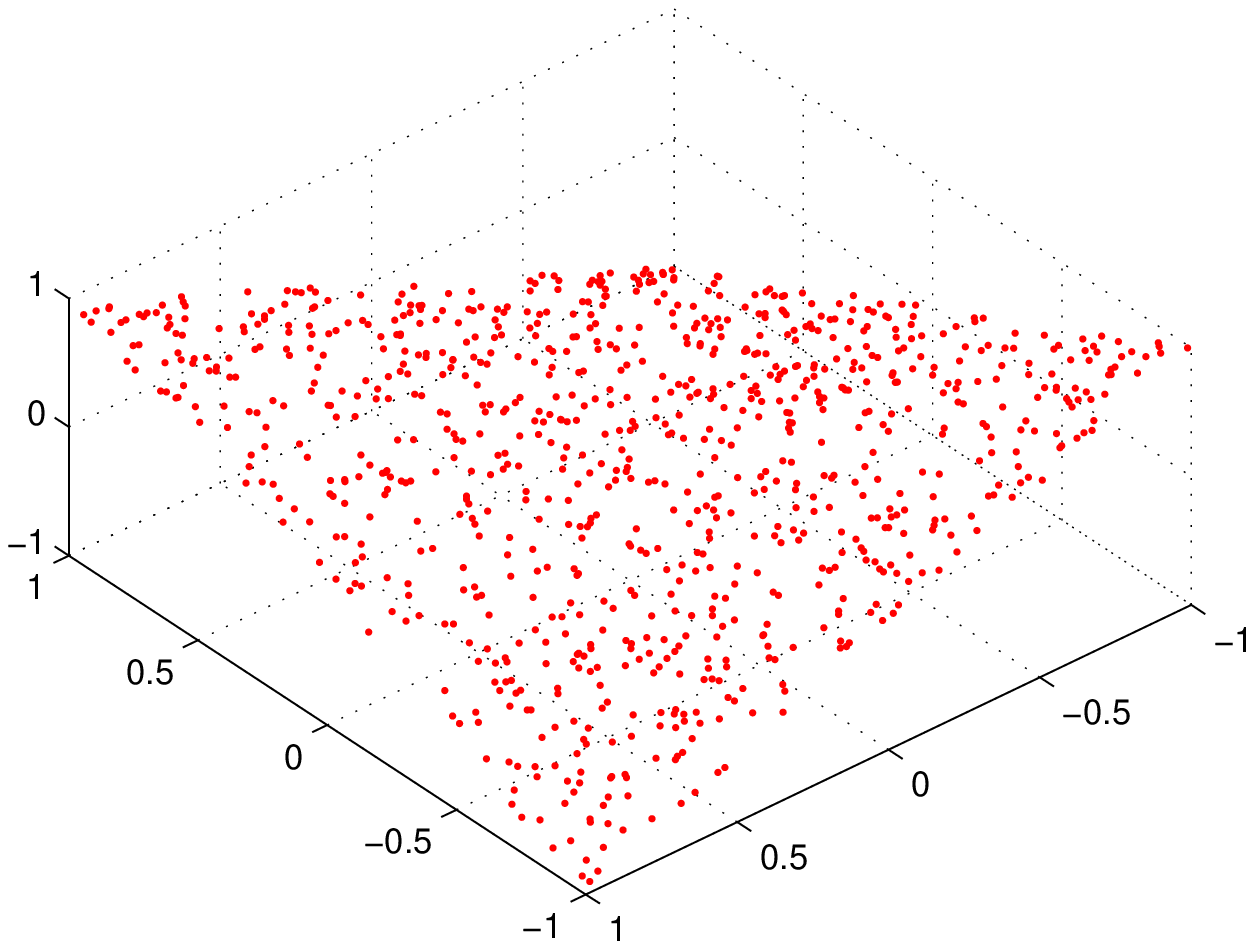}}  
  \subfigure[The method in \cite{Lee}.]{
  \includegraphics[width=0.2\textwidth]{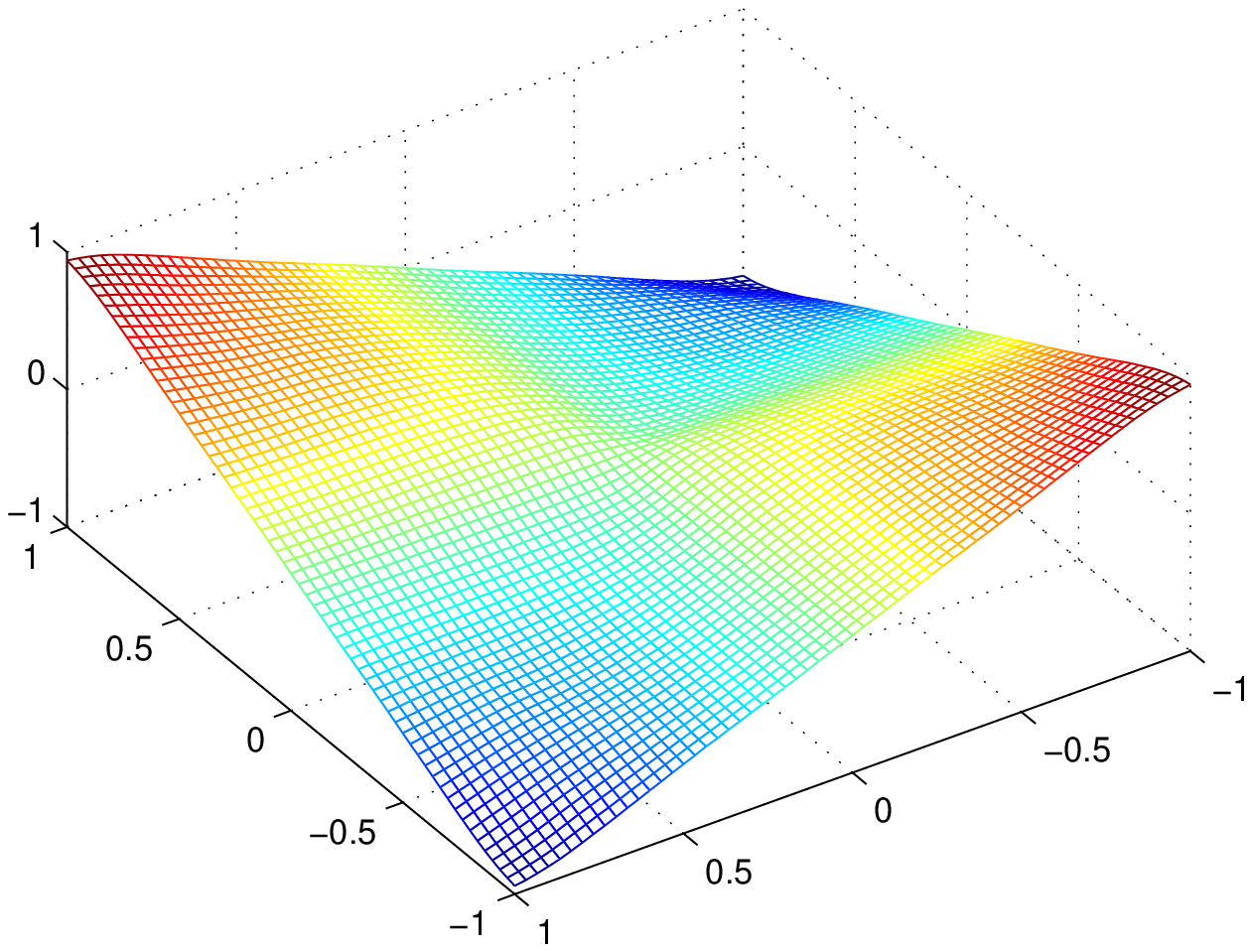}}
  \subfigure[$J=4$ for MLASSO: $\lambda_{1}=\lambda_{2}=\lambda_{3}=\lambda_{4}=0.01$.]{
  \includegraphics[width=0.2\textwidth]{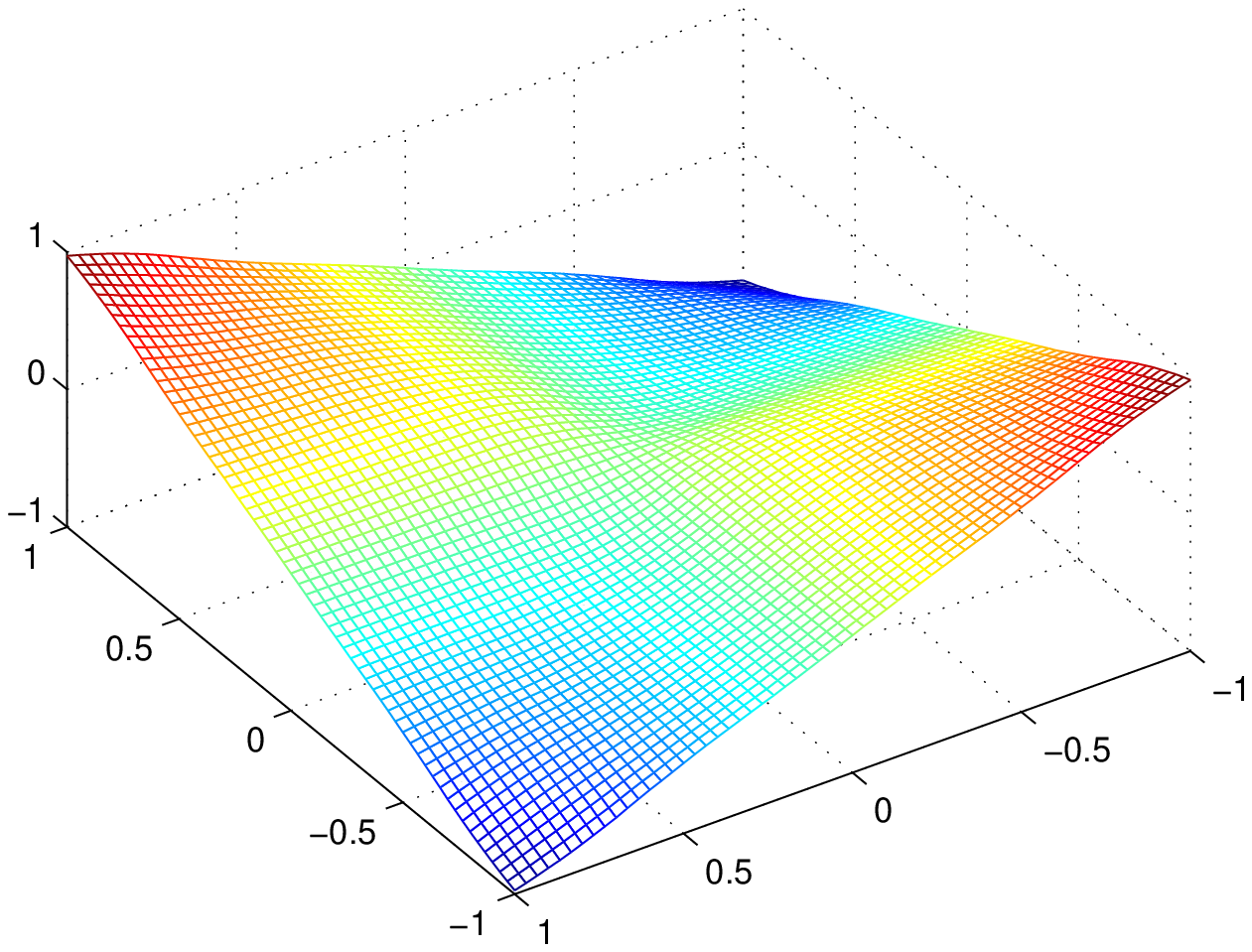}}
  \subfigure[$J=4$ for AGLASSO: $\lambda_{1}=\lambda_{2}=\lambda_{3}=\lambda_{4}=0.01$.]{
  \includegraphics[width=0.2\textwidth]{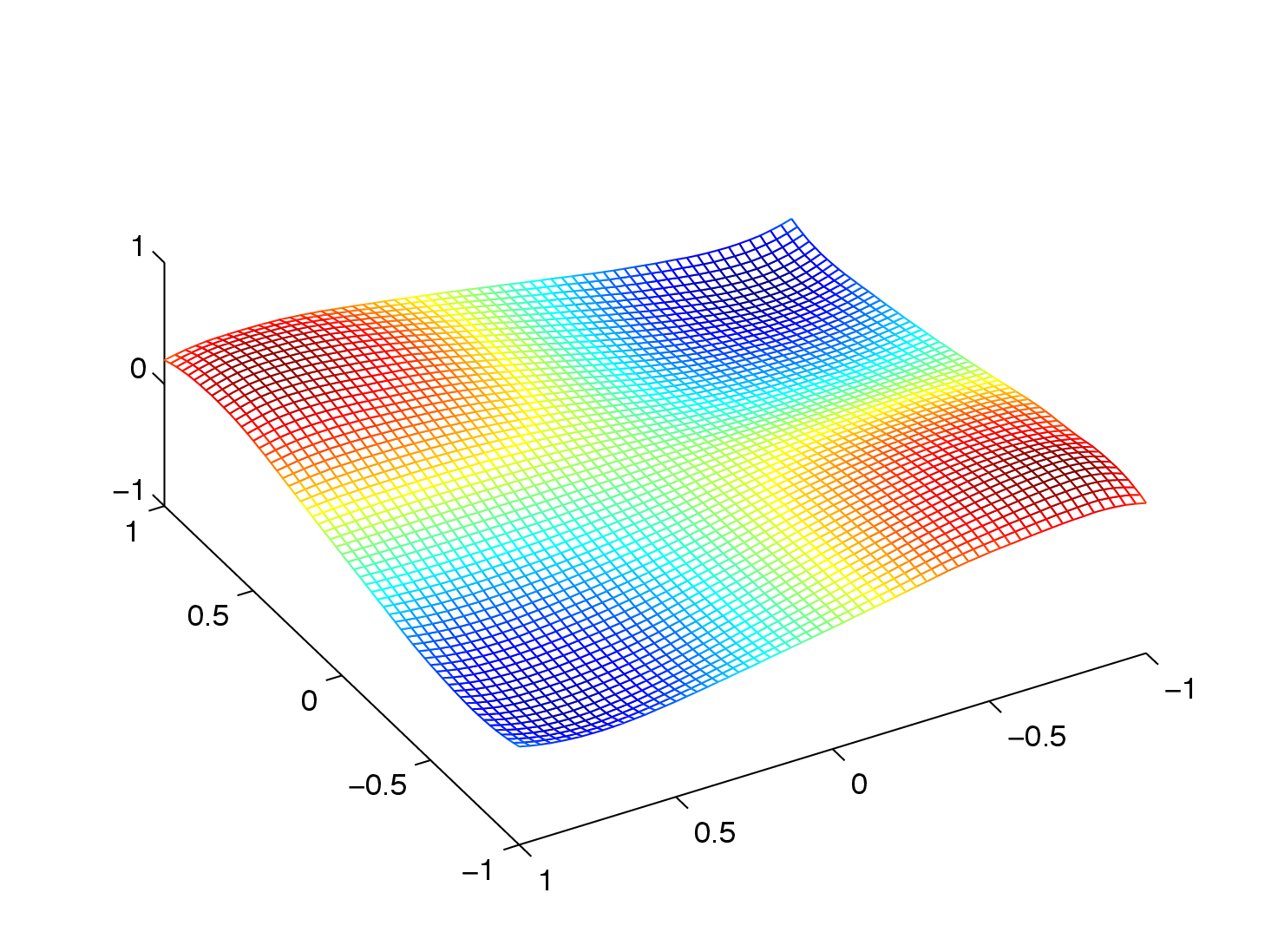}}
  \subfigure[900 randomly scattered data points for $f_{3}$.]{
  \includegraphics[width=0.2\textwidth]{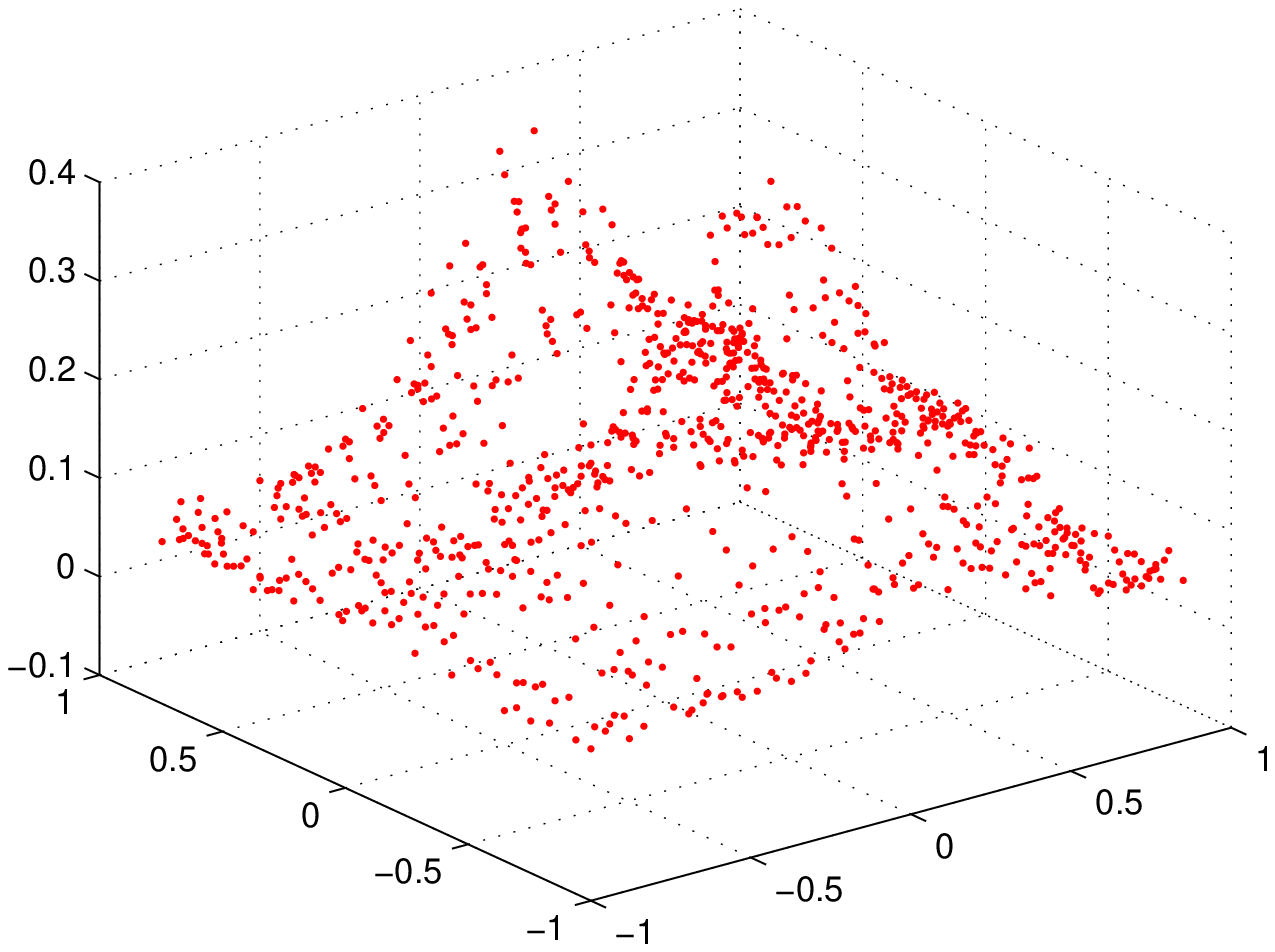}}  
  \subfigure[The method in \cite{Lee}.]{
  \includegraphics[width=0.2\textwidth]{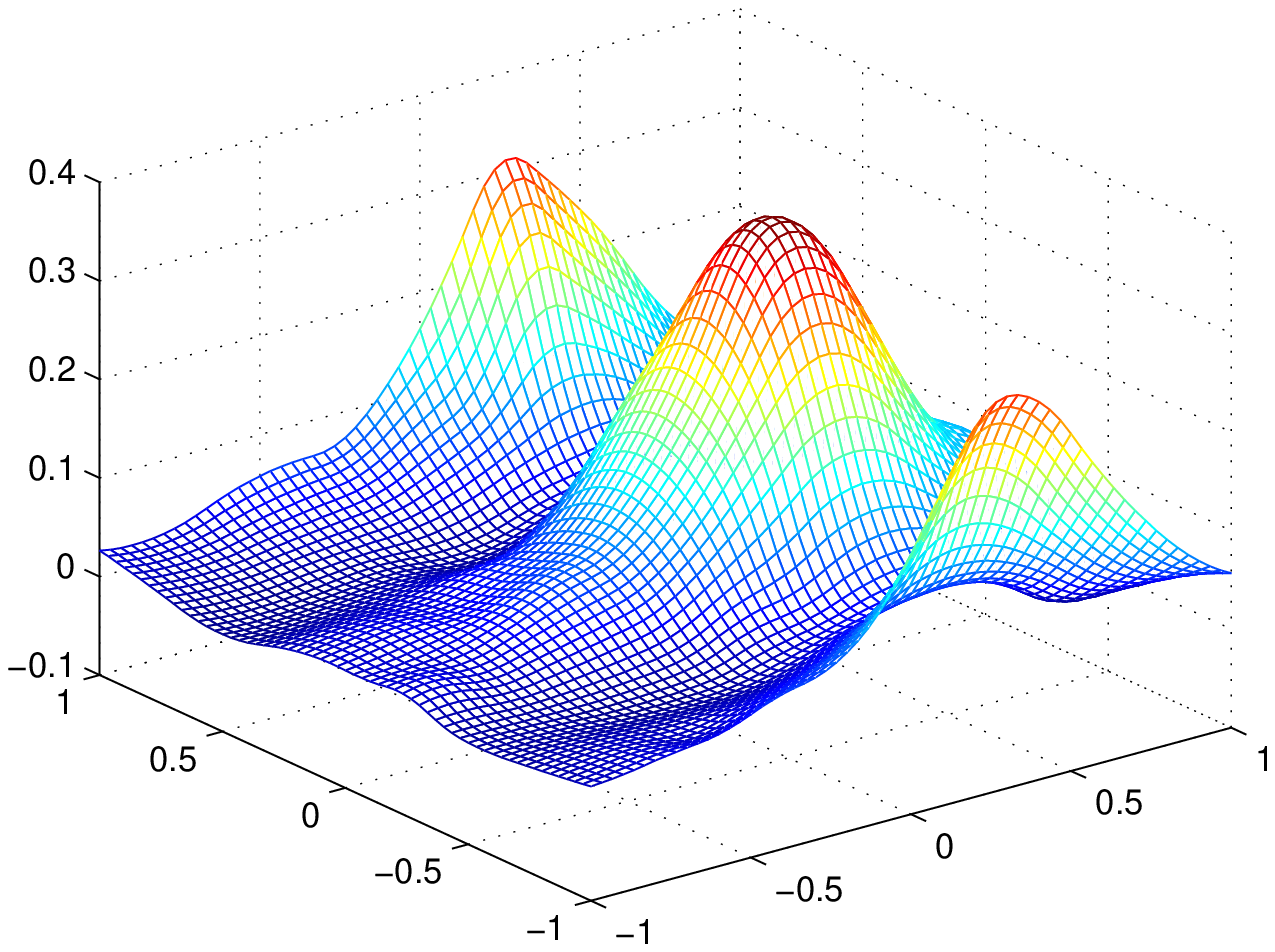}}
  \subfigure[$J=4$ for MLASSO: $\lambda_{1}=\lambda_{2}=\lambda_{3}=\lambda_{4}=0.01$.]{
  \includegraphics[width=0.2\textwidth]{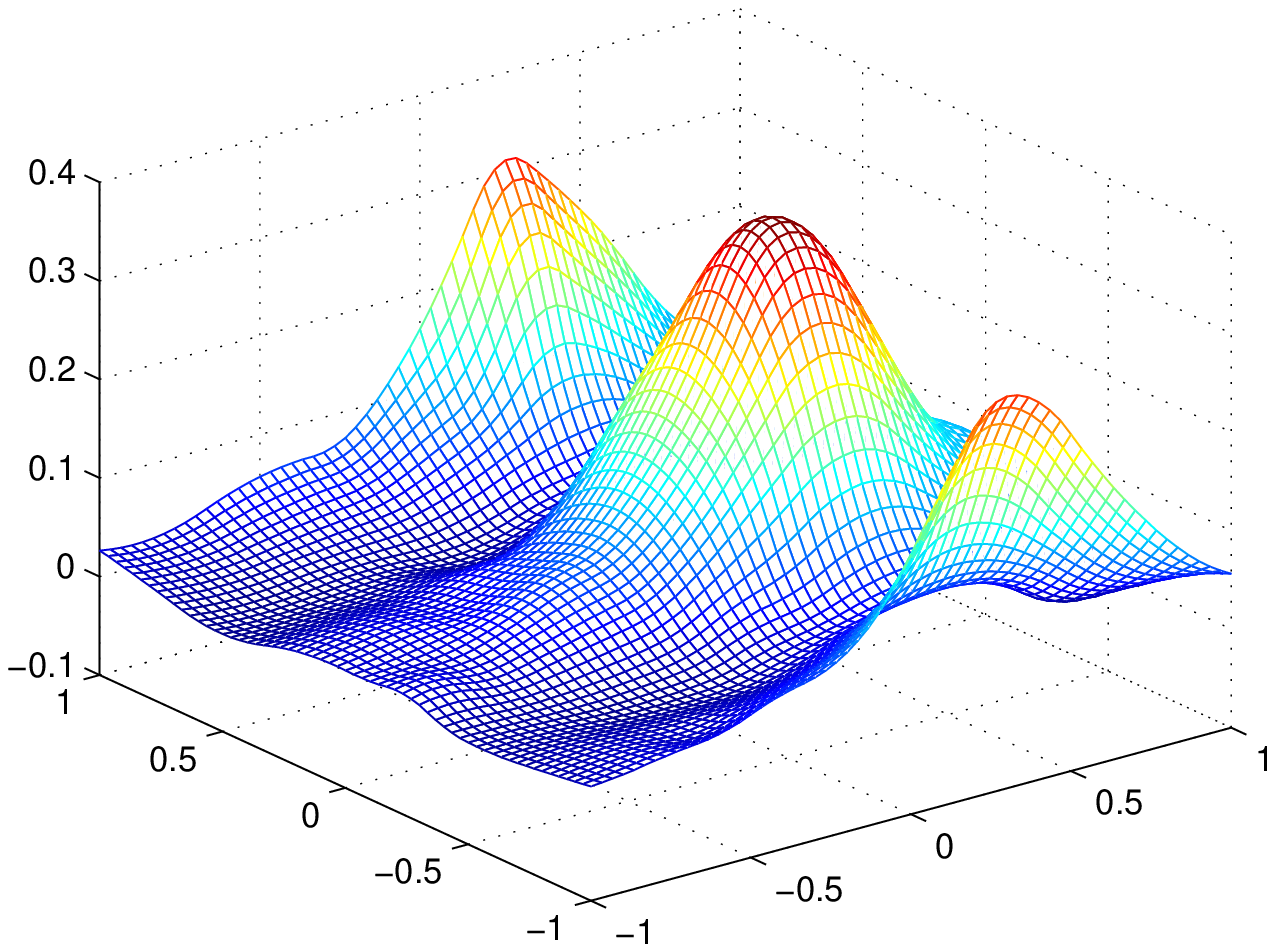}}
  \subfigure[$J=4$ for AGLASSO: $\lambda_{1}=\lambda_{2}=\lambda_{3}=\lambda_{4}=0.01$.]{
  \includegraphics[width=0.2\textwidth]{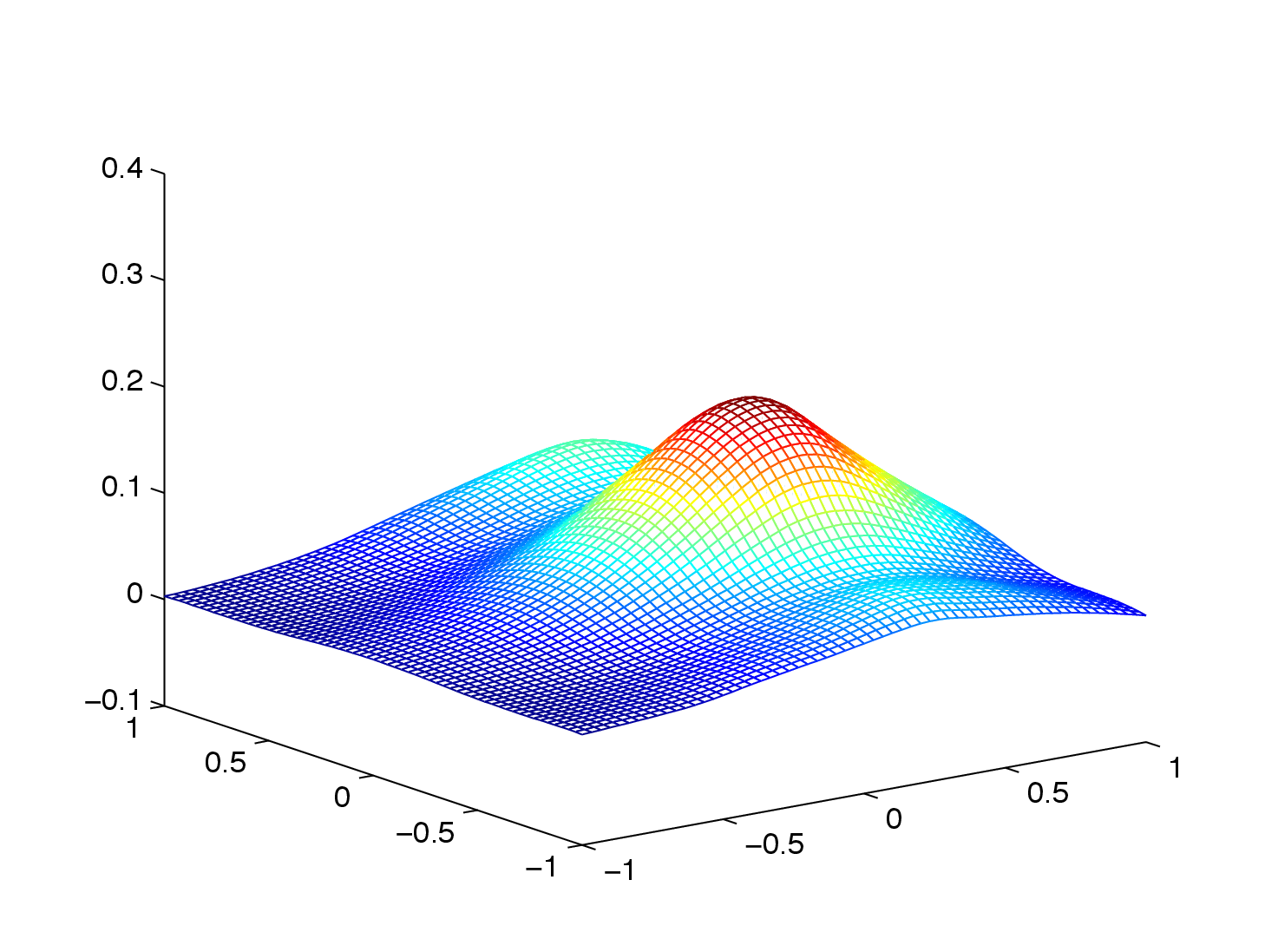}}
  \subfigure[900 randomly scattered data points for $f_{4}$.]{
  \includegraphics[width=0.2\textwidth]{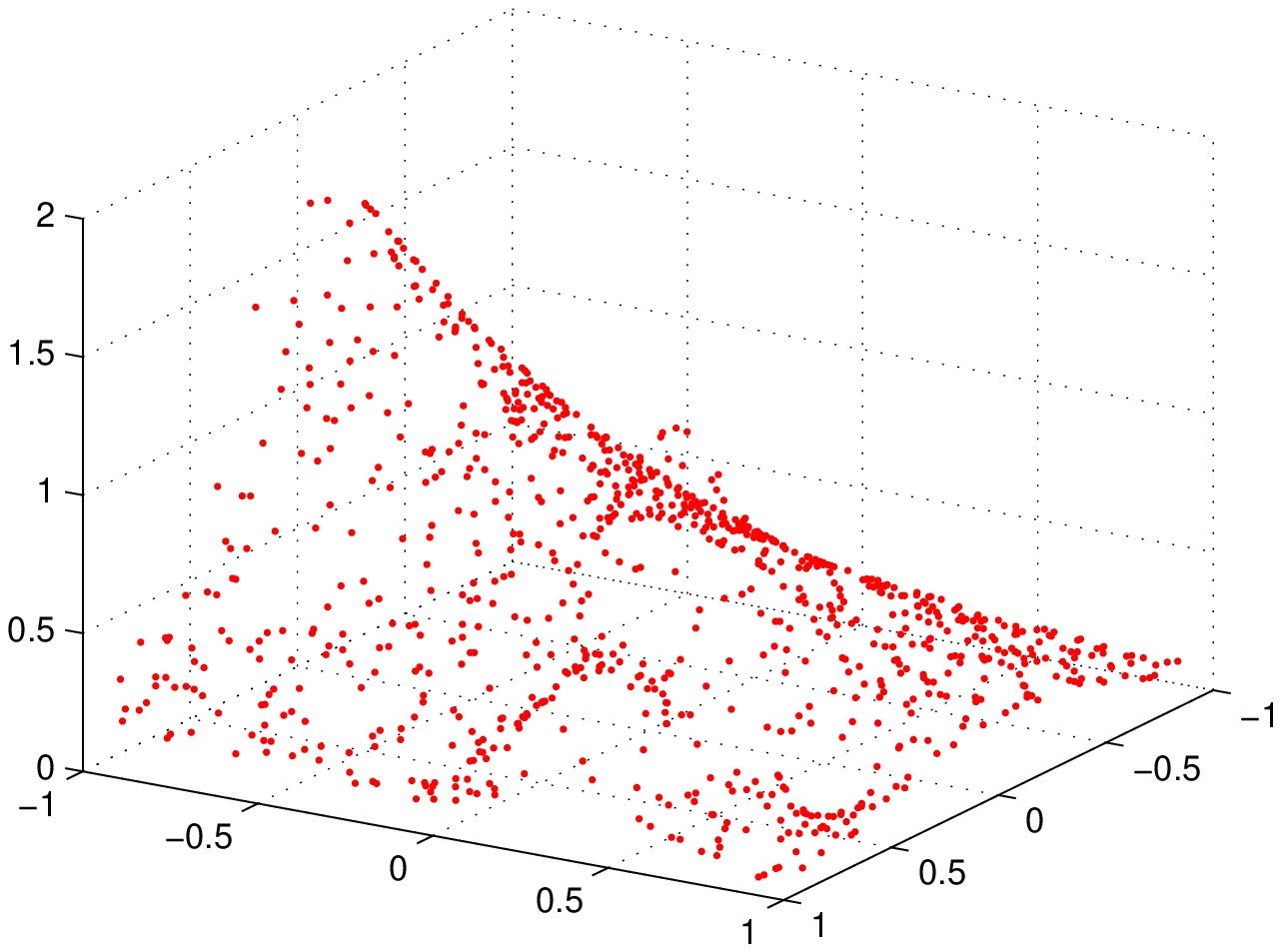}}
  \subfigure[The method in \cite{Lee}.]{
  \includegraphics[width=0.2\textwidth]{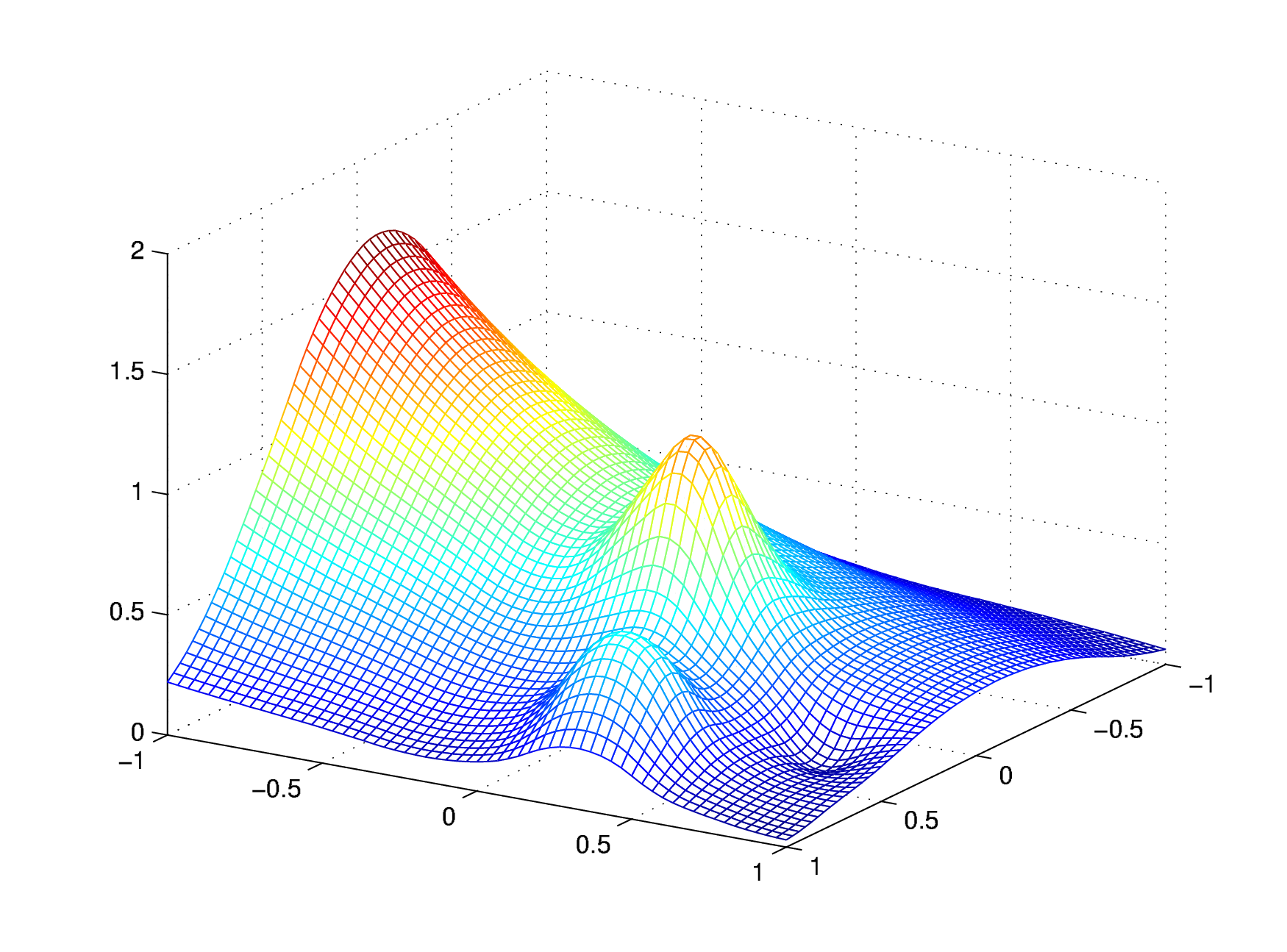}}
  \subfigure[$J=4$ for MLASSO: $\lambda_{1}=0.03,\lambda_{2}=0.02,\lambda_{3}=0.02,\lambda_{4}=0.04$.]{
  \includegraphics[width=0.2\textwidth]{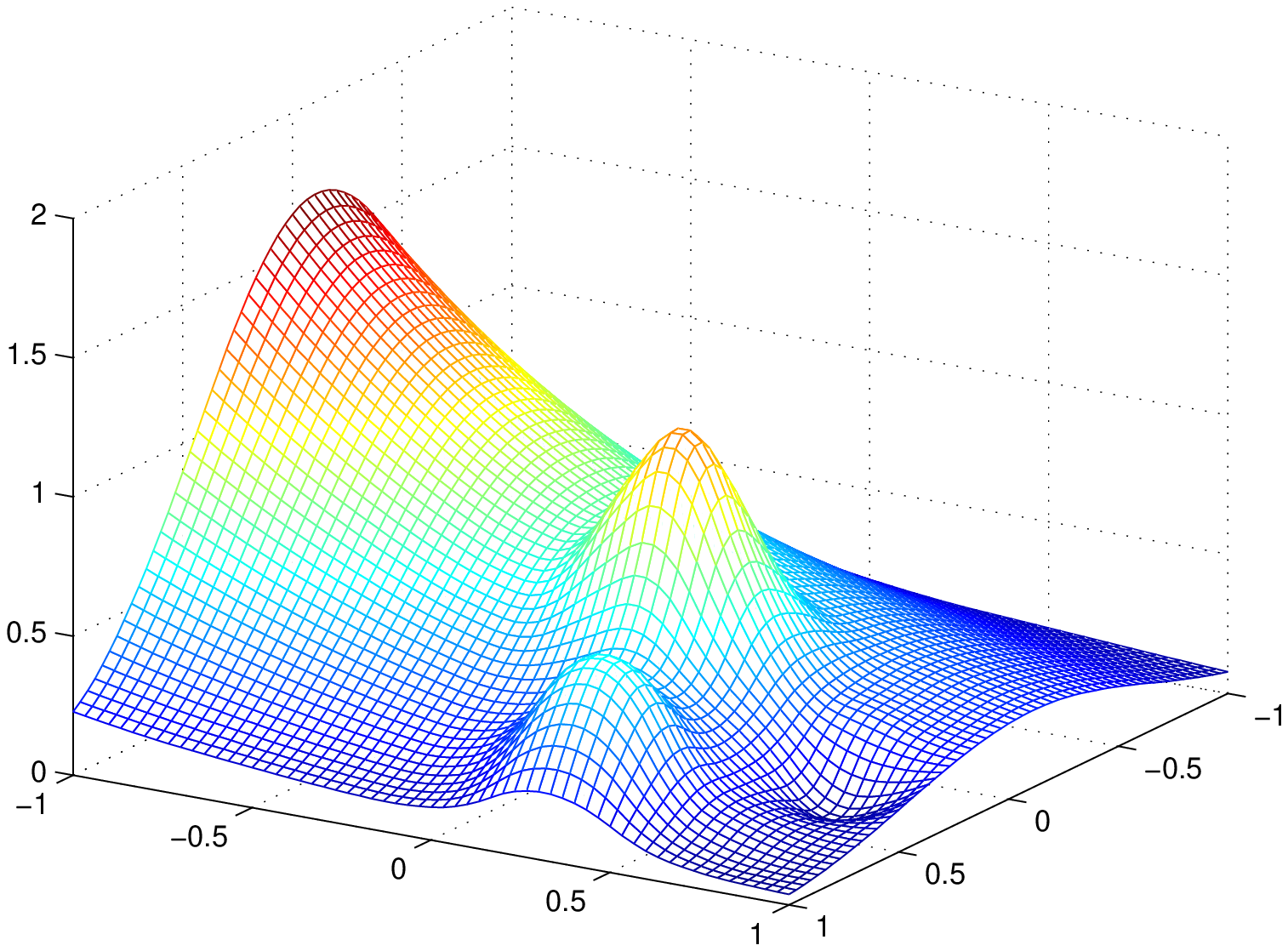}}
  \subfigure[$J=4$ for AGLASSO: $\lambda_{1}=0.03,\lambda_{2}=0.02,\lambda_{3}=0.02,\lambda_{4}=0.04$.]{
  \includegraphics[width=0.2\textwidth]{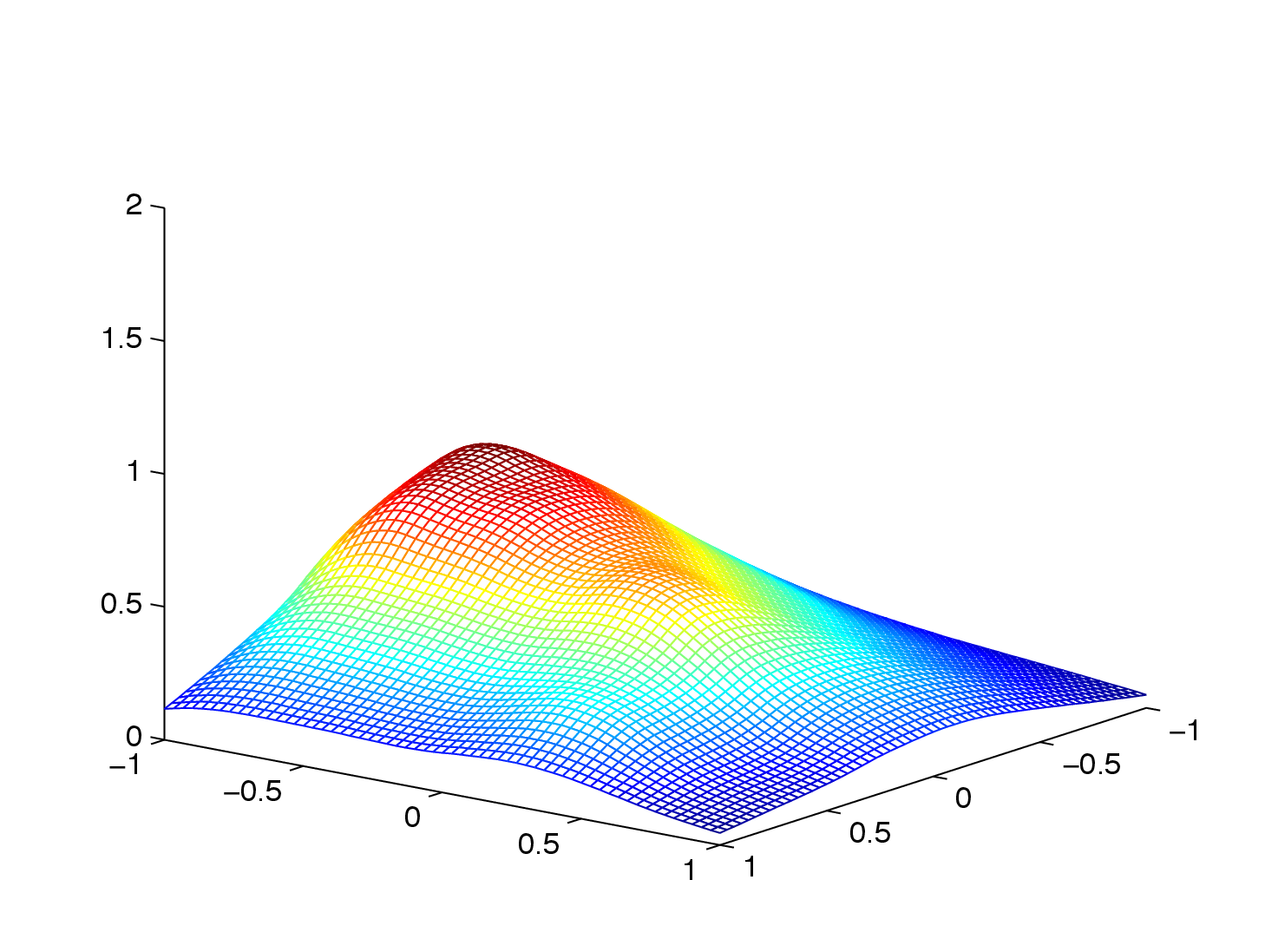}}
  \caption{The scattered data points and the corresponding approximation surfaces: (a)-(d) for $f_{1}$, (e)-(h) for $f_{2}$, (i)-(l) for $f_{3}$, (m)-(p) for $f_{4}$.}
  \label{11}
\end{figure}

\subsection{Numerical experiments}
Given a test function $f(x,y)$ with $\Omega =[-1,1]\times [-1,1]$, we first sample data points with certain noises from it, i.e.,
$\{(x_{i},y_{i},f(x_{i},y_{i})+\varepsilon_{i})\}$. The error vector
$\textbf{n}=(\varepsilon_{i})_{i}$, whose entries consist of the pseudorandom values drawn
from the standard uniform distribution on the open interval $(-\max\limits_{i}\frac{|f(x_{i},y_{i})|}{10},\max\limits_{i}\frac{|f(x_{i},y_{i})|}{10})$.
Then apply different methods to obtain the approximation function $g$. The difference between $f$ and $g$ is
measured by computing the normalized RMS (root mean square) error \cite{Lee}. That is,
$$RMS=\sqrt{\frac{\sum_{i,j=1}^{M_{1},N_{1}}(g(\tilde{x}_{i},\tilde{y}_{j})-f(\tilde{x}_{i},\tilde{y}_{j}))^{2}}{M_{1}N_{1}}},$$
where $\tilde{x}_{i}=-1+\frac{2i}{M_{1}-1}, i=0,1,\ldots,M_{1}-1, \,\tilde{y}_{j}=-1+\frac{2j}{N_{1}-1}, j=0,1,\ldots,N_{1}-1$, $M_{1}=N_{1}=50$.
Moreover, denote
$$Error=\sqrt{\frac{\sum_{i=1}^{N}(g(x_{i},y_{i})-f(x_{i},y_{i}))^{2}}{N}}$$
as the fitting error of the given scattered points $\{(x_{i},y_{i})\}_{i=1}^{N}\subseteq \Omega$.
To demonstrate the accuracy of the proposed algorithm, we
perform experiments with four functions: a discontinuous function $f_{1}$, a non-smooth function $f_{2}$, a smooth function $f_{3}$ and
the Franke test function $f_{4}$ as follows.
$$f_{1}(x,y)=\left\{
  \begin{array}{cc}
   \frac{x^{2}y}{x^{2}+y^{2}}, & x^{2}+y^{2}\leq 1, \\
   x+y, & x^{2}+y^{2}>1.
  \end{array}
\right.$$
$$f_{2}(x,y)=\left\{
  \begin{array}{cc}
   \frac{xy}{\sqrt{x^{2}+y^{2}}}, & x^{2}+y^{2}\leq 1, \\
   xy, & x^{2}+y^{2}> 1.
  \end{array}
\right.$$
$$f_{3}(x,y)=\frac{1.25+cos(5.4y)}{6+6(3x-1)^{2}}.$$
\begin{eqnarray*}
f_{4}(x,y)&=&0.75\exp[-\frac{(9x-2)^{2}+(9y-2)^{2}}{4}]+0.75\exp[-\frac{(9x+1)^{2}}{49}-\frac{(9y+1)^{2}}{10}]\cr
          &+&0.5\exp[-\frac{(9x-7)^{2}+(9y-3)^{2}}{4}]-0.2\exp[-(9x-4)^{2}-(9y-7)^{2}].
\end{eqnarray*}

\begin{figure}[htbp]
\renewcommand{\figurename}{Fig.}
 \centering
  \subfigure[$J=3$ for MLASSO: $\lambda_{1}=\lambda_{2}=\lambda_{3}=0.001$.]{
  \includegraphics[width=0.3\textwidth]{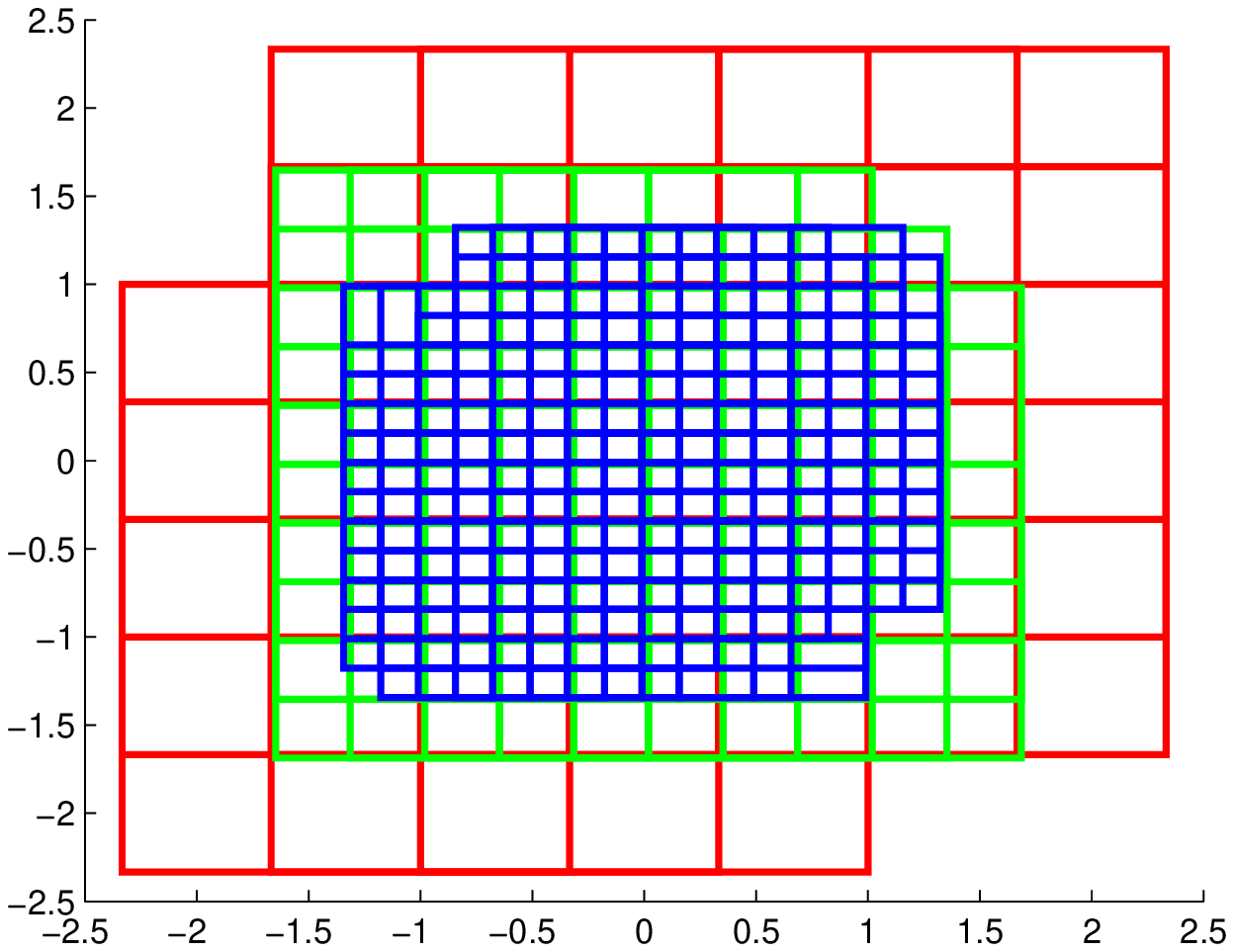}}
  \subfigure[$J=3$ for MLASSO: $\lambda_{1}=\lambda_{2}=\lambda_{3}=0.01$.]{
  \includegraphics[width=0.3\textwidth]{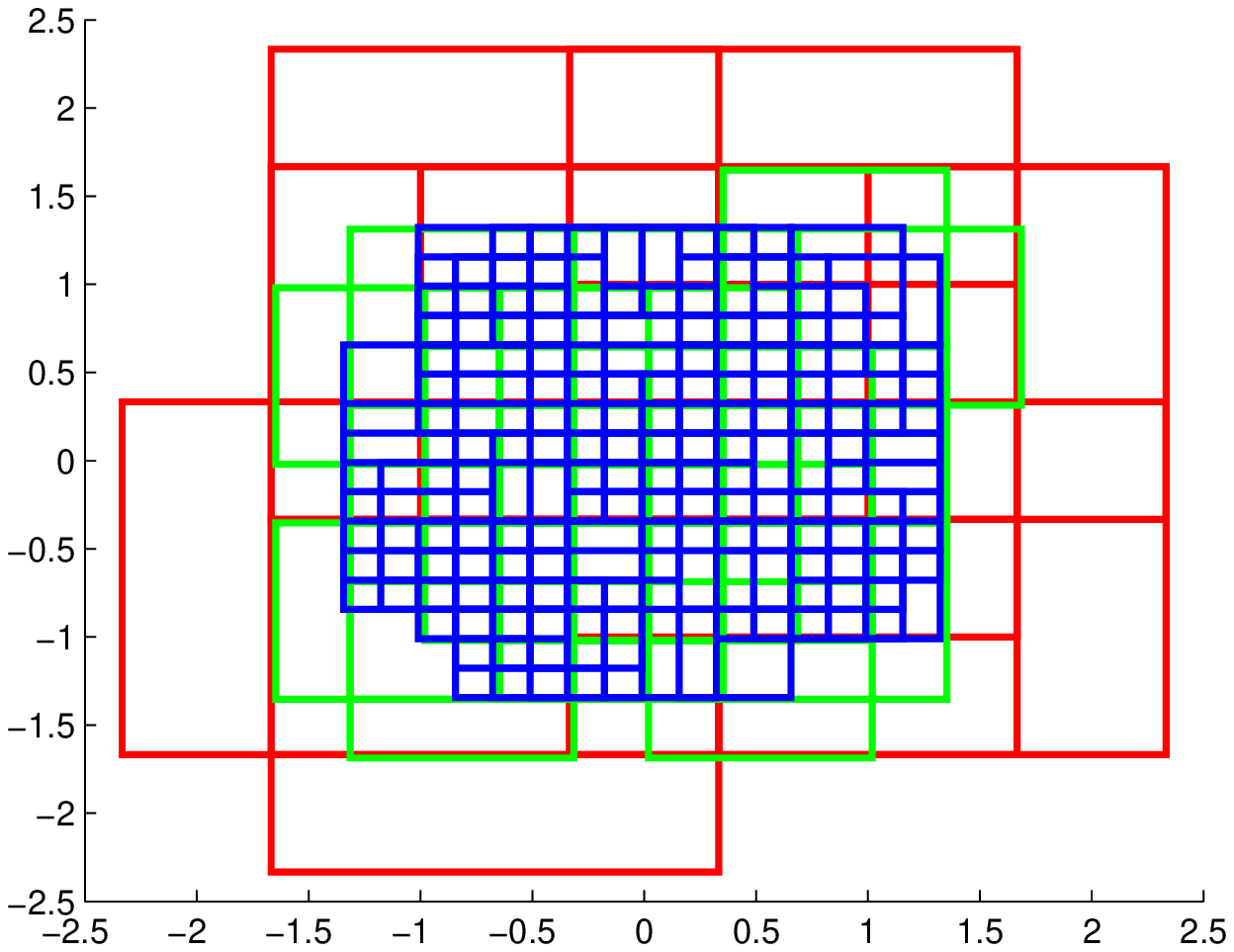}}
  \subfigure[$J=3$ for MLASSO: $\lambda_{1}=0.03, \lambda_{2}=0.01, \lambda_{3}=0.02$.]{
  \includegraphics[width=0.3\textwidth]{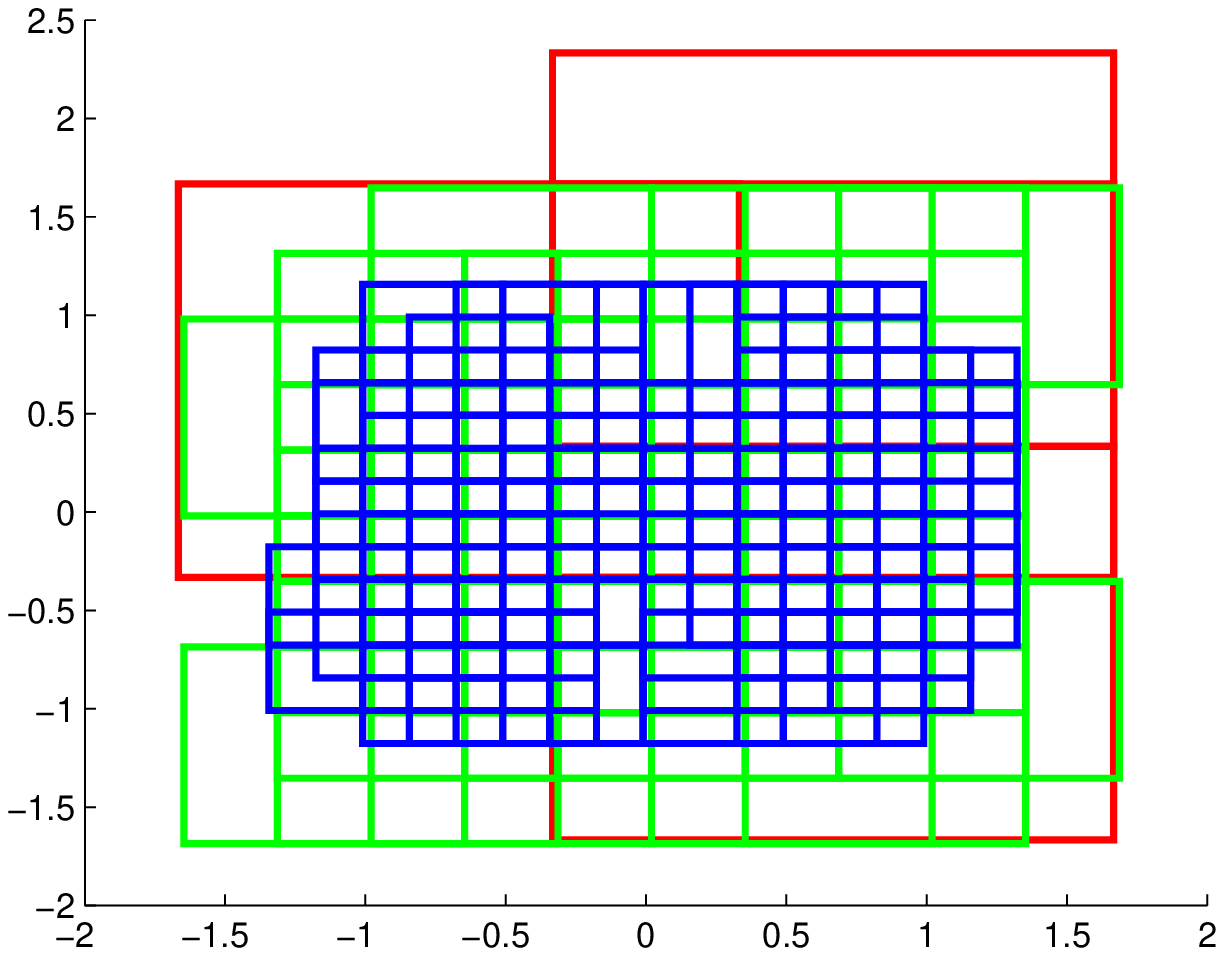}}
  \subfigure[$J=3$ for MLASSO: $\lambda_{1}=0.02,\lambda_{2}=0.01,\lambda_{3}=0.03$.]{
  \includegraphics[width=0.3\textwidth]{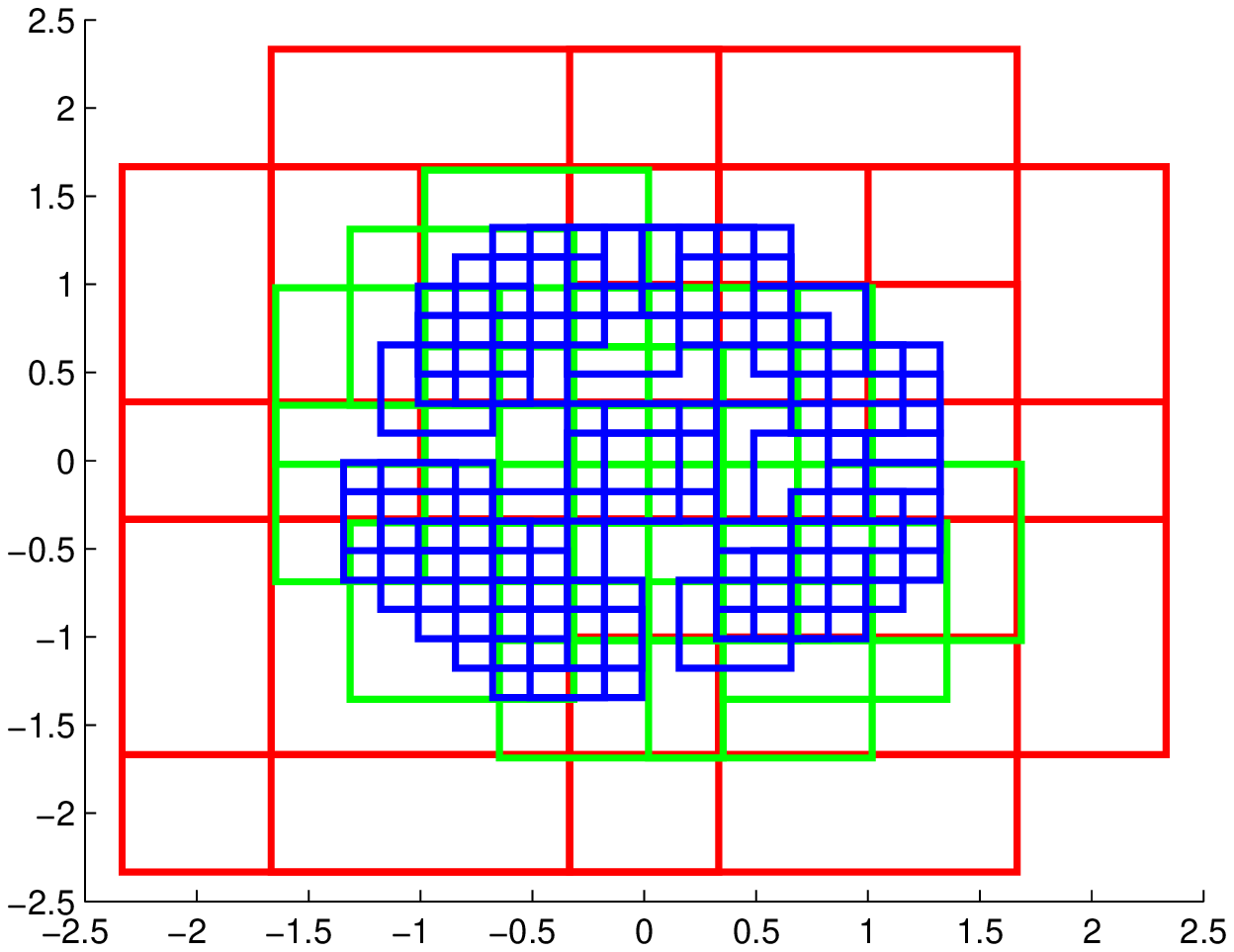}}
  \subfigure[$J=3$ for AGLASSO: $\lambda_{1}=\lambda_{2}=\lambda_{3}=0.001$.]{
  \includegraphics[width=0.3\textwidth]{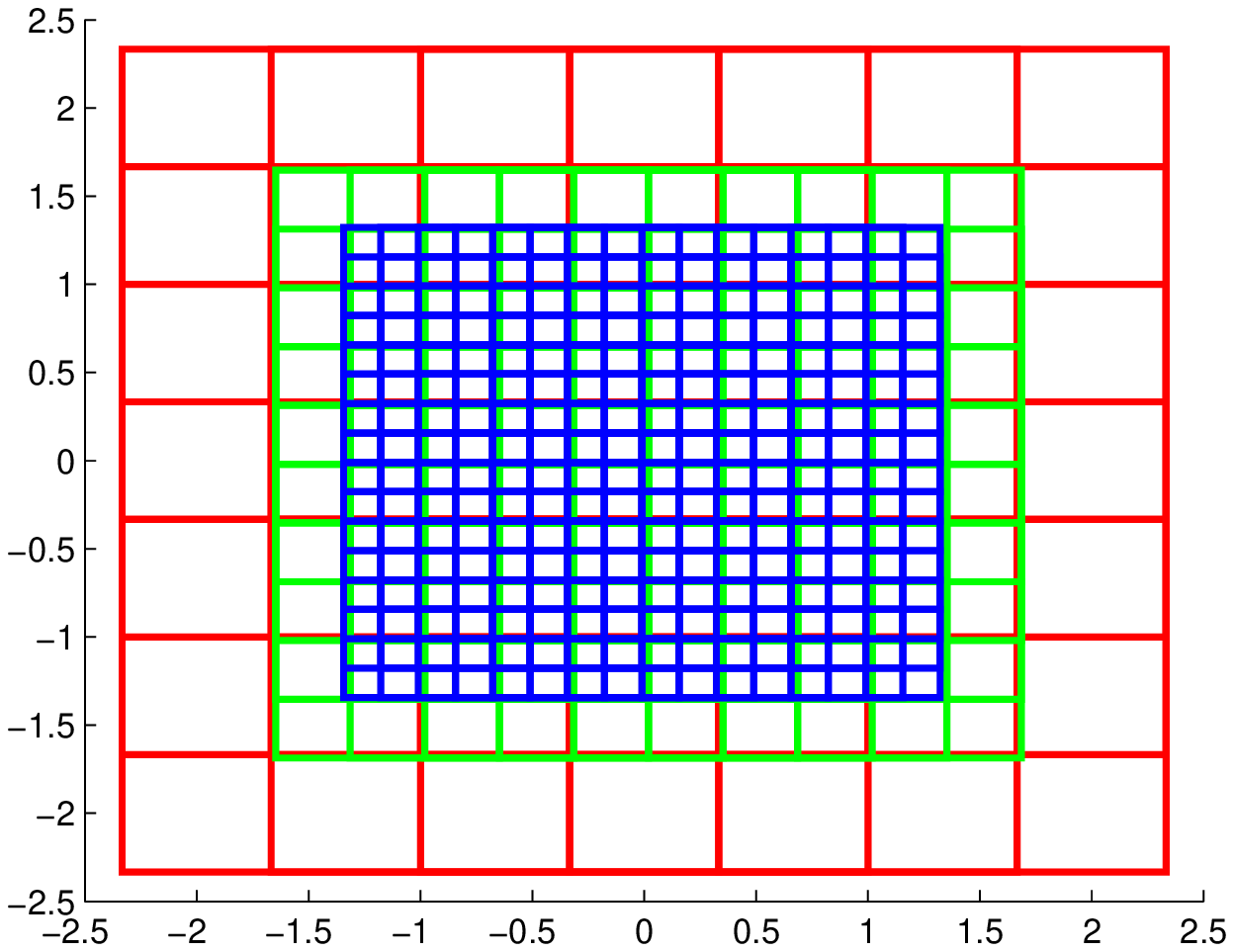}}
  \subfigure[$J=3$ for AGLASSO: $\lambda_{1}=\lambda_{2}=\lambda_{3}=0.01$.]{
  \includegraphics[width=0.3\textwidth]{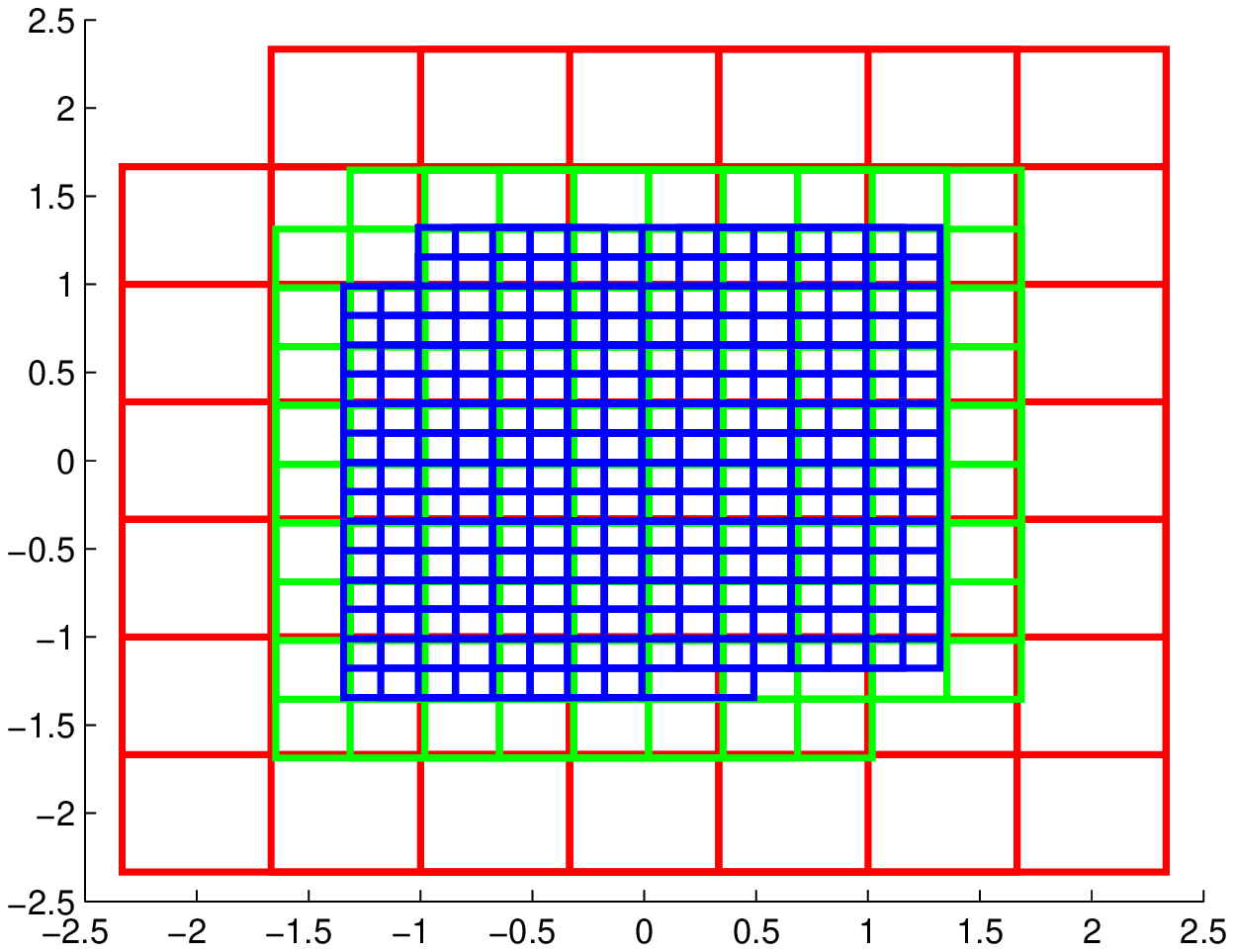}}
  \subfigure[$J=3$ for AGLASSO: $\lambda_{1}=0.03, \lambda_{2}=0.01, \lambda_{3}=0.02$.]{
  \includegraphics[width=0.3\textwidth]{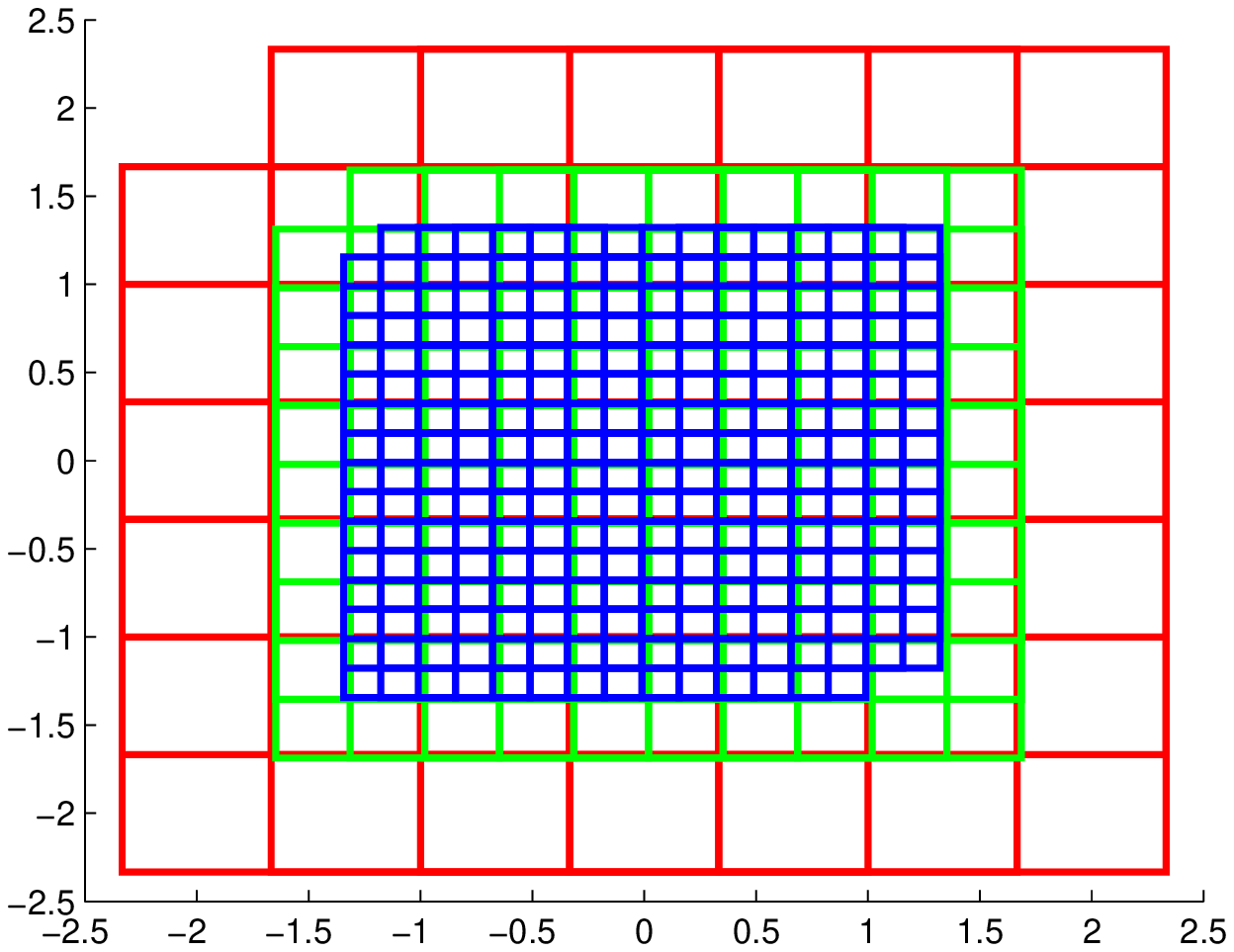}}
  \subfigure[$J=3$ for AGLASSO: $\lambda_{1}=0.02,\lambda_{2}=0.01,\lambda_{3}=0.03$.]{
  \includegraphics[width=0.3\textwidth]{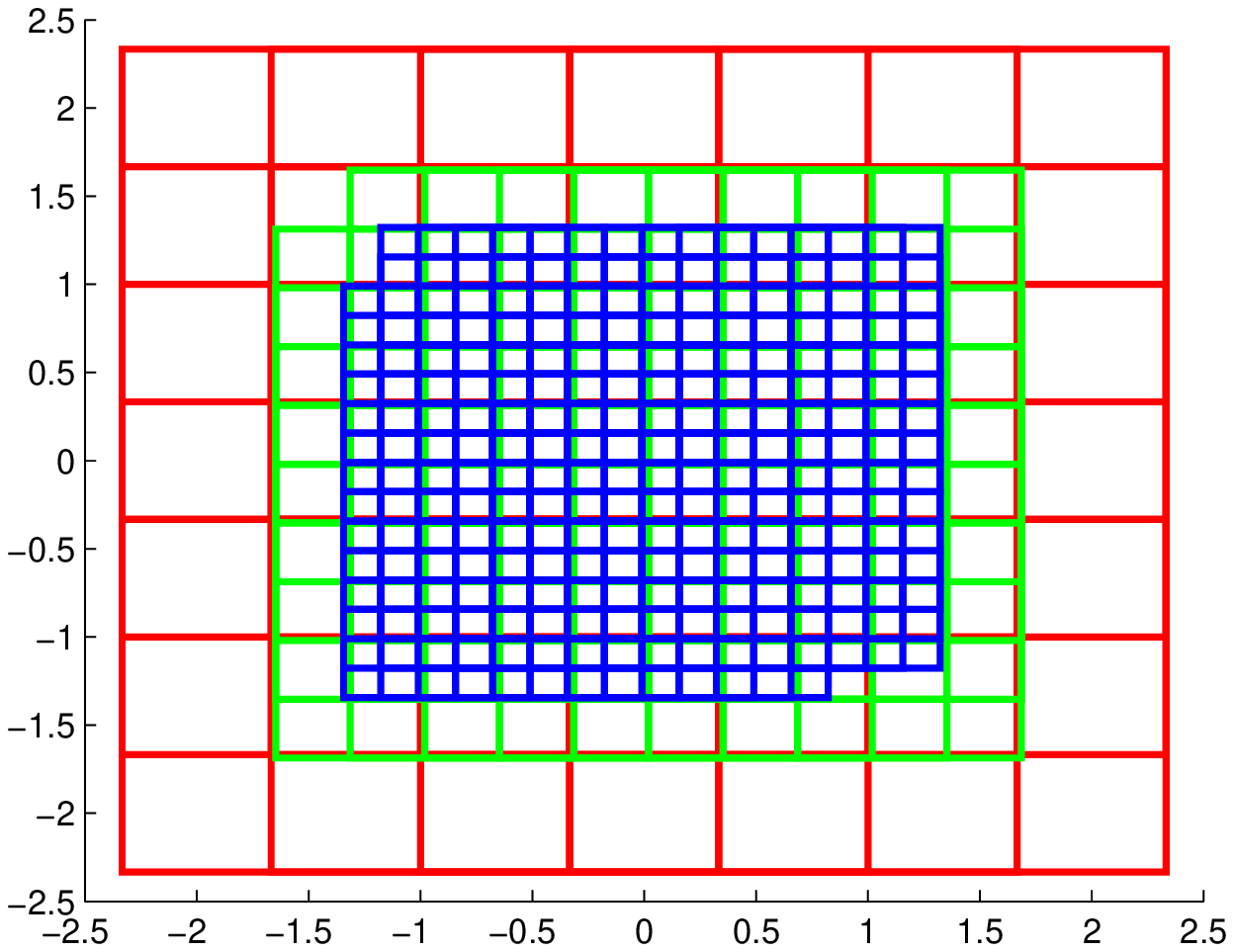}}
  \caption{The distribution of the support for $f_{1}$ with $J=3$.}
  \label{13}
\end{figure}

\begin{figure}[htbp]
\renewcommand{\figurename}{Fig.}
 \centering
  \subfigure[$J=4$ for MLASSO: $\lambda_{1}=\lambda_{2}=\lambda_{3}=\lambda_{4}=0.001$.]{
  \includegraphics[width=0.3\textwidth]{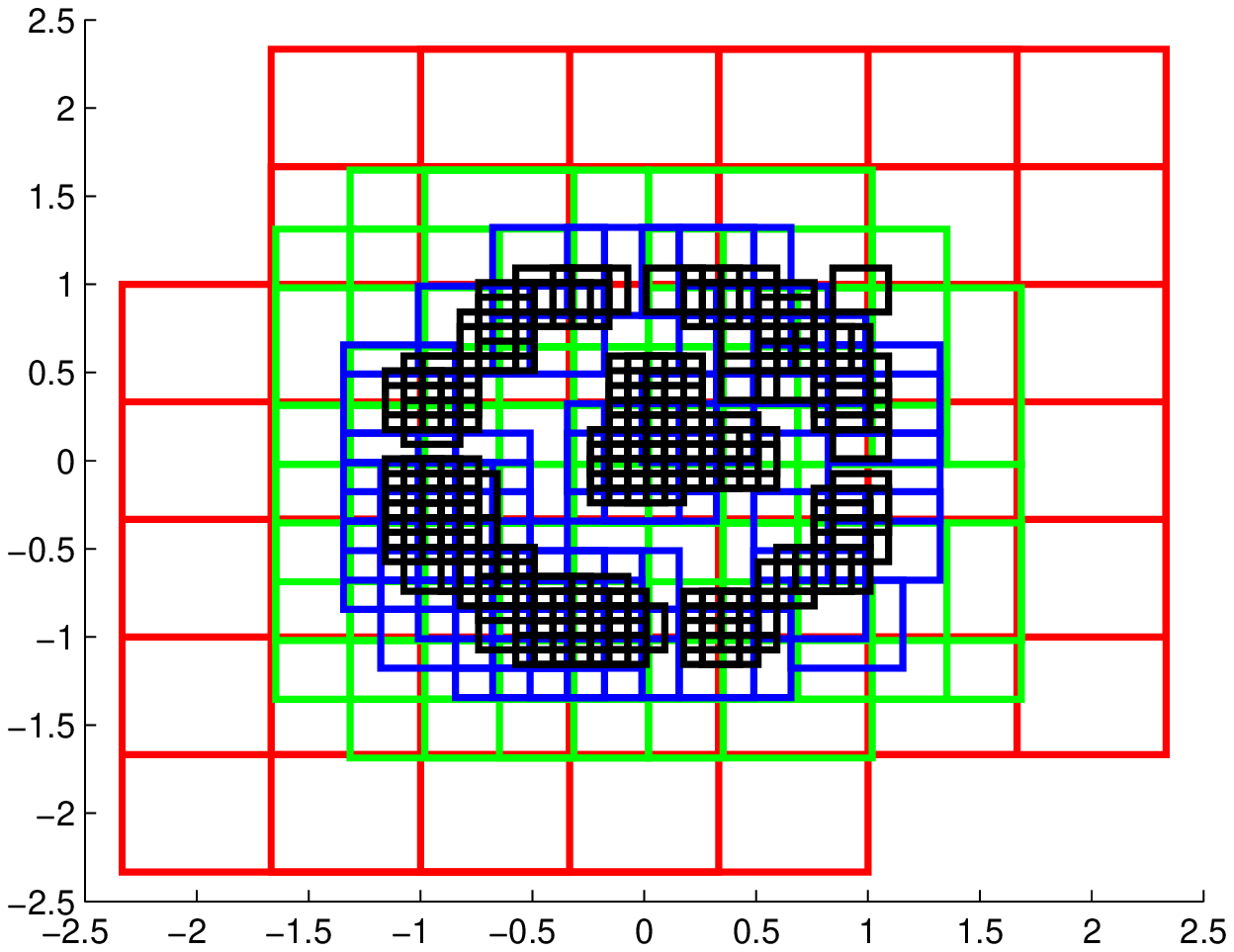}}
  \subfigure[$J=4$ for MLASSO: $\lambda_{1}=\lambda_{2}=\lambda_{3}=\lambda_{4}=0.01$.]{
  \includegraphics[width=0.3\textwidth]{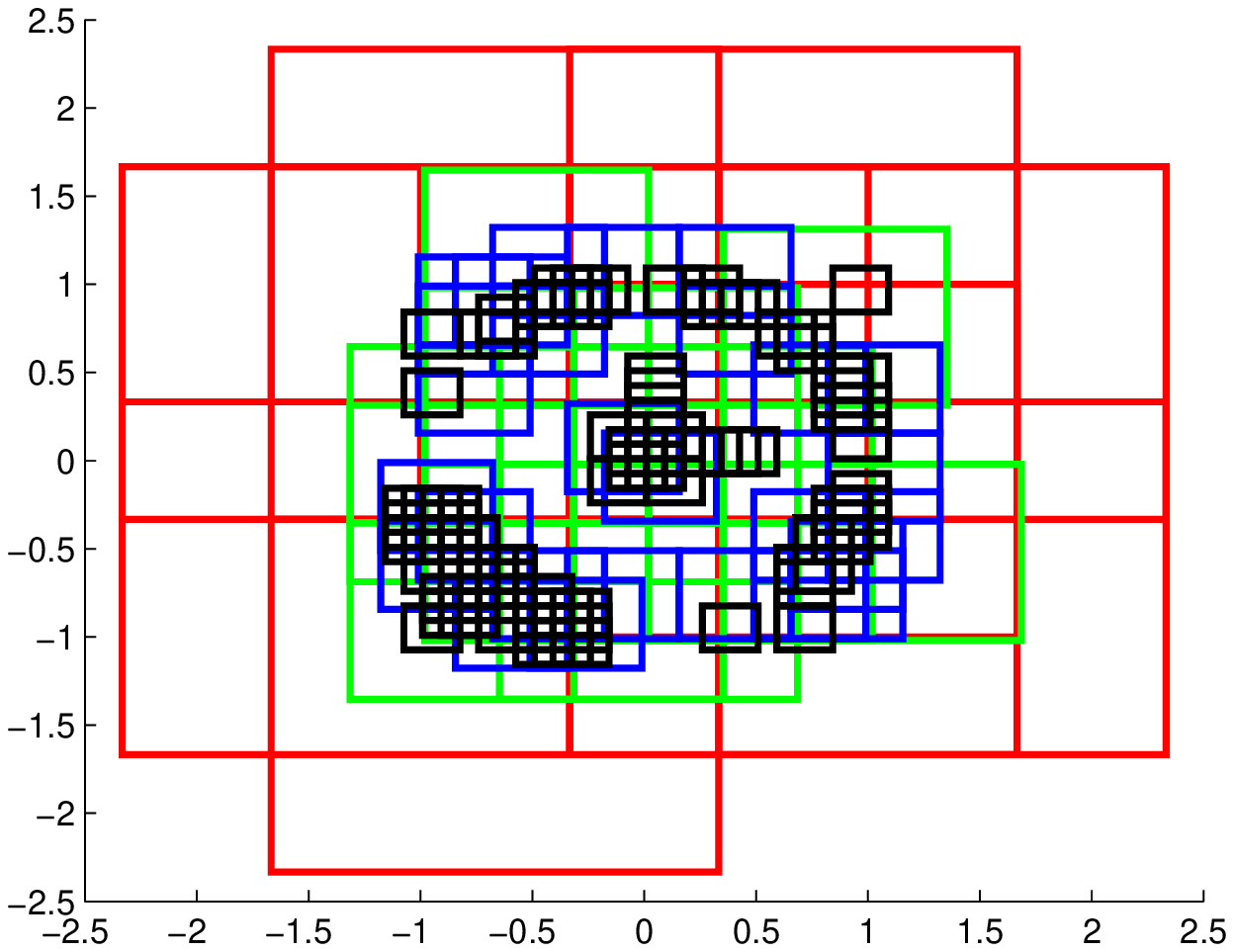}}
  \subfigure[$J=4$ for MLASSO: $\lambda_{1}=\lambda_{4}=0.01, \lambda_{2}=\lambda_{3}=0.001$.]{
  \includegraphics[width=0.3\textwidth]{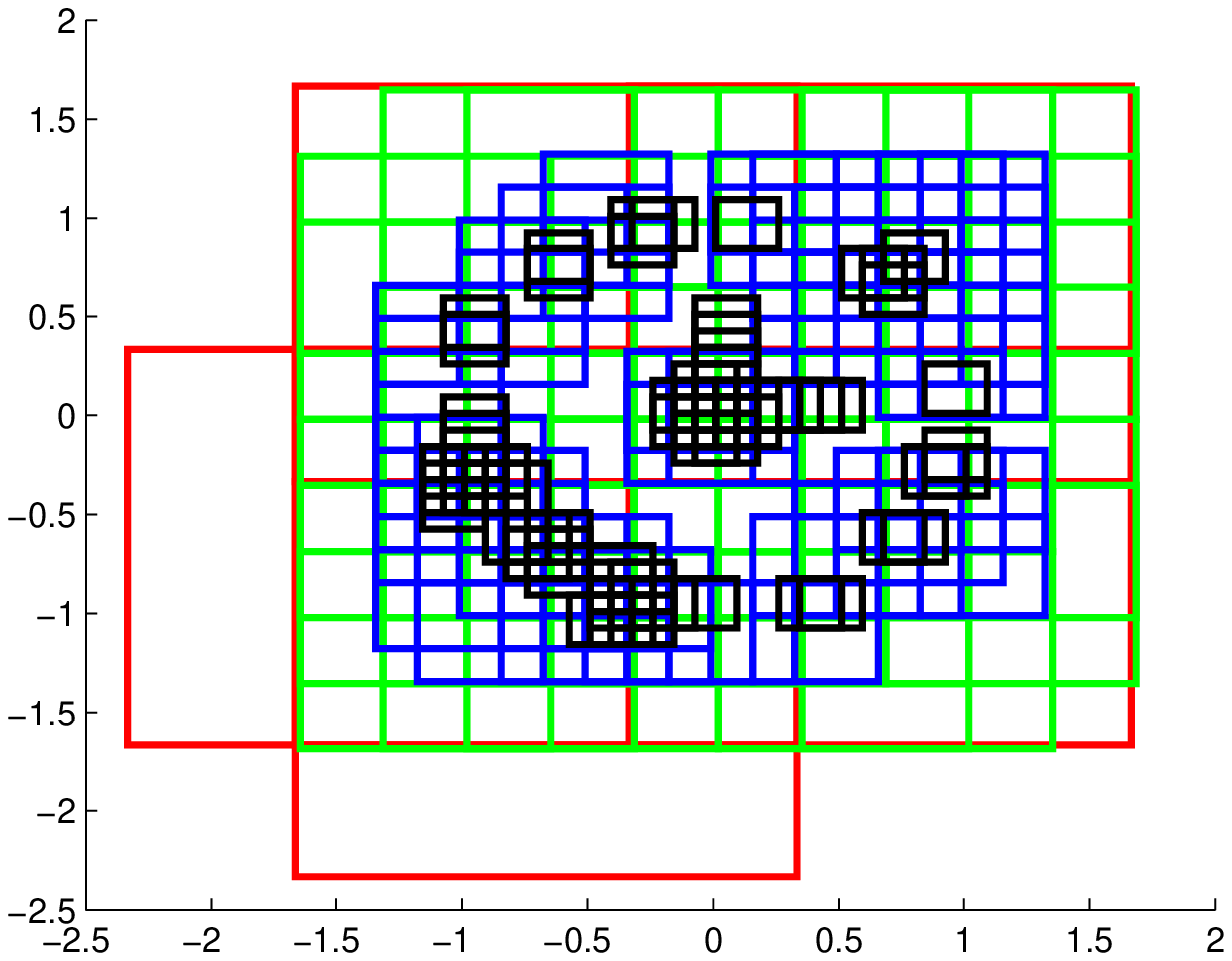}}
  \subfigure[$J=4$ for MLASSO: $\lambda_{1}=\lambda_{4}=0.02, \lambda_{2}=0.001,\lambda_{3}=0.01$.]{
  \includegraphics[width=0.3\textwidth]{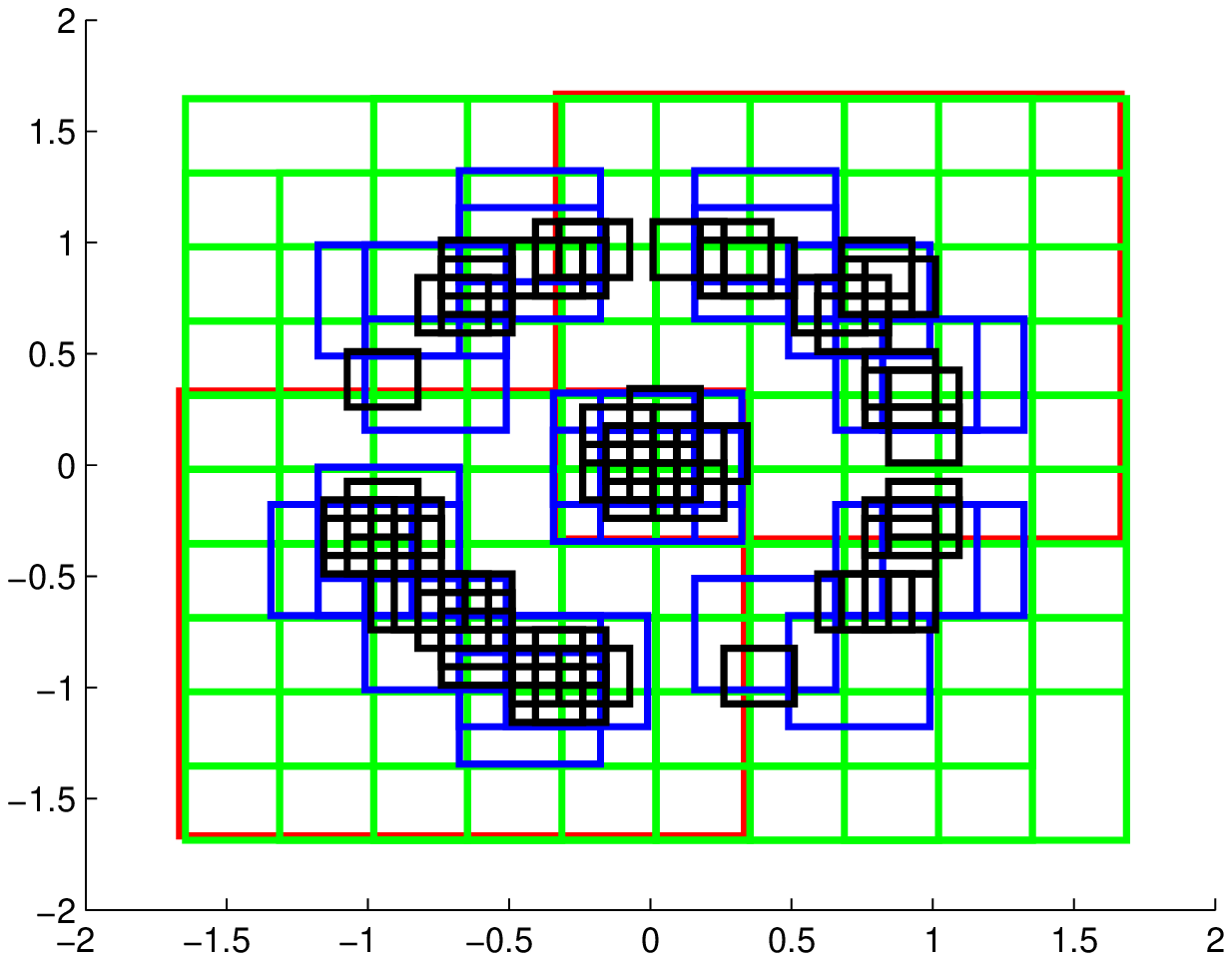}}
  \subfigure[$J=4$ for AGLASSO: $\lambda_{1}=\lambda_{2}=\lambda_{3}=\lambda_{4}=0.001$.]{
  \includegraphics[width=0.3\textwidth]{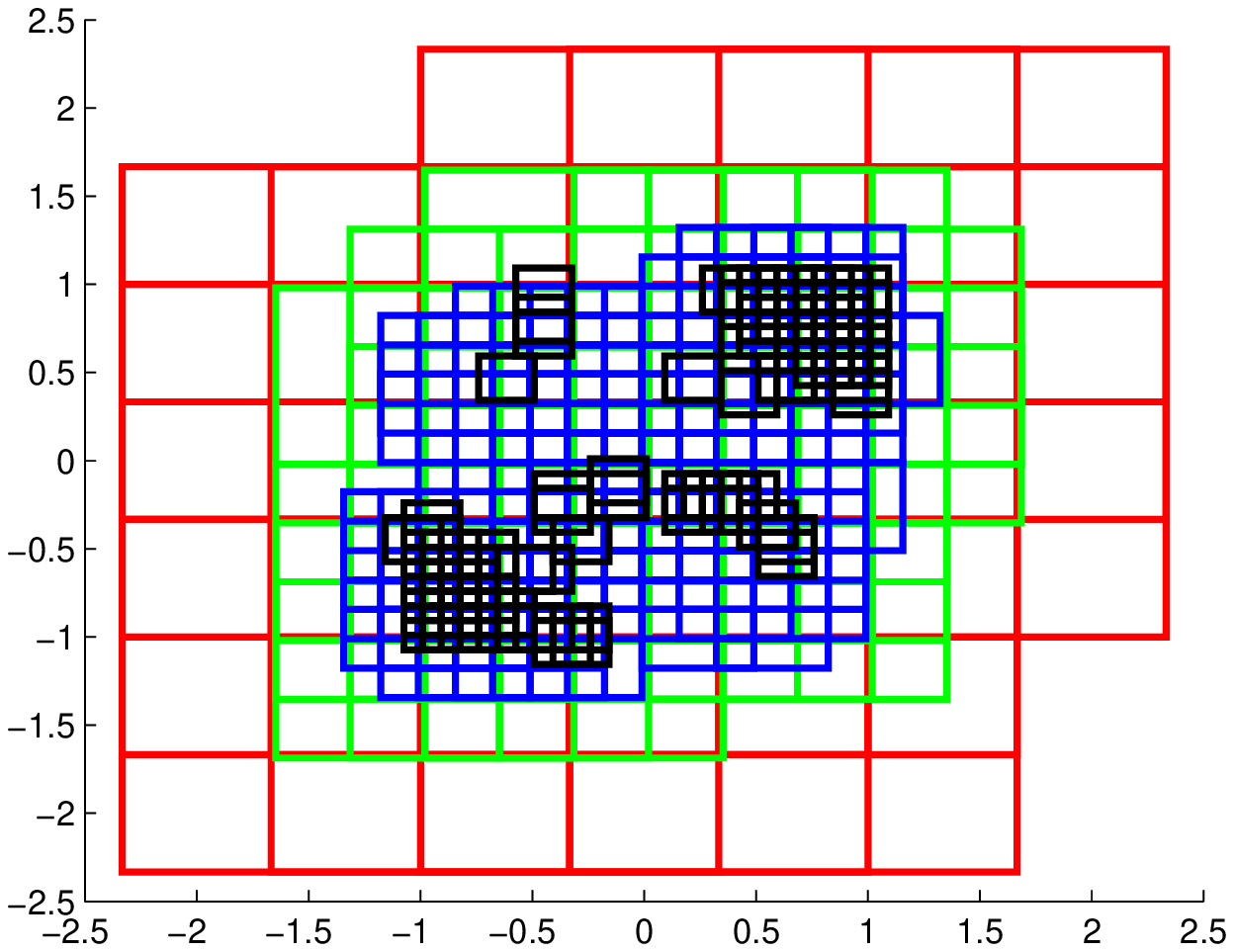}}
  \subfigure[$J=4$ for AGLASSO: $\lambda_{1}=\lambda_{2}=\lambda_{3}=\lambda_{4}=0.01$.]{
  \includegraphics[width=0.3\textwidth]{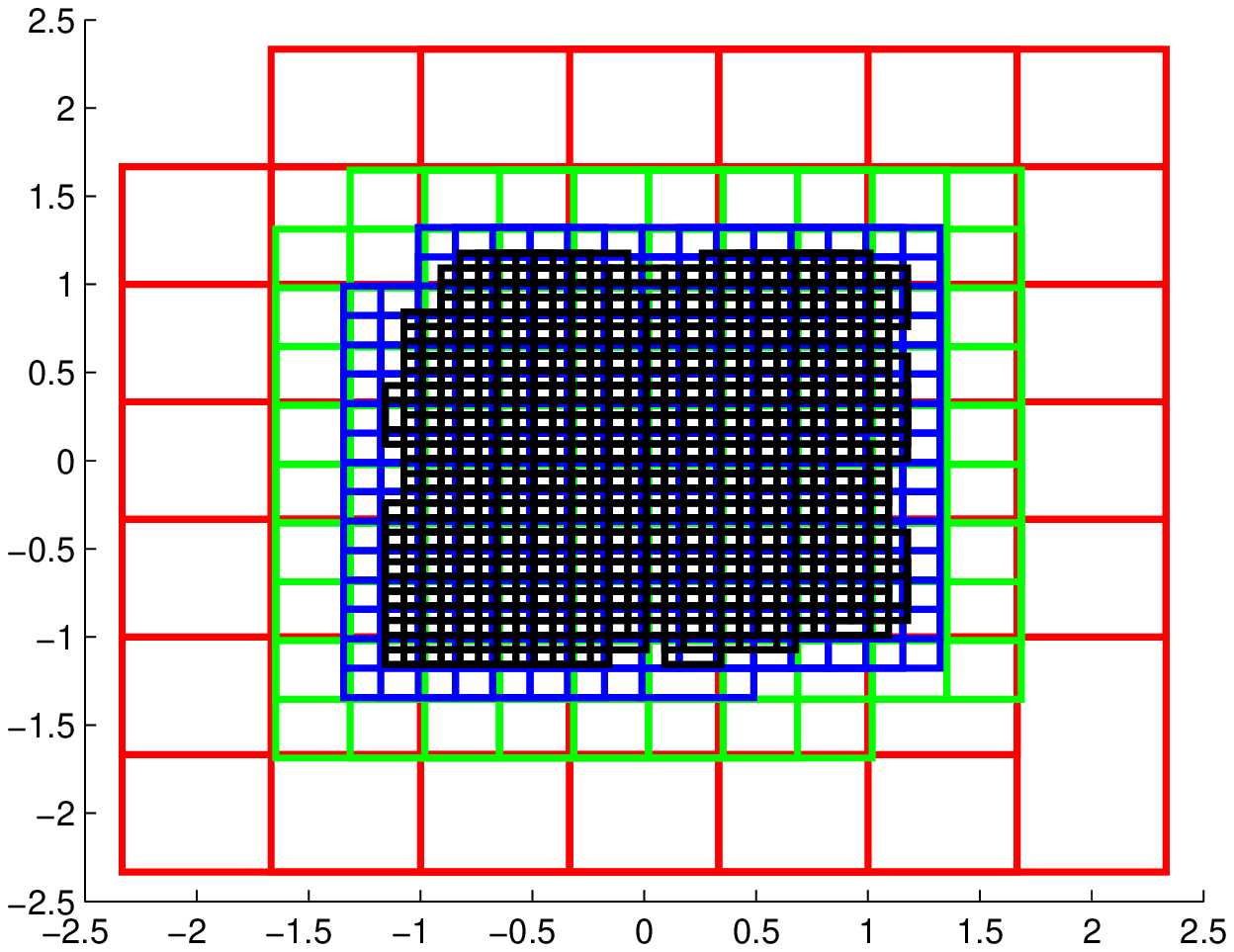}}
  \subfigure[$J=4$ for AGLASSO: $\lambda_{1}=\lambda_{4}=0.01, \lambda_{2}=\lambda_{3}=0.001$.]{
  \includegraphics[width=0.3\textwidth]{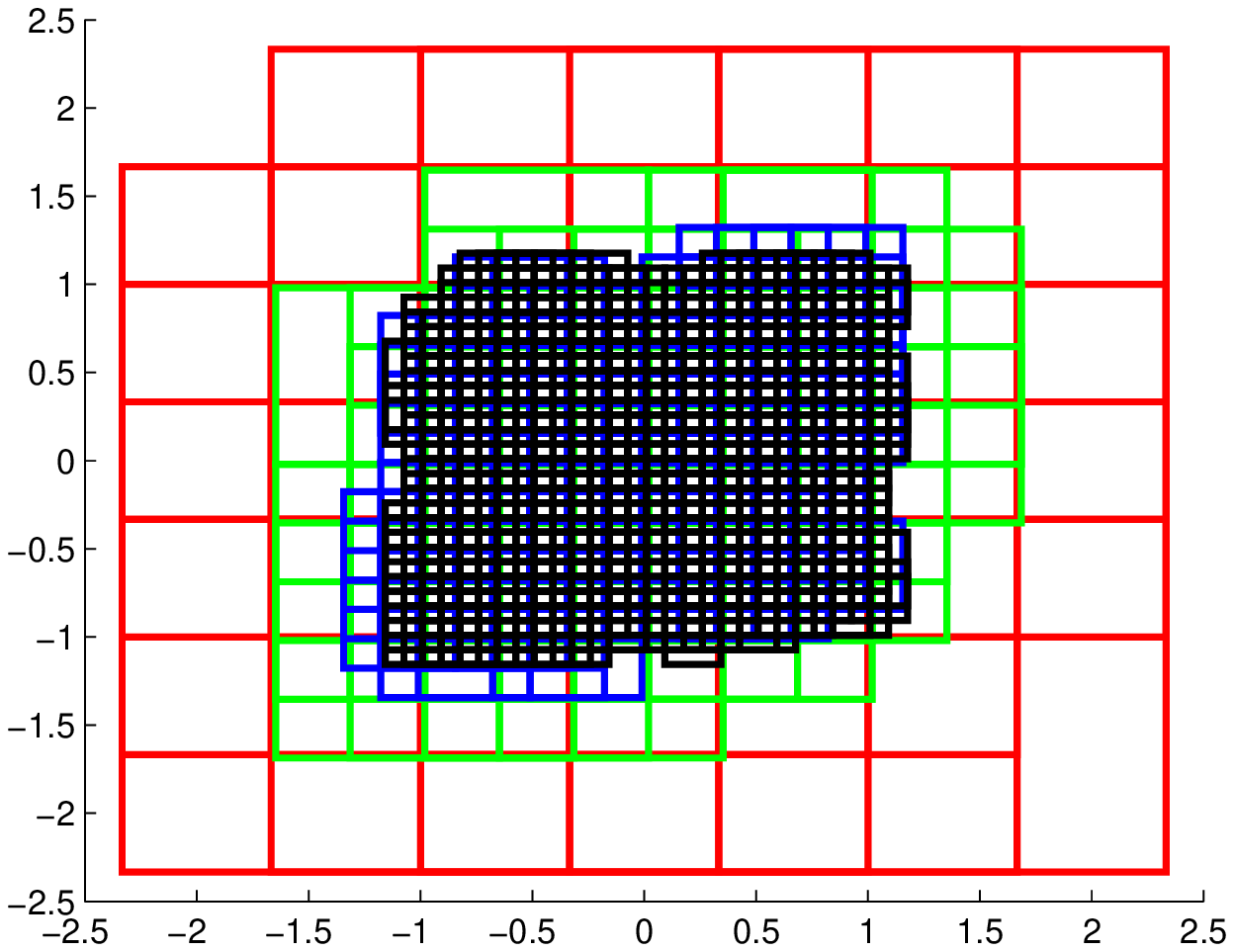}}
  \subfigure[$J=4$ for AGLASSO: $\lambda_{1}=\lambda_{4}=0.02, \lambda_{2}=0.001,\lambda_{3}=0.01$.]{
  \includegraphics[width=0.3\textwidth]{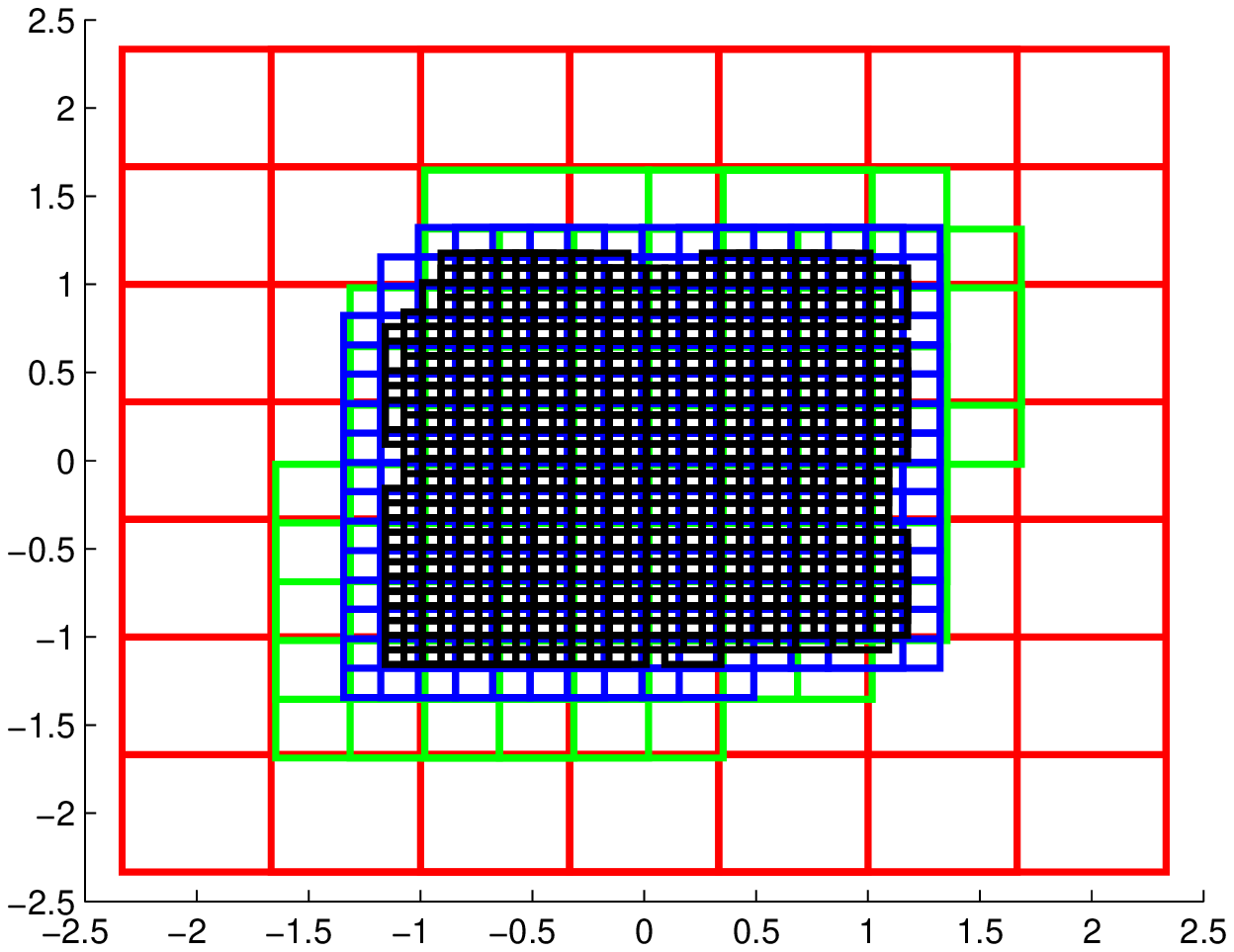}}
  \caption{The distribution of the support for $f_{1}$ with $J=4$ . }
  \label{14}
\end{figure}

\begin{table}[htbp]
\center
\caption{The $l_{0}$ norm and the approximation errors for different parameters shown in Fig.\ref{13} for $f_{1}$ with $J=3$.}\label{Tab:13}
\begin{tabular}{cccccc}\hline
                   & $(l_{0}(\textbf{X}_{1}),l_{0}(\textbf{X}_{2}),l_{0}(\textbf{X}_{3}))$  & Error & RMS & Iterations  & Time(sec) \\ \hline
\cite{Lee}         &   (24, 64, 193)      &   2.1241e-3    &  2.3512e-3   & 3      & 0.0241   \\
Fig.\ref{13}(a)    &   (18, 53, 143)      &   2.2538e-3    &  4.3462e-3   & 17     & 0.0431   \\
Fig.\ref{13}(b)    &   (13, 29, 84)       &   3.8215e-3    &  5.6138e-3   & 1322   & 2.0032    \\
Fig.\ref{13}(c)    &   (4, 41, 88)        &   5.1433e-3    &  7.3768e-3   & 1841   & 2.4456    \\
Fig.\ref{13}(d)    &   (15, 22, 49)       &   5.3102e-3    &  8.7471e-3   & 2436   & 3.2345    \\
Fig.\ref{13}(e)    &   (20, 51, 116)      &   5.9973e-1    &  6.2737e-1   & 1      & 0.3917    \\
Fig.\ref{13}(f)    &   (24, 61, 181)      &   5.0691e-1    &  5.4132e-1   & 1      & 0.2961    \\
Fig.\ref{13}(g)    &   (24, 63, 190)      &   4.2005e-1    &  4.6263e-1   & 1      & 0.3524    \\
Fig.\ref{13}(h)    &   (25, 62, 189)      &   3.9211e-1    &  4.3247e-1   & 1      & 0.3012    \\
\hline
\end{tabular}
\end{table}

\begin{table}[htbp]
\center
\caption{The $l_{0}$ norm and the approximation errors for different parameters shown in Fig.\ref{14} for $f_{1}$ with $J=4$.}\label{Tab:14}
\begin{tabular}{cccccc}\hline
                   & $(l_{0}(\textbf{X}_{1}),l_{0}(\textbf{X}_{2}),l_{0}(\textbf{X}_{3}),l_{0}(\textbf{X}_{4}))$  & Error & RMS & Iterations  & Time(sec) \\ \hline
\cite{Lee}         &   (25, 62, 192, 674)     &   1.1185e-3    &  3.3243e-4   & 4      & 0.2191   \\
Fig.\ref{14}(a)    &   (17, 37, 50, 138)      &   1.4298e-3    &  3.5873e-3   & 267    & 4.3145   \\
Fig.\ref{14}(b)    &   (15, 21, 27, 88)       &   2.1073e-3    &  4.2481e-3   & 1209   & 11.1452    \\
Fig.\ref{14}(c)    &   (7, 41, 69, 63)        &   3.1132e-3    &  5.6871e-3   & 1307   & 13.0139    \\
Fig.\ref{14}(d)    &   (2, 53, 27, 59)        &   3.9413e-3    &  6.1678e-3   & 1183   & 12.7346    \\
Fig.\ref{14}(e)    &   (20, 52, 114, 90)      &   5.2471e-1    &  6.0124e-1   & 1      & 5.1328    \\
Fig.\ref{14}(f)    &   (23, 61, 181, 537)     &   5.0109e-1    &  4.7453e-1   & 1      & 5.0129    \\
Fig.\ref{14}(g)    &   (24, 48, 112, 544)     &   5.1902e-1    &  5.4287e-1   & 1      & 4.9879    \\
Fig.\ref{14}(h)    &   (25, 45, 176, 581)     &   3.8716e-1    &  5.2634e-1   & 1      & 4.3472    \\
\hline
\end{tabular}
\end{table}

\begin{figure}[htbp]
\renewcommand{\figurename}{Fig.}
  \centering
  \subfigure[$J=3$ for MLASSO: $\lambda_{1}=\lambda_{2}=\lambda_{3}=0.001$.]{
  \includegraphics[width=0.3\textwidth]{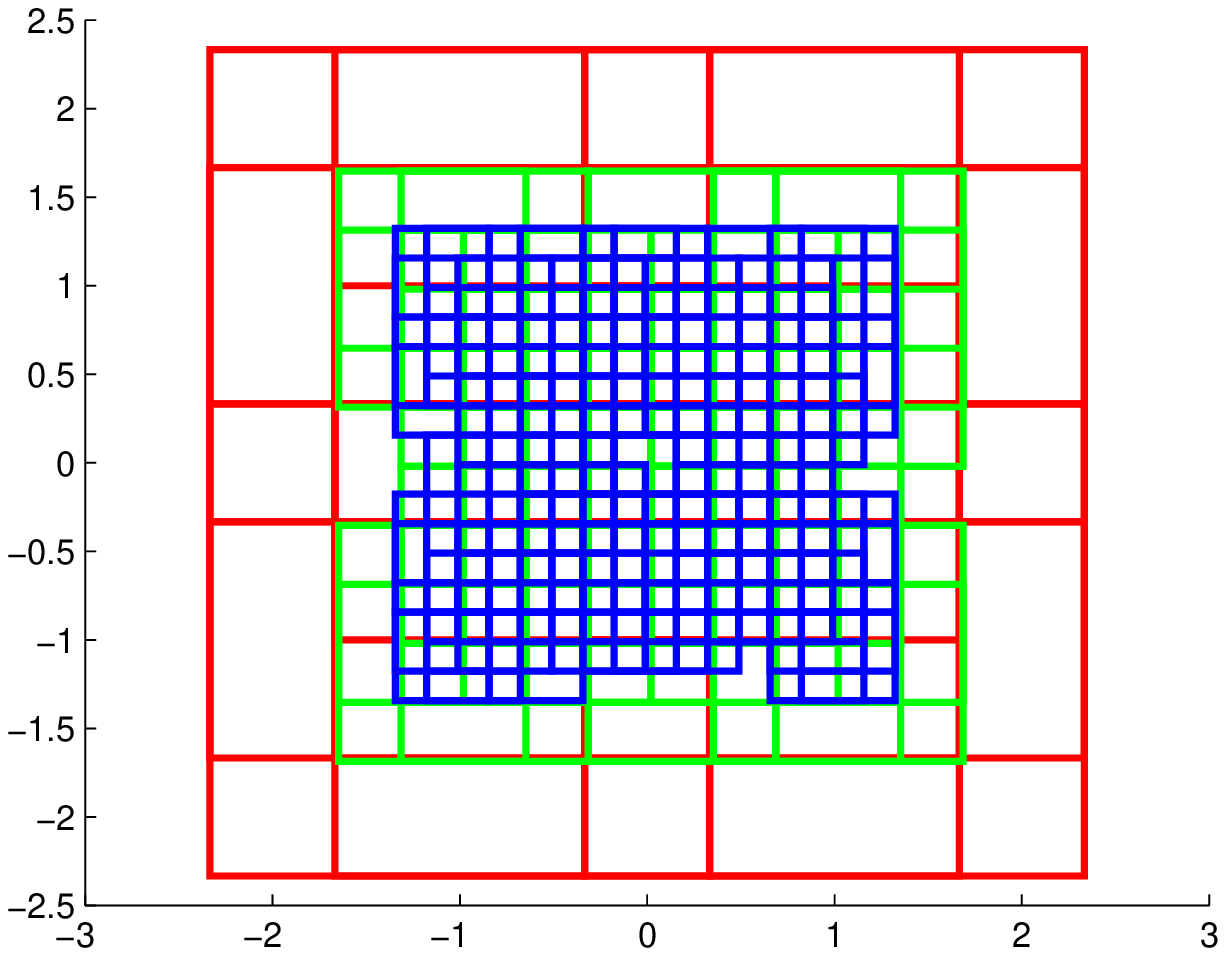}}
  \subfigure[$J=3$ for MLASSO: $\lambda_{1}=\lambda_{2}=\lambda_{3}=0.01$.]{
  \includegraphics[width=0.3\textwidth]{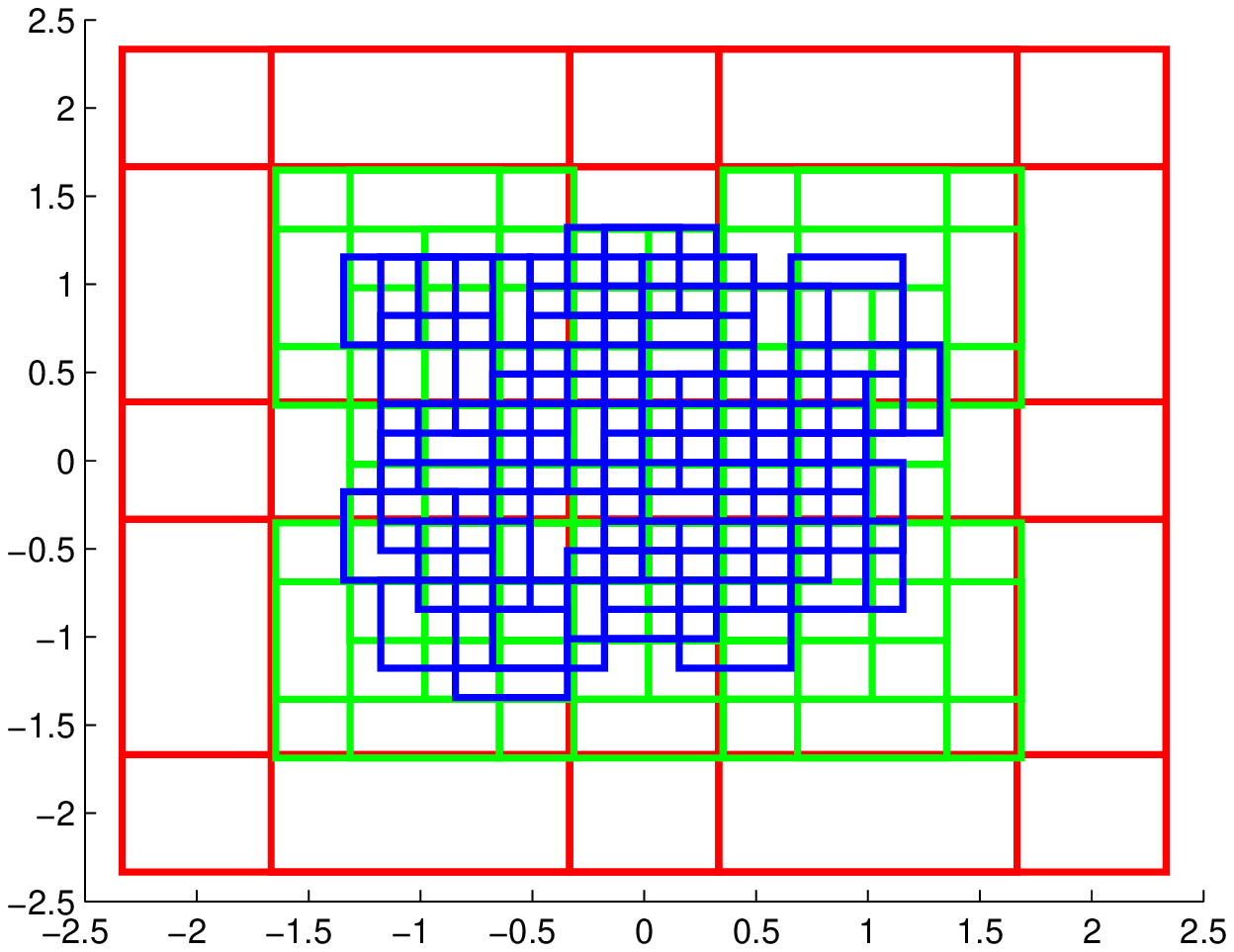}}
  \subfigure[$J=3$ for MLASSO: $\lambda_{1}=0.02,\lambda_{2}=0.01,\lambda_{3}=0.03$.]{
  \includegraphics[width=0.3\textwidth]{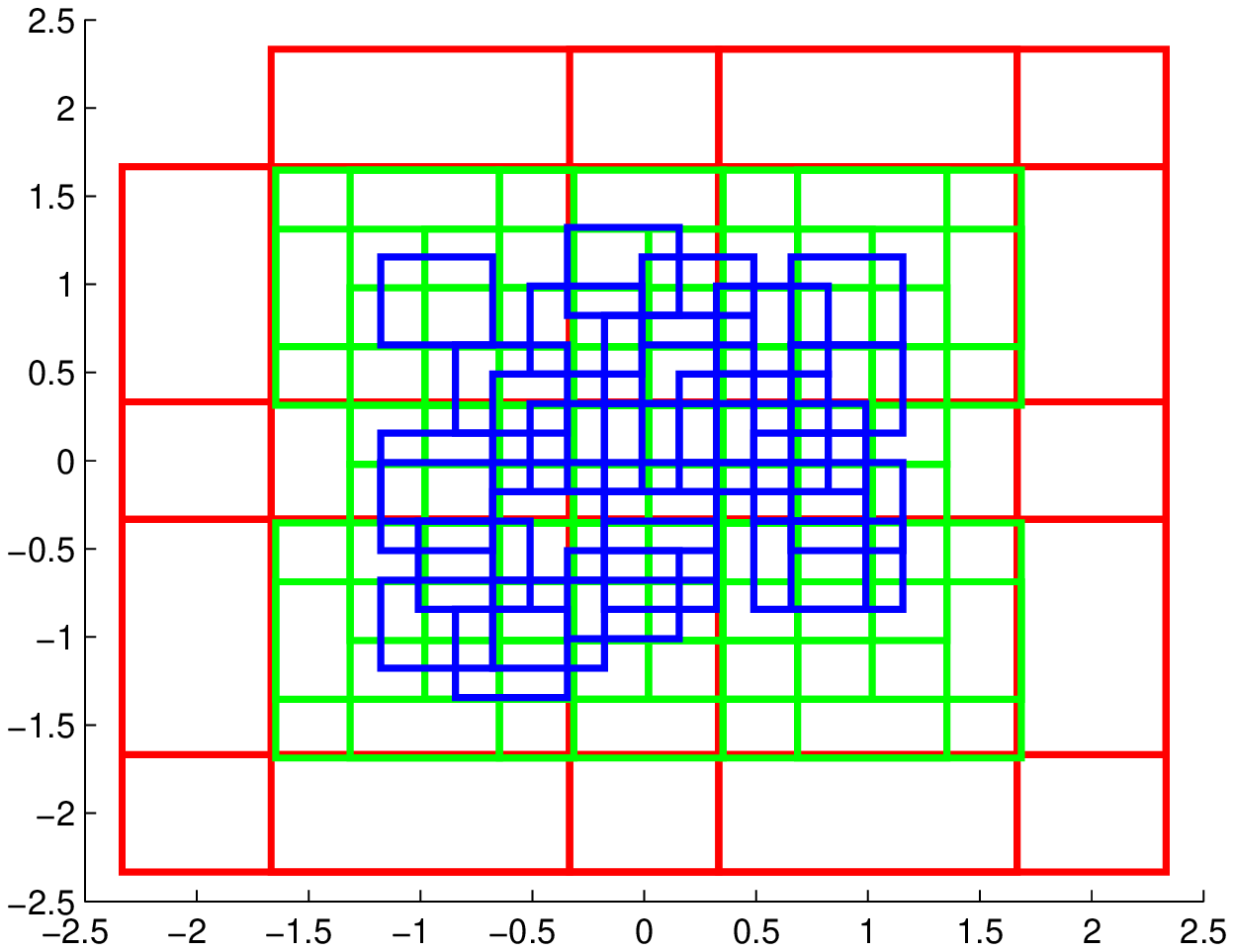}}
  \subfigure[$J=3$ for MLASSO: $\lambda_{1}=0.03, \lambda_{2}=0.01,\lambda_{3}=0.02$.]{
  \includegraphics[width=0.3\textwidth]{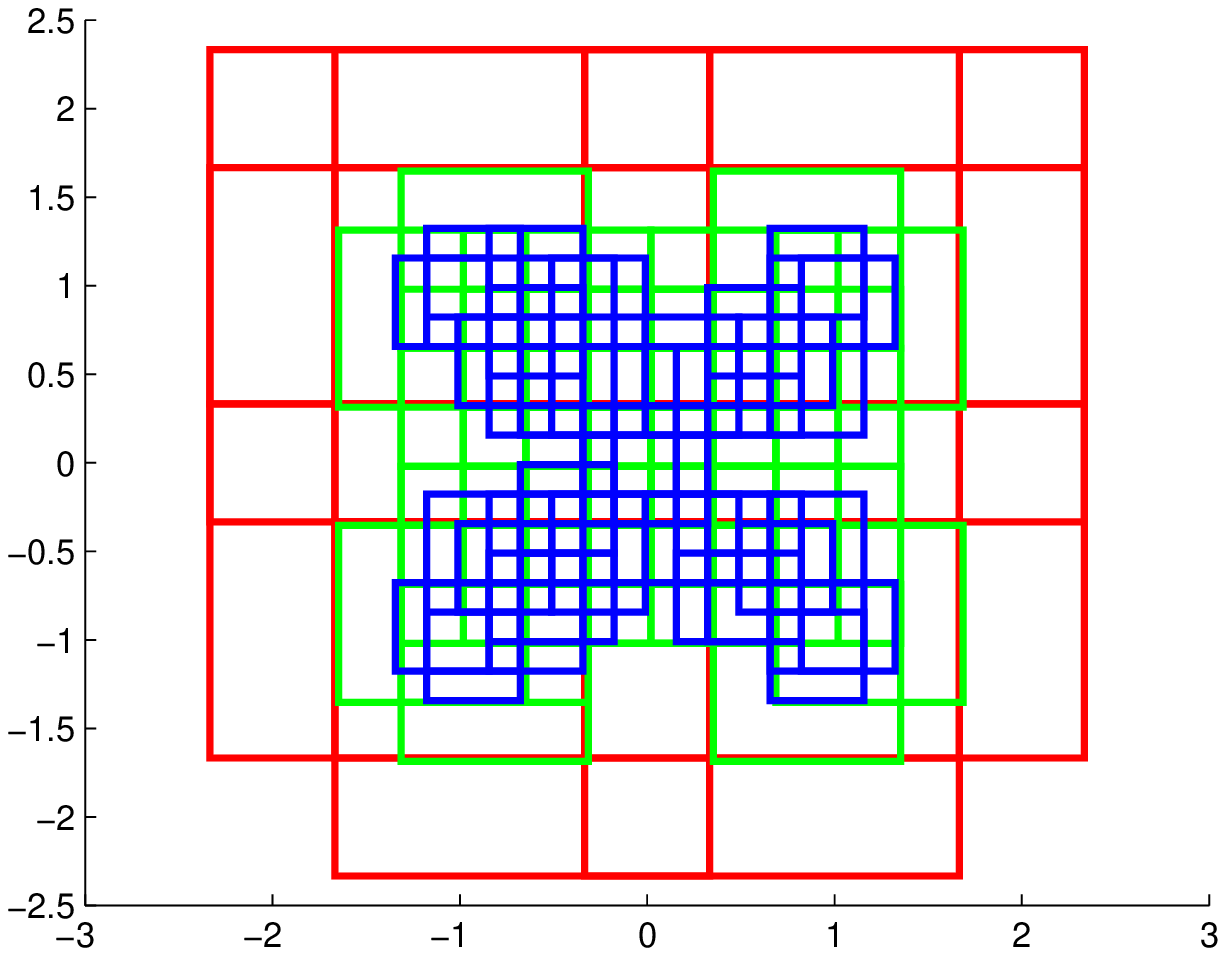}}
  \subfigure[$J=3$ for AGLASSO: $\lambda_{1}=\lambda_{2}=\lambda_{3}=0.001$.]{
 \includegraphics[width=0.3\textwidth]{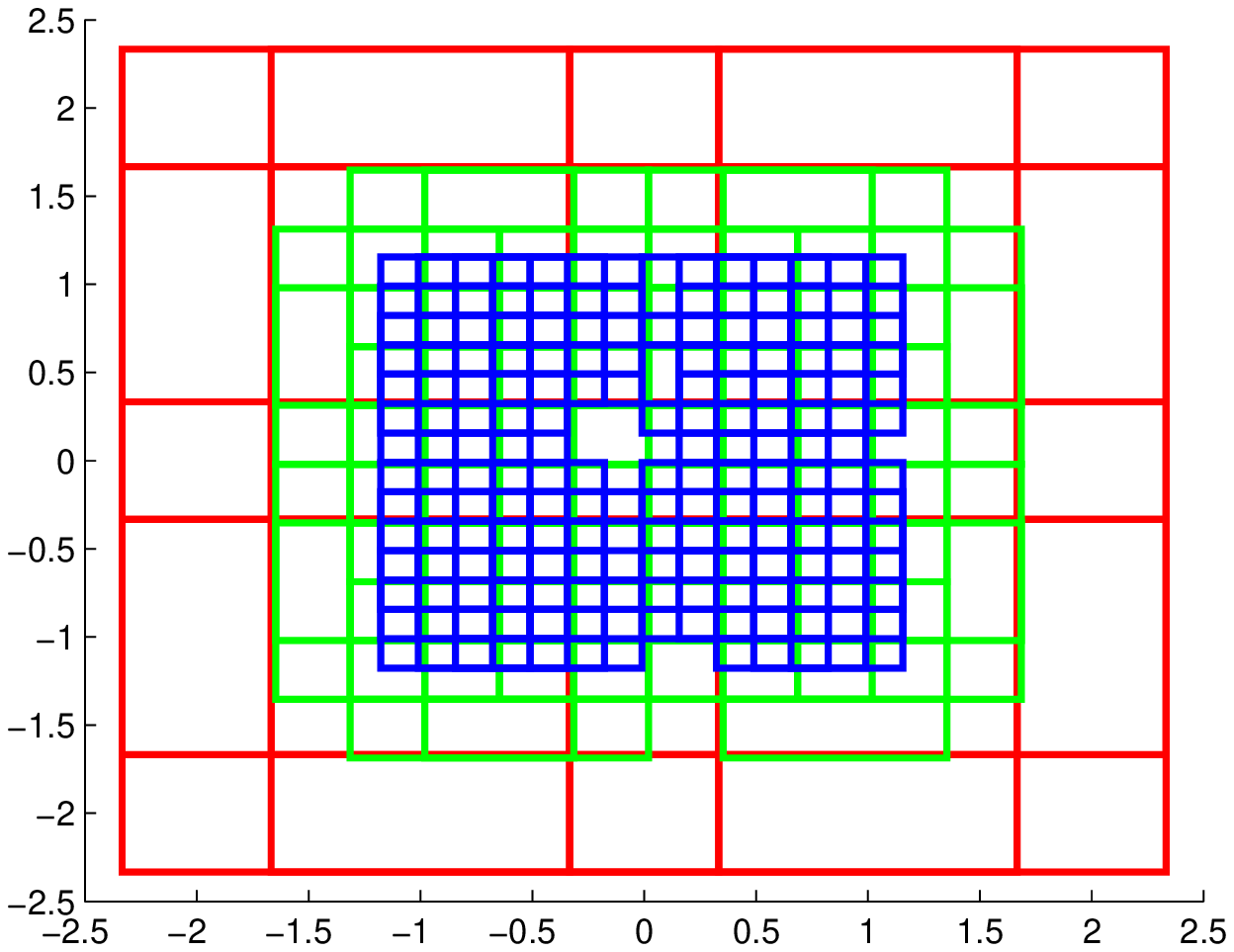}}
  \subfigure[$J=3$ for AGLASSO: $\lambda_{1}=\lambda_{2}=\lambda_{3}=0.01$.]{
  \includegraphics[width=0.3\textwidth]{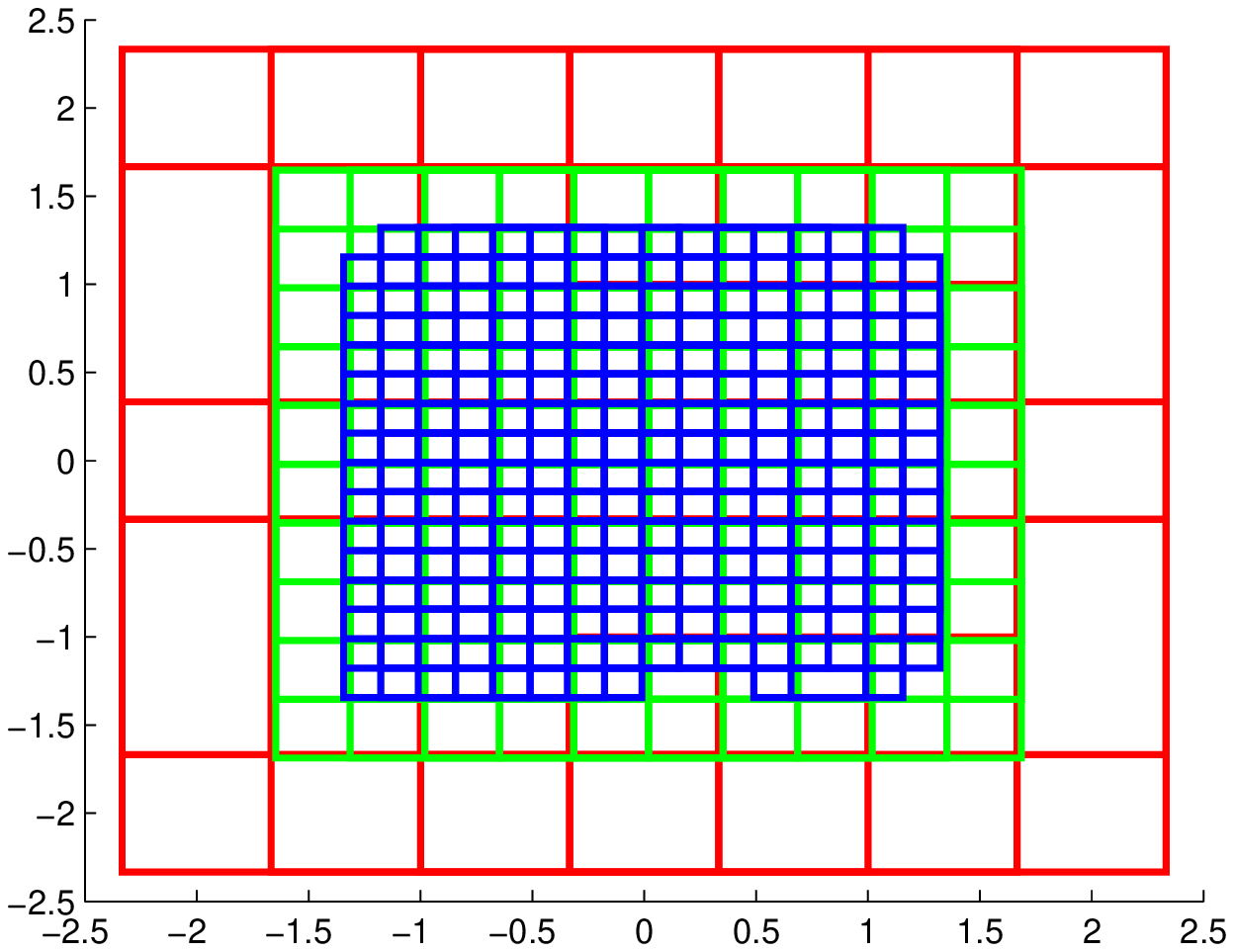}}
  \subfigure[$J=3$ for AGLASSO: $\lambda_{1}=0.02,\lambda_{2}=0.01,\lambda_{3}=0.03$.]{
  \includegraphics[width=0.3\textwidth]{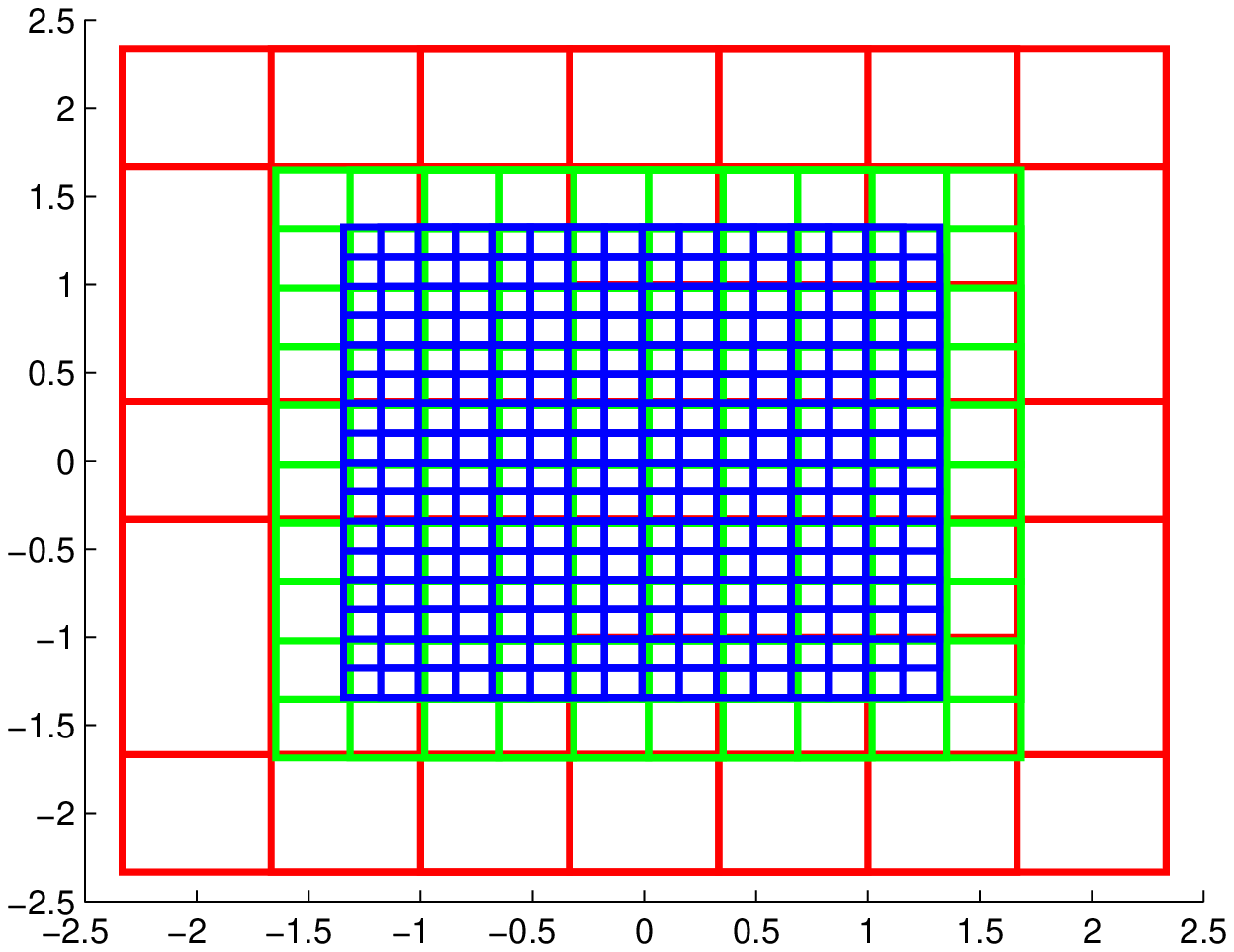}}
  \subfigure[$J=3$ for AGLASSO: $\lambda_{1}=0.03, \lambda_{2}=0.01,\lambda_{3}=0.02$.]{
  \includegraphics[width=0.3\textwidth]{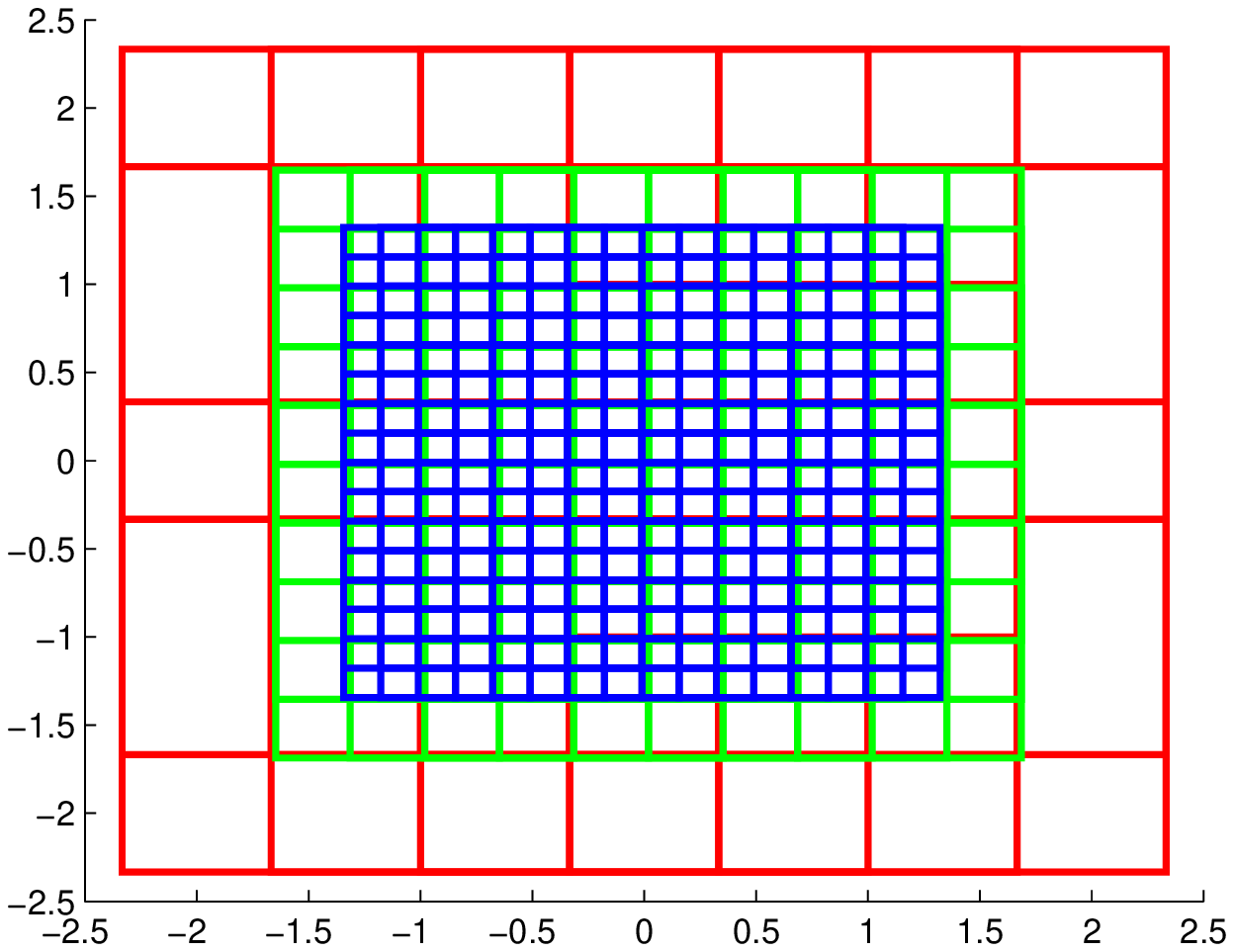}}
  \caption{The distribution of the support for $f_{2}$ with $J=3$. }
  \label{23}
\end{figure}

\begin{table}[htbp]
\center
\caption{The $l_{0}$ norm and the approximation errors for different parameters shown in Fig.\ref{23} for $f_{2}$ with $J=3$.}\label{Tab:23}
\begin{tabular}{cccccc}\hline
                   & $(l_{0}(\textbf{X}_{1}),l_{0}(\textbf{X}_{2}),l_{0}(\textbf{X}_{3}))$  & Error & RMS & Iterations  & Time(sec) \\ \hline
\cite{Lee}         &   (22, 64, 194)      &   3.7541e-3    &  7.3274e-3   & 3      & 0.1102   \\
Fig.\ref{23}(a)    &   (18, 41, 115)      &   4.1637e-3    &  8.3761e-3   & 232    & 0.4257   \\
Fig.\ref{23}(b)    &   (16, 41, 53)       &   5.0749e-3    &  9.0652e-3   & 678    & 2.3091    \\
Fig.\ref{23}(c)    &   (15, 40, 31)       &   4.7049e-3    &  9.3981e-3   & 1103   & 1.6027    \\
Fig.\ref{23}(d)    &   (14, 33, 49)       &   5.2564e-3    &  9.2971e-3   & 1089   & 1.5453    \\
Fig.\ref{23}(e)    &   (15, 51, 91)       &   3.3797e-1    &  3.5205e-1   & 1      & 1.0274    \\
Fig.\ref{23}(f)    &   (20, 62, 183)      &   2.2143e-1    &  2.6308e-1   & 1      & 0.7348    \\
Fig.\ref{23}(g)    &   (21, 64, 194)      &   2.0213e-1    &  2.1453e-1   & 1      & 0.5762    \\
Fig.\ref{23}(h)    &   (21, 64, 195)      &   1.0254e-1    &  1.2190e-1   & 1      & 0.4209    \\
\hline
\end{tabular}
\end{table}

\begin{figure}[htbp]
\renewcommand{\figurename}{Fig.}
  \centering
  \subfigure[$J=4$ for MLASSO: $\lambda_{1}=\lambda_{2}=\lambda_{3}=\lambda_{4}=0.001$.]{
  \includegraphics[width=0.3\textwidth]{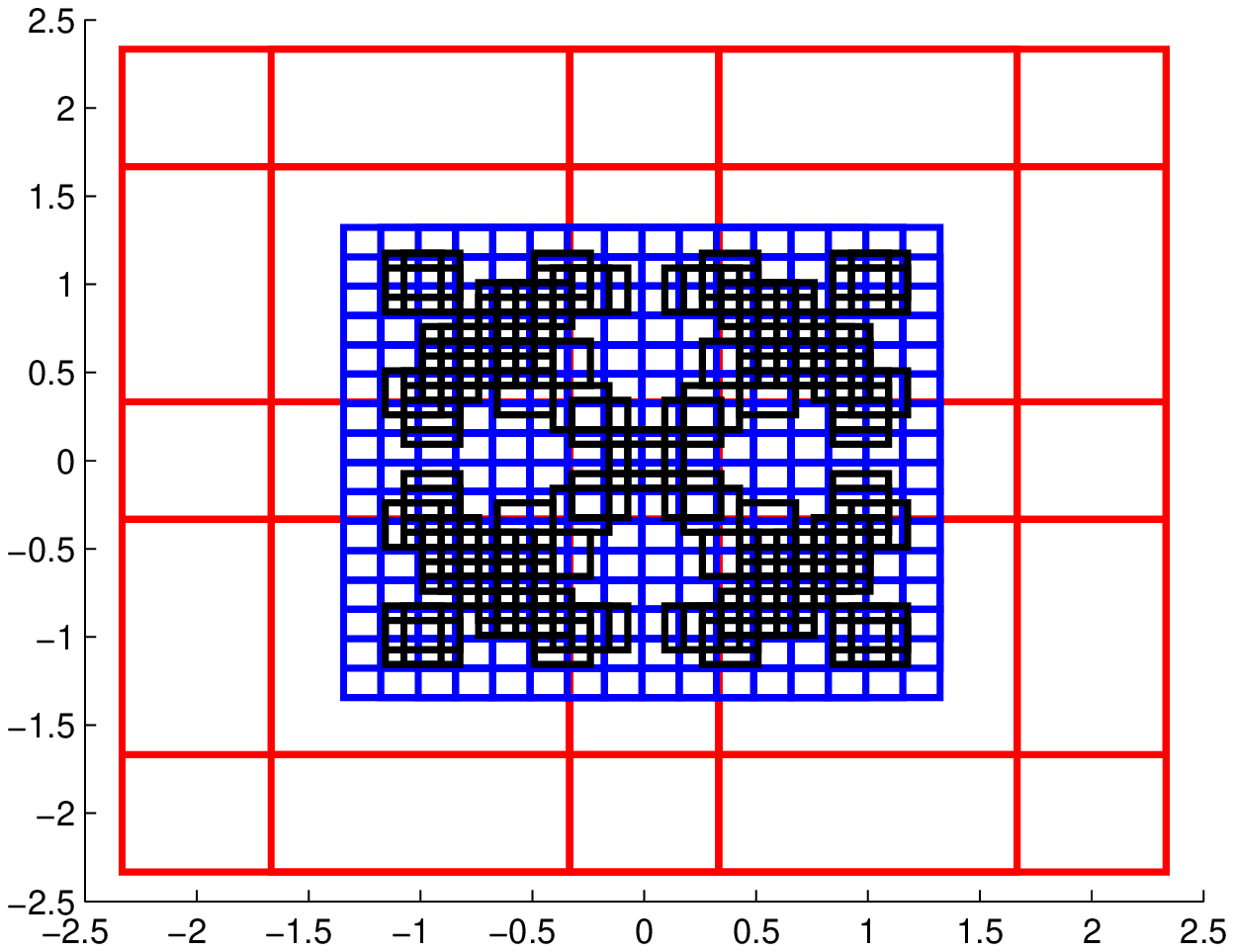}}
  \subfigure[$J=4$ for MLASSO: $\lambda_{1}=\lambda_{2}=\lambda_{3}=\lambda_{4}=0.01$.]{
  \includegraphics[width=0.3\textwidth]{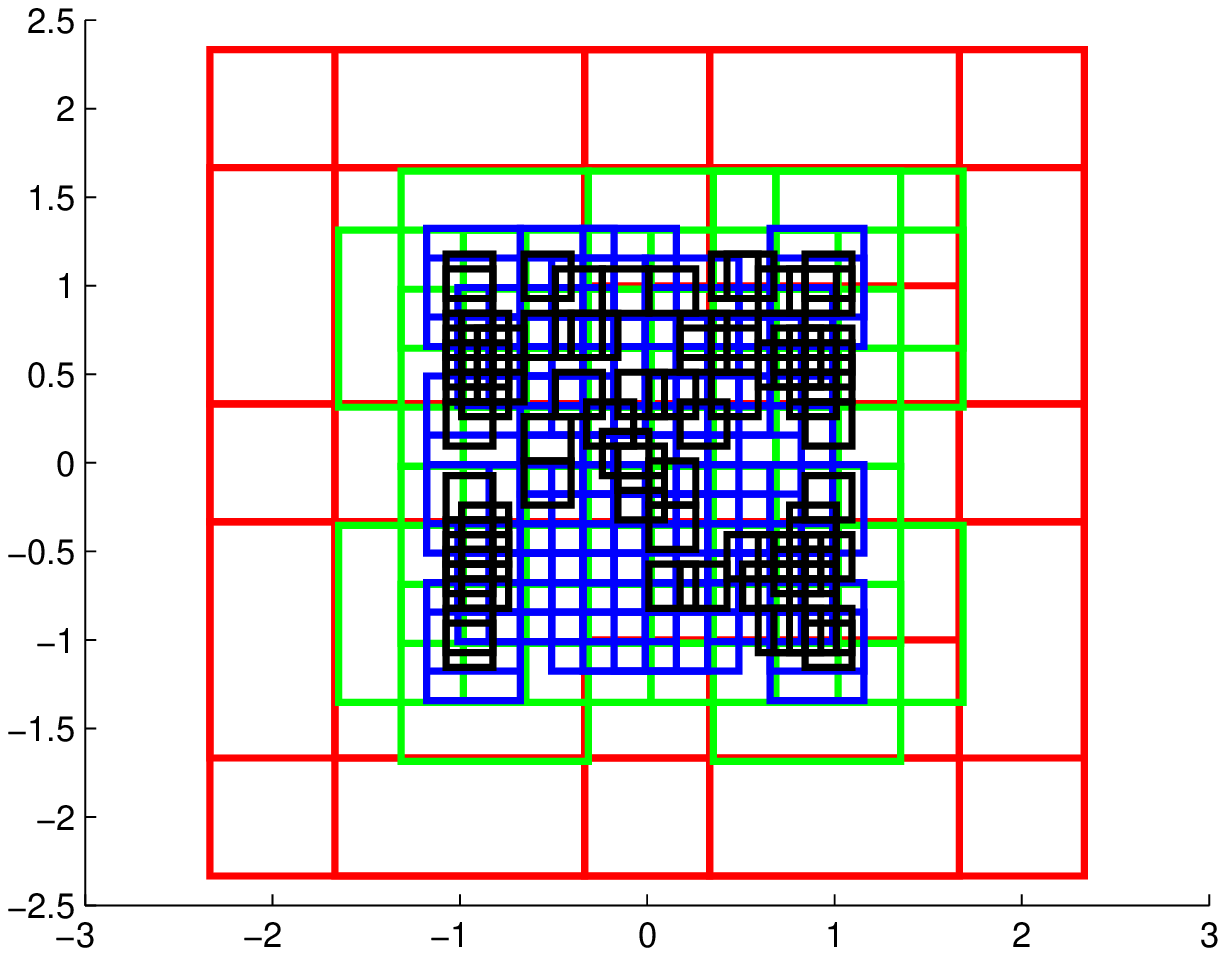}}
  \subfigure[$J=4$ for MLASSO: $\lambda_{1}=0.03,\lambda_{2}=0.02,\lambda_{3}=0.01,\lambda_{4}=0.04$.]{
  \includegraphics[width=0.3\textwidth]{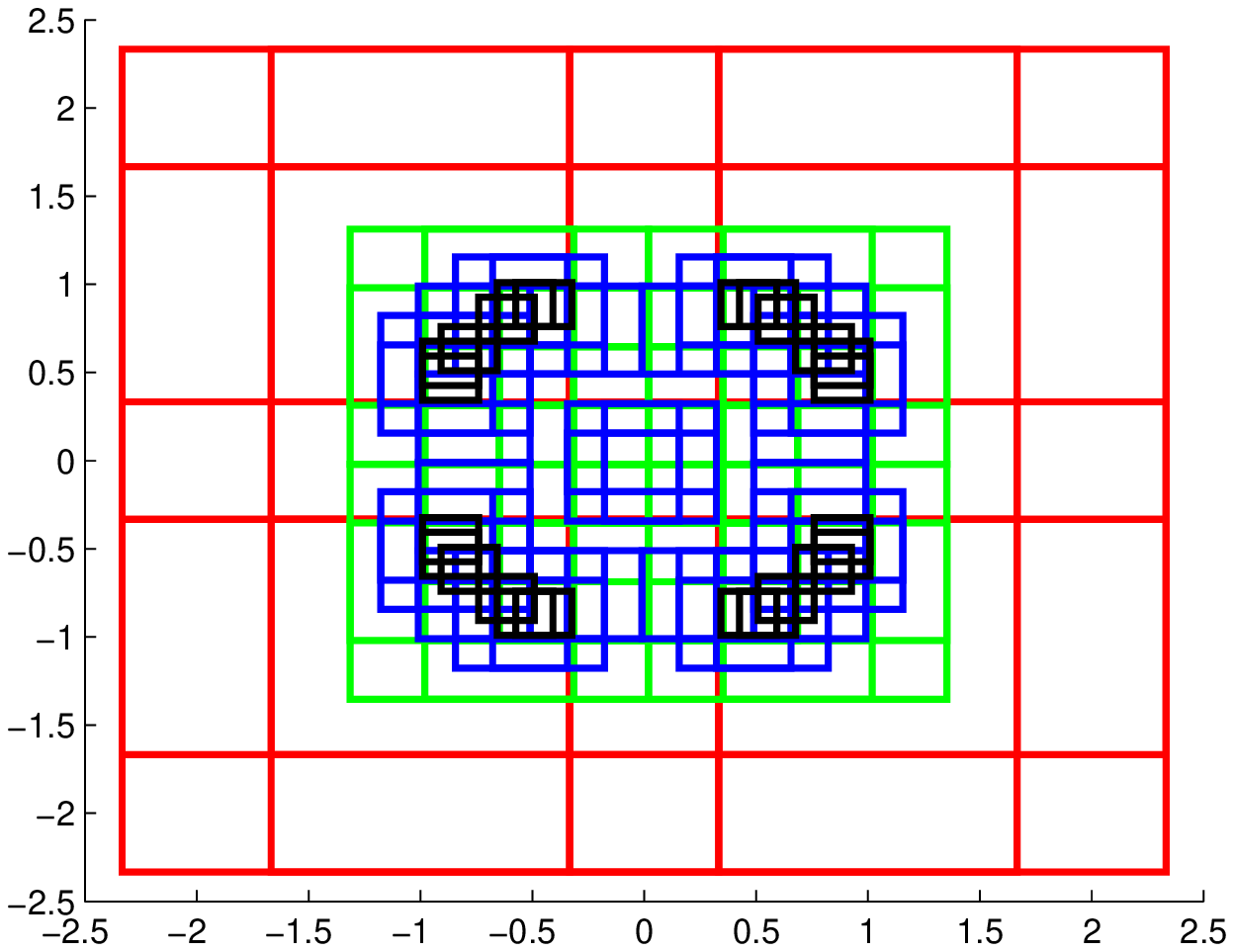}}
  \subfigure[$J=4$ for MLASSO: $\lambda_{1}=0.05,\lambda_{2}=0.01,\lambda_{3}=0.02,\lambda_{4}=0.05$.]{
  \includegraphics[width=0.3\textwidth]{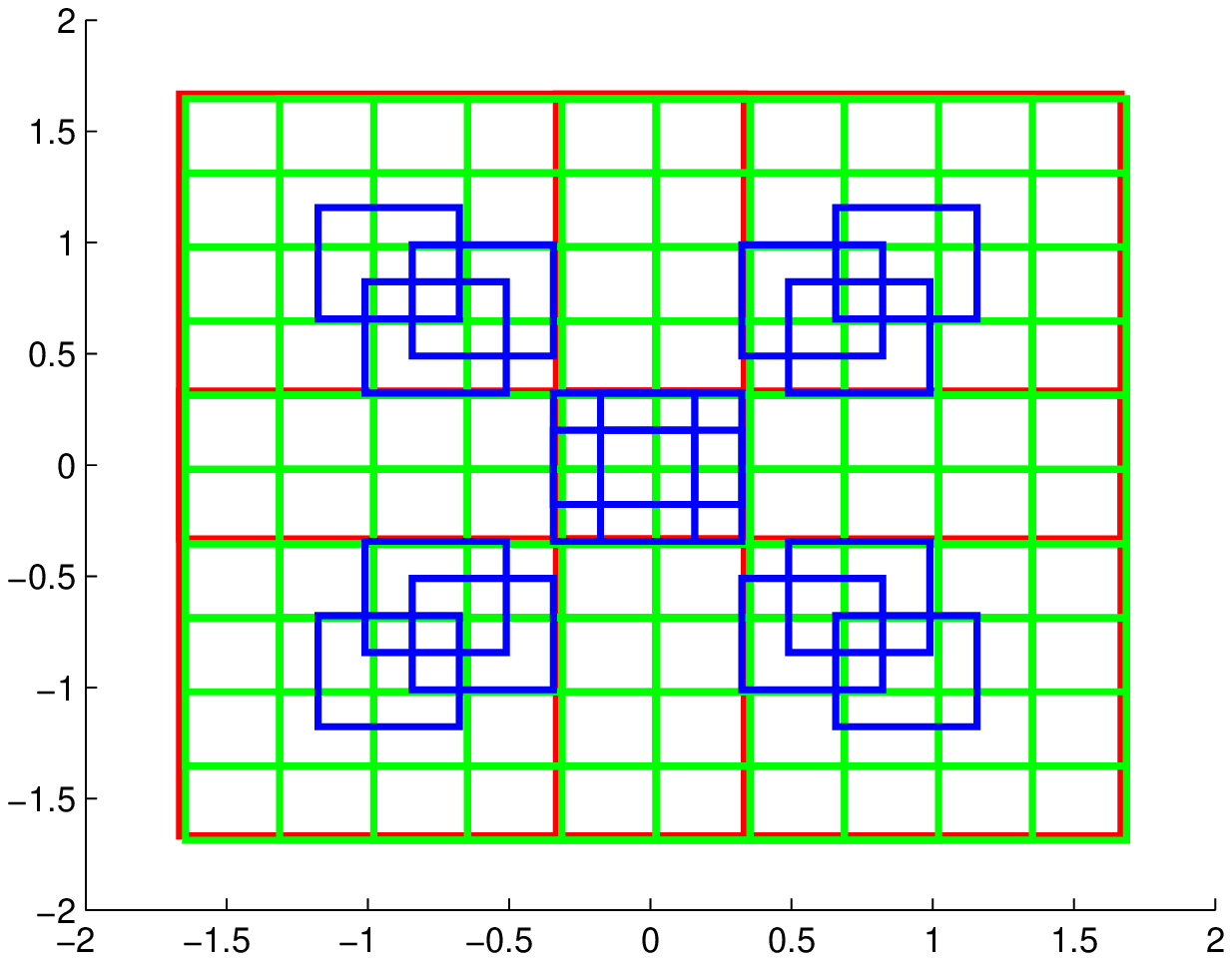}}
  \subfigure[$J=4$ for AGLASSO: $\lambda_{1}=\lambda_{2}=\lambda_{3}=\lambda_{4}=0.001$.]{
 \includegraphics[width=0.3\textwidth]{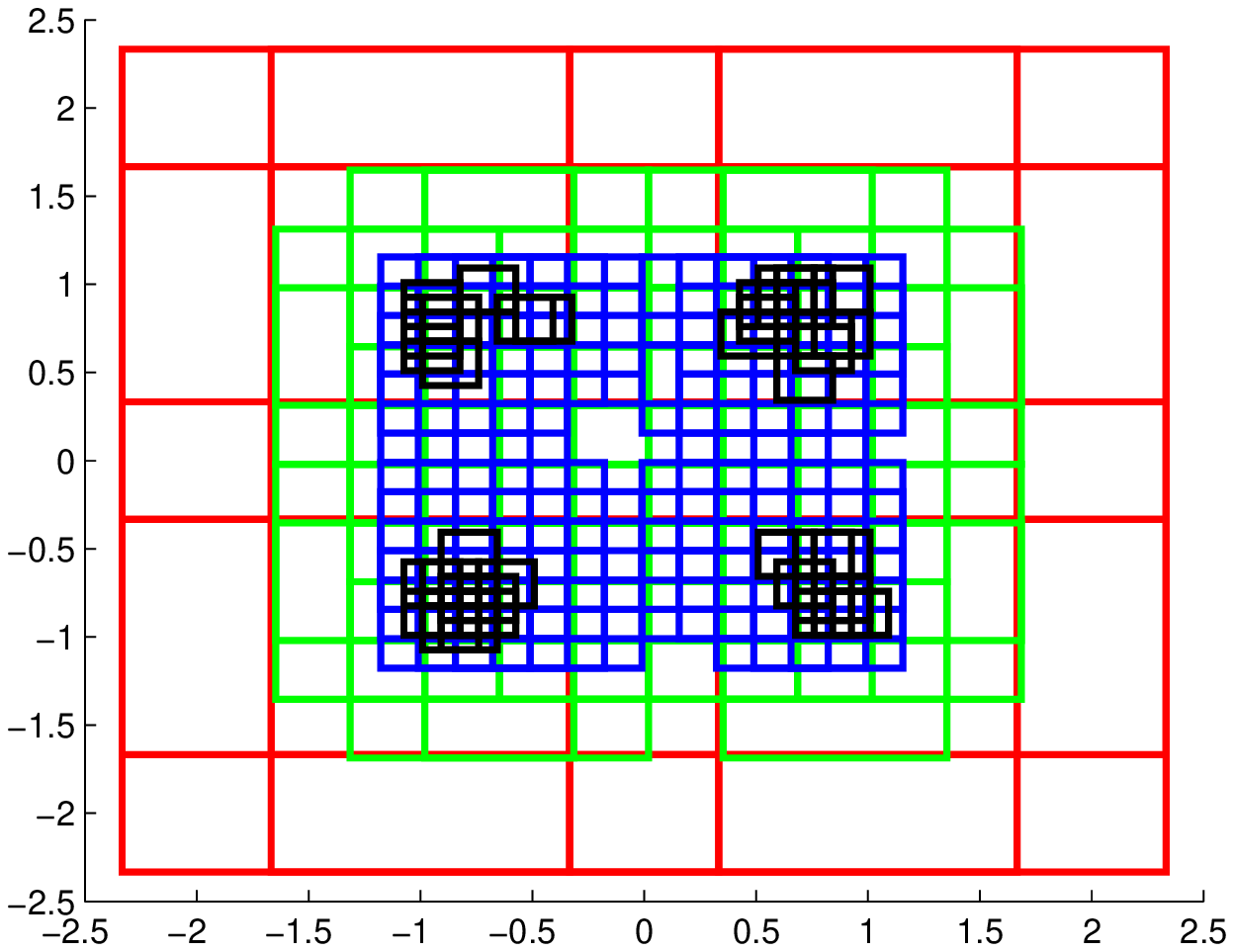}}
  \subfigure[$J=4$ for AGLASSO: $\lambda_{1}=\lambda_{2}=\lambda_{3}=\lambda_{4}=0.01$.]{
  \includegraphics[width=0.3\textwidth]{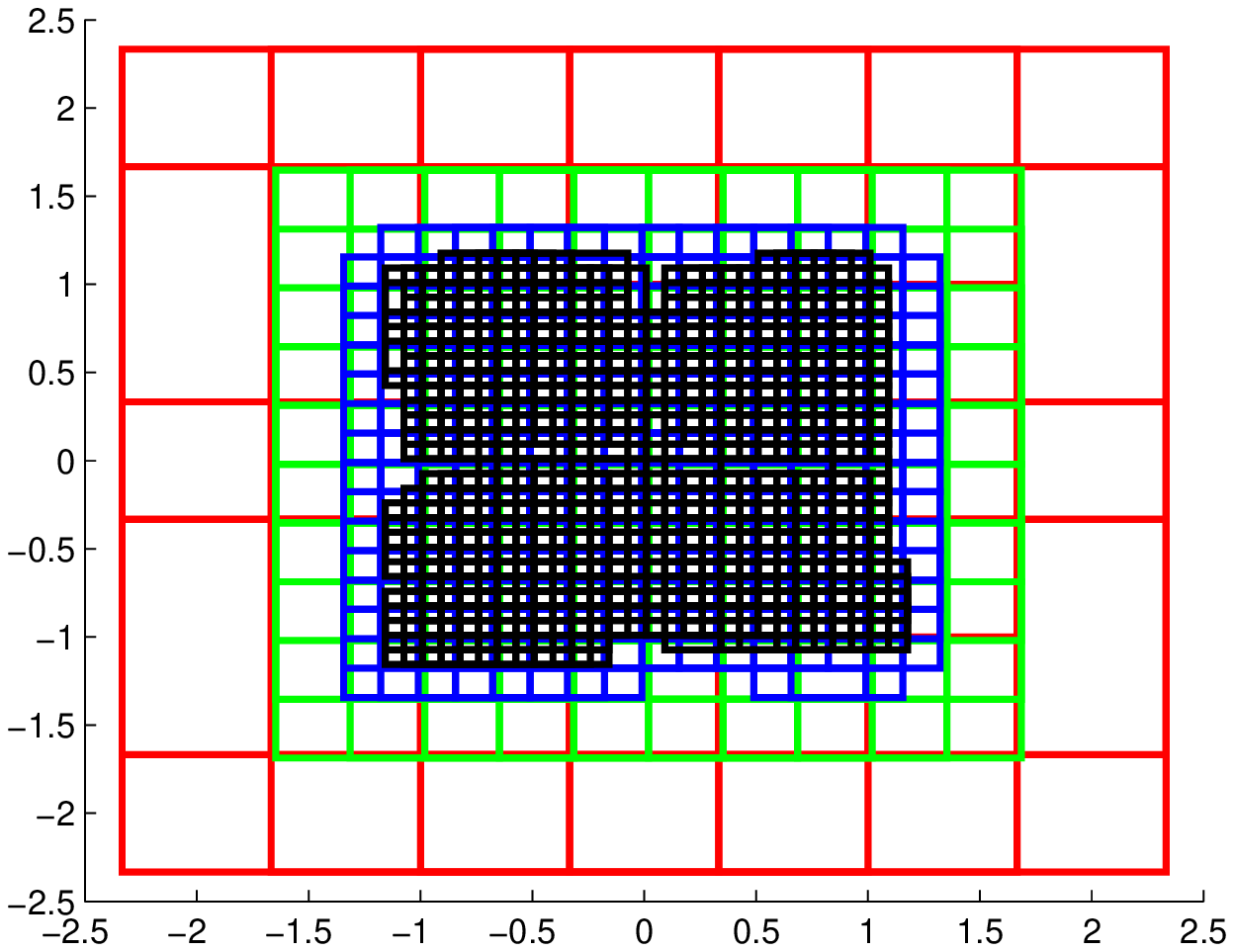}}
  \subfigure[$J=4$ for AGLASSO: $\lambda_{1}=0.03,\lambda_{2}=0.02,\lambda_{3}=0.01,\lambda_{4}=0.04$.]{
  \includegraphics[width=0.3\textwidth]{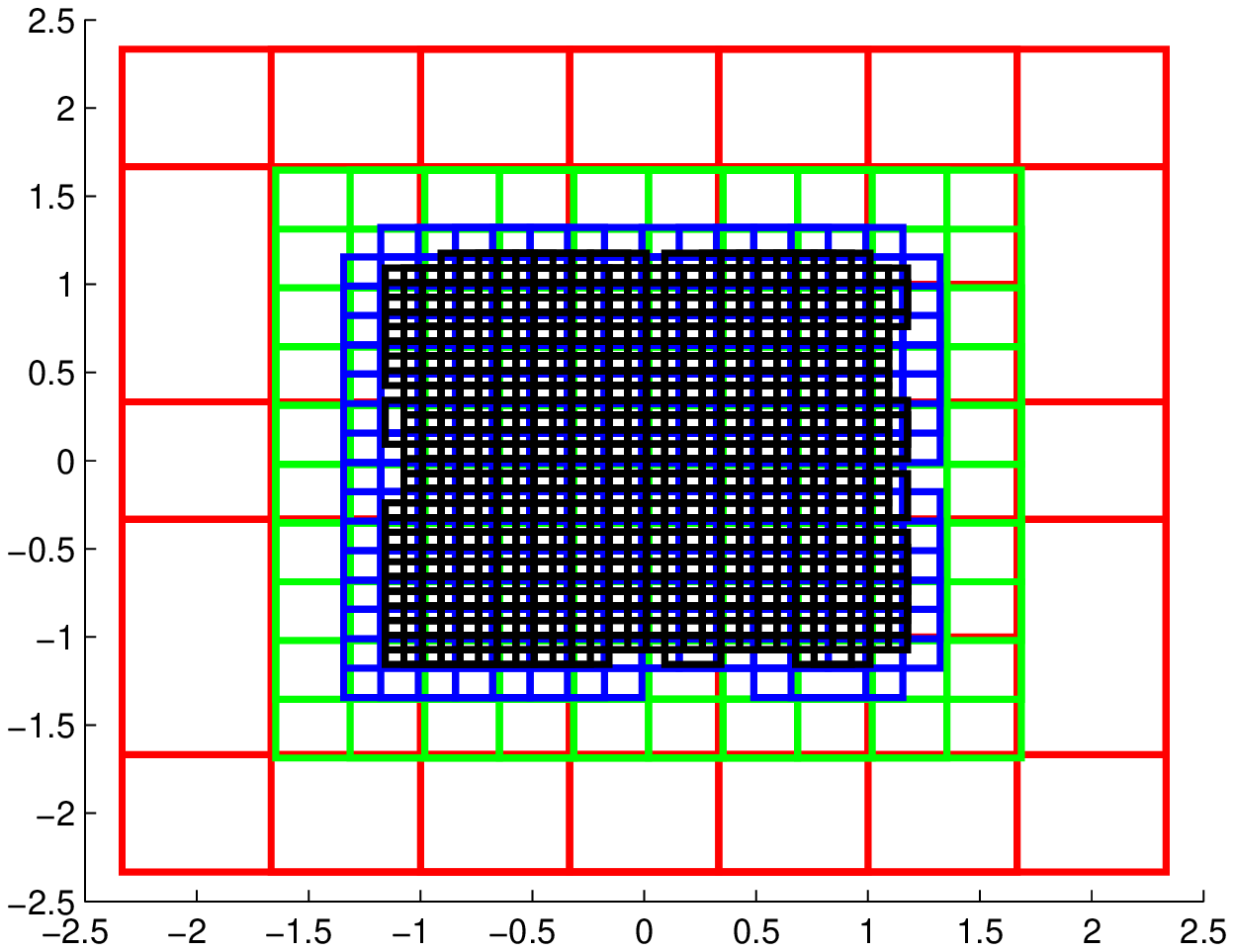}}
  \subfigure[$J=4$ for AGLASSO: $\lambda_{1}=0.05,\lambda_{2}=0.01,\lambda_{3}=0.02,\lambda_{4}=0.05$.]{
  \includegraphics[width=0.3\textwidth]{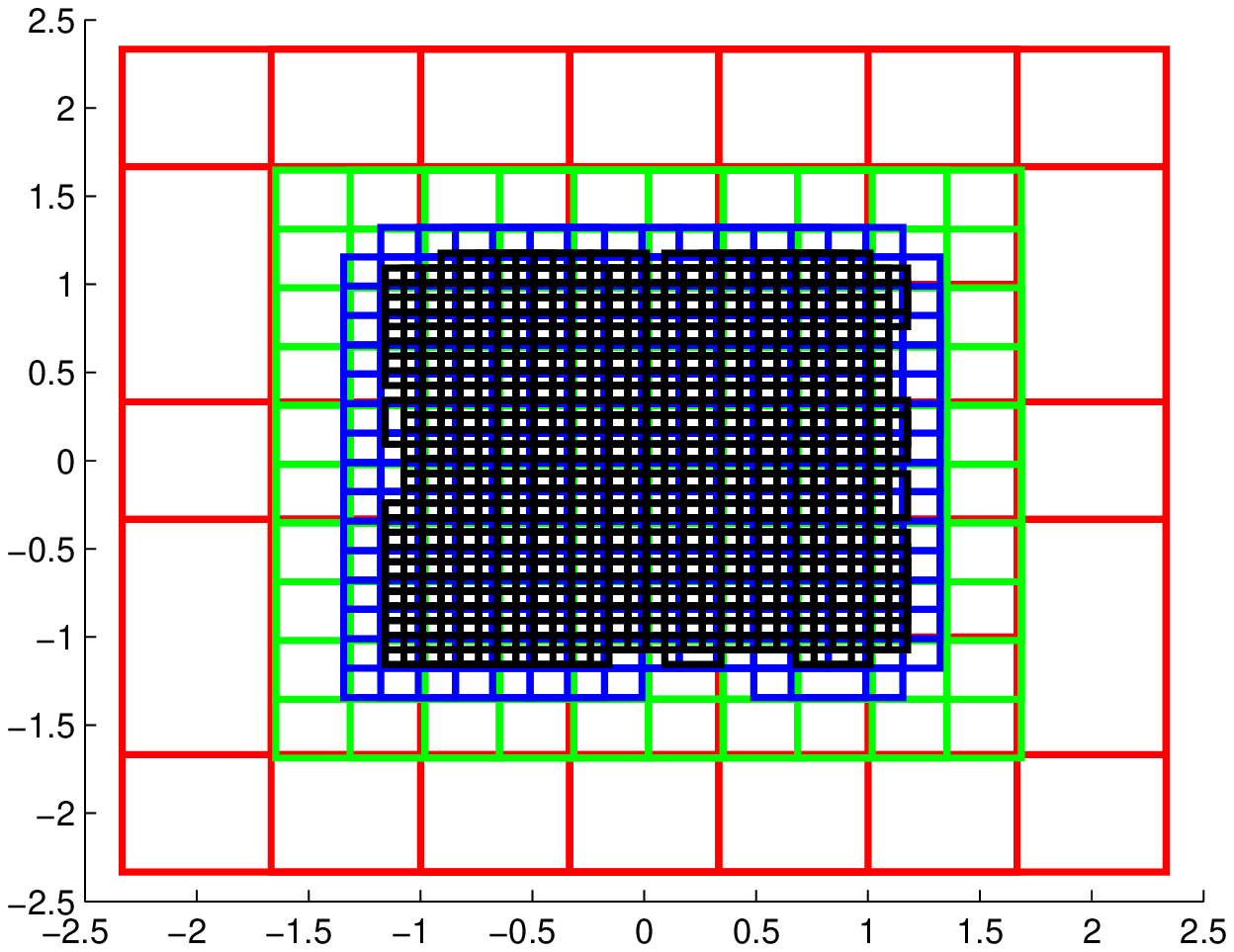}}
  \caption{The distribution of the support for $f_{2}$ with $J=4$ . }
  \label{24}
\end{figure}

\begin{table}[htbp]
\center
\caption{The $l_{0}$ norm and the approximation errors for different parameters shown in Fig.\ref{24} for $f_{2}$ with $J=4$.}\label{Tab:24}
\begin{tabular}{cccccc}\hline
                   & $(l_{0}(\textbf{X}_{1}),l_{0}(\textbf{X}_{2}),l_{0}(\textbf{X}_{3}),l_{0}(\textbf{X}_{4}))$  & Error & RMS & Iterations  & Time(sec) \\ \hline
\cite{Lee}         &   (22, 63, 195, 673)     &   2.1342e-3    &  3.2496e-3   & 4      & 0.7326   \\
Fig.\ref{24}(a)    &   (16, 28, 45, 73)       &   3.0054e-3    &  3.4124e-3   & 378    & 5.0312   \\
Fig.\ref{24}(b)    &   (17, 31, 77, 54)       &   2.3321e-3    &  3.0061e-3   & 1214   & 9.7821   \\
Fig.\ref{24}(c)    &   (15, 25, 56, 24)       &   4.1072e-3    &  5.1132e-3   & 2601   & 13.8871    \\
Fig.\ref{24}(d)    &   (4, 61, 18, 0)         &   5.7461e-3    &  6.3242e-3   & 3261   & 15.3487    \\
Fig.\ref{24}(e)    &   (16, 52, 91, 40)       &   1.5002e-1    &  2.1392e-1   & 1      & 5.1902    \\
Fig.\ref{24}(f)    &   (20, 64, 185, 471)     &   1.2301e-1    &  2.0342e-1   & 1      & 4.8106    \\
Fig.\ref{24}(g)    &   (21, 64, 180, 589)     &   2.0039e-1    &  2.1759e-1   & 1      & 4.7951    \\
Fig.\ref{24}(h)    &   (20, 64, 185, 596)     &   3.1402e-1    &  3.1132e-1   & 1      & 4.6283    \\
\hline
\end{tabular}
\end{table}

\begin{figure}[htbp]
\renewcommand{\figurename}{Fig.}
  \centering
  \subfigure[$J=3$ for MLASSO: $\lambda_{1}=\lambda_{2}=\lambda_{3}=0.001$.]{
  \includegraphics[width=0.3\textwidth]{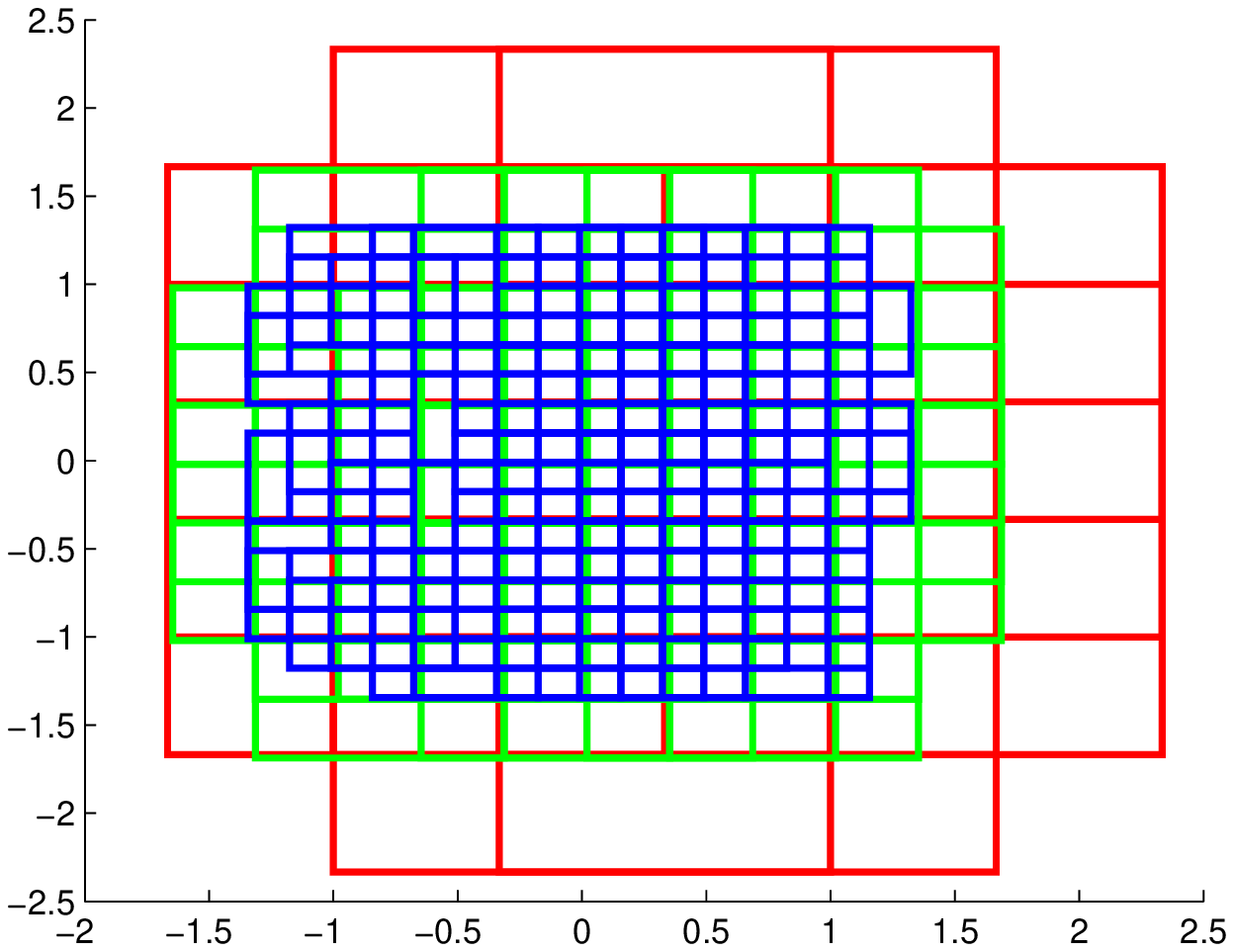}}
  \subfigure[$J=3$ for MLASSO: $\lambda_{1}=\lambda_{2}=\lambda_{3}=0.01$.]{
  \includegraphics[width=0.3\textwidth]{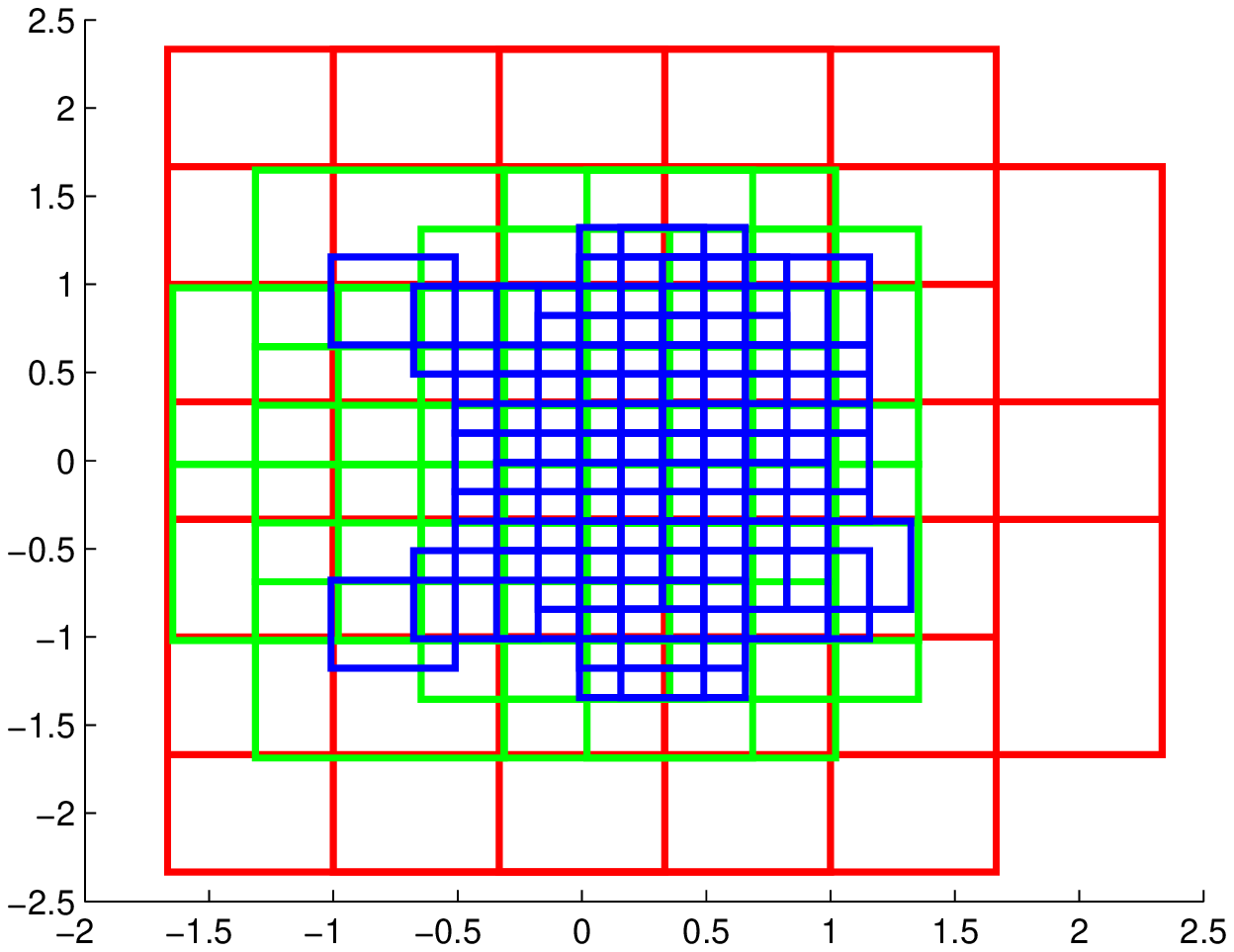}}
  \subfigure[$J=3$ for MLASSO: $\lambda_{1}=0.02,\lambda_{2}=0.01,\lambda_{3}=0.03$.]{
  \includegraphics[width=0.3\textwidth]{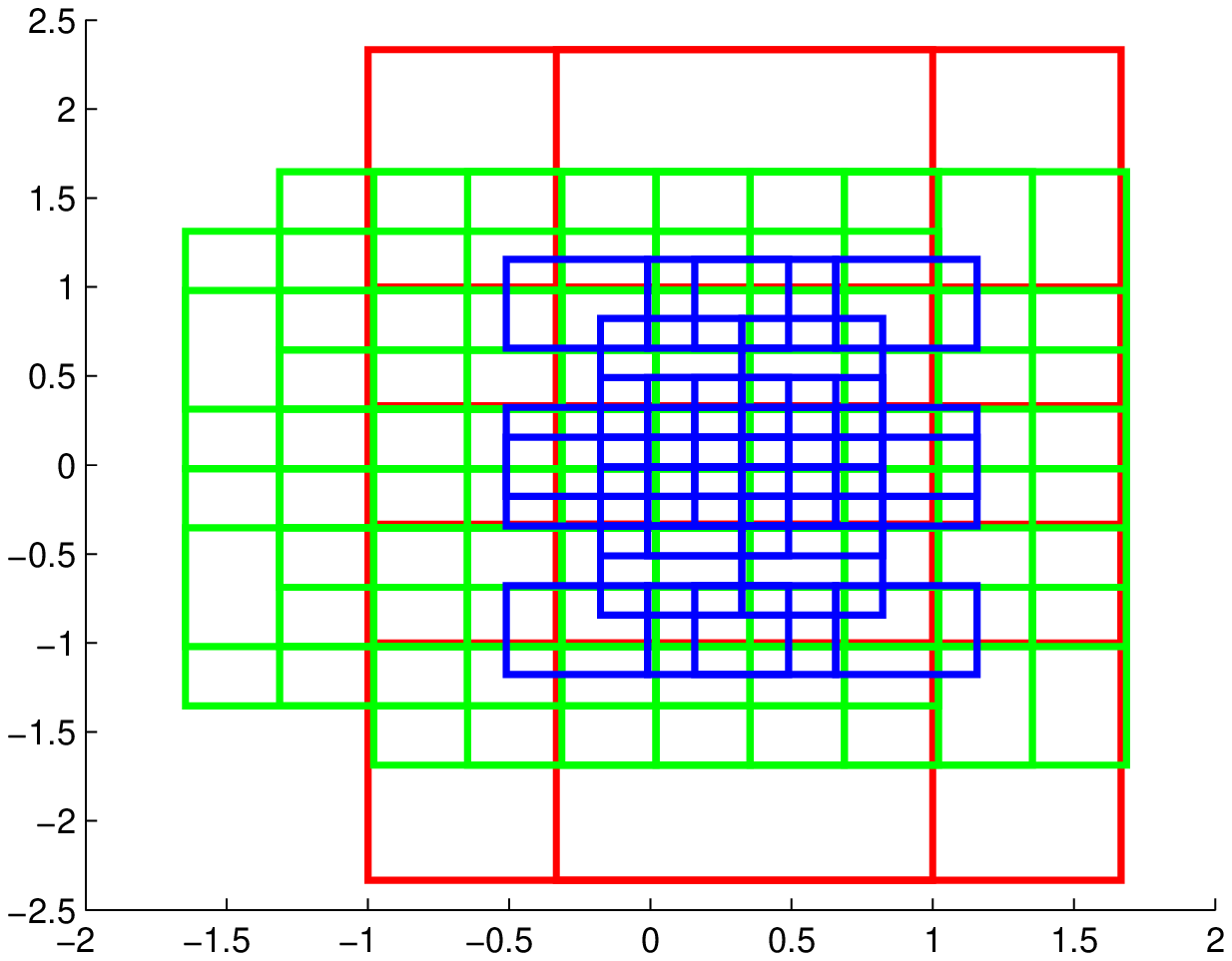}}
  \subfigure[$J=3$ for MLASSO: $\lambda_{1}=0.03, \lambda_{2}=0.02,\lambda_{3}=0.05$.]{
  \includegraphics[width=0.3\textwidth]{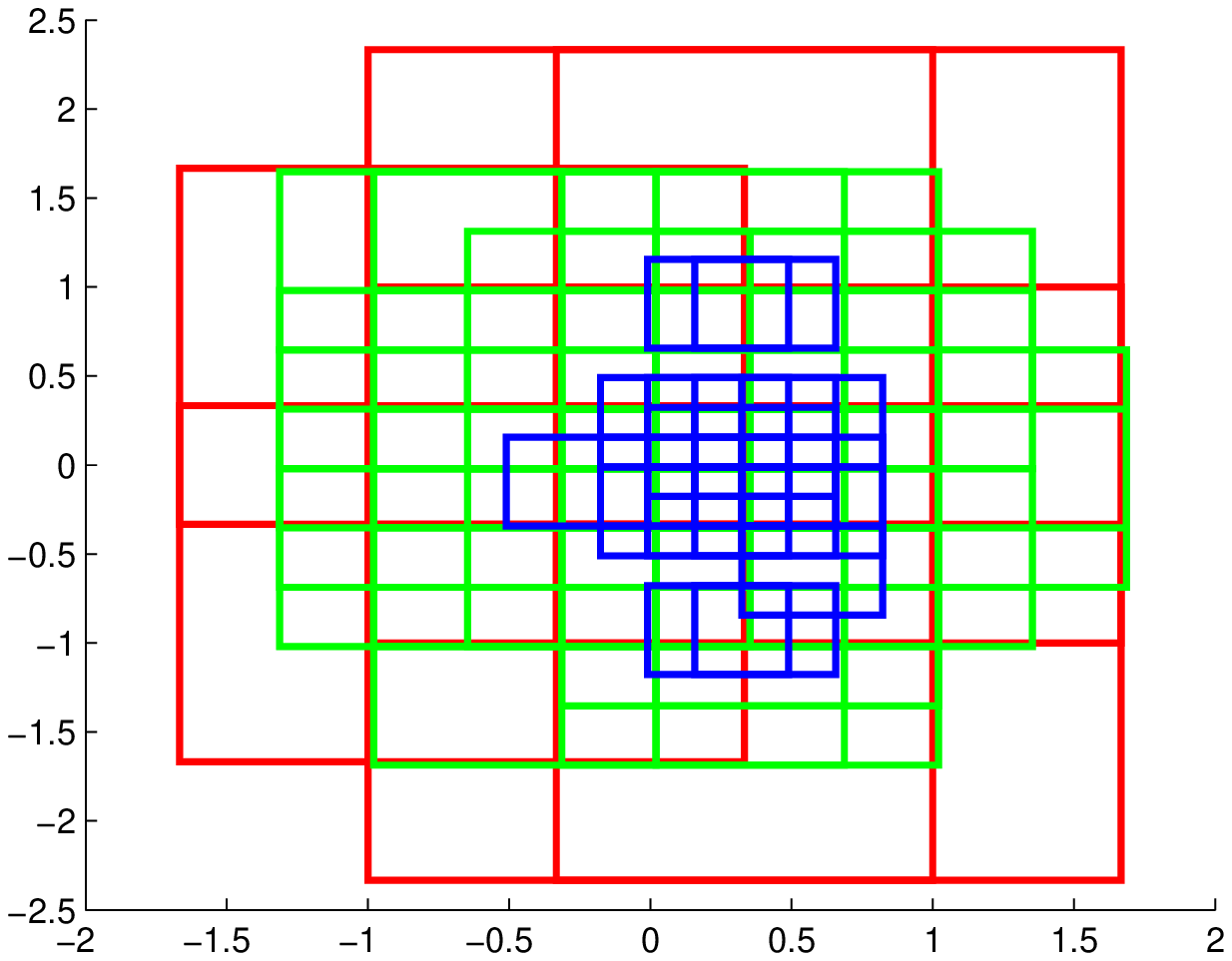}}
  \subfigure[$J=3$ for AGLASSO: $\lambda_{1}=\lambda_{2}=\lambda_{3}=0.001$.]{
  \includegraphics[width=0.3\textwidth]{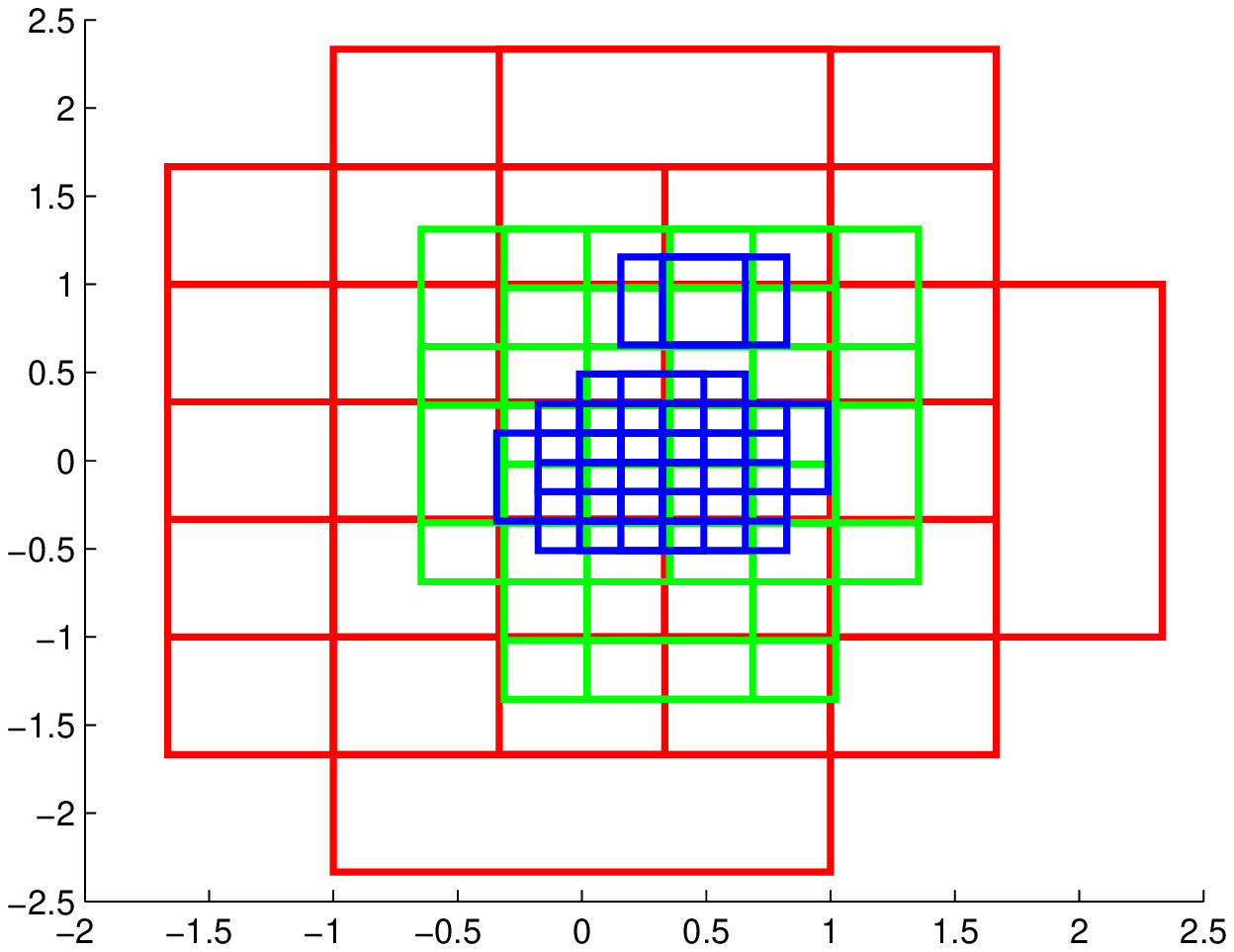}}
  \subfigure[$J=3$ for AGLASSO: $\lambda_{1}=\lambda_{2}=\lambda_{3}=0.01$.]{
  \includegraphics[width=0.3\textwidth]{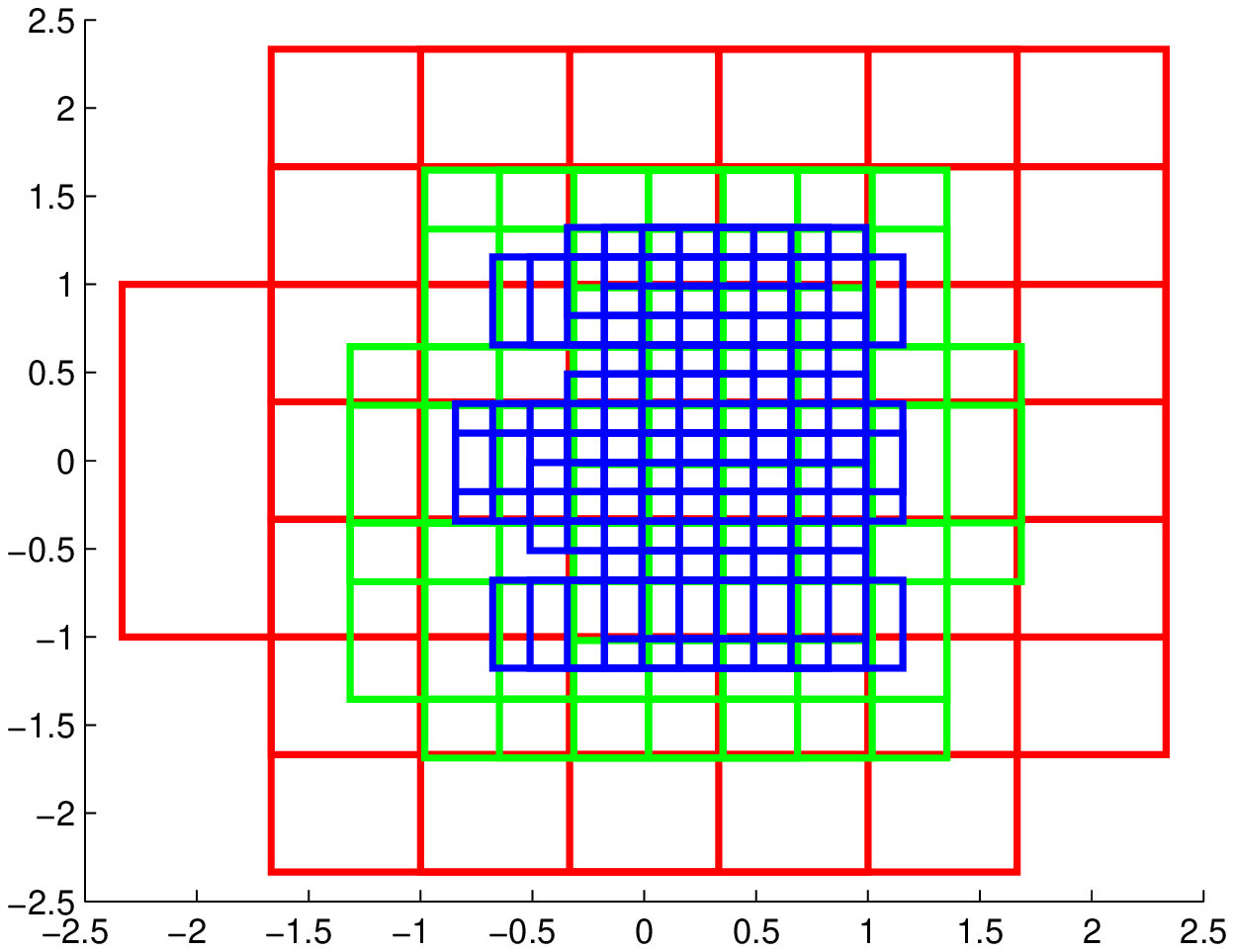}}
  \subfigure[$J=3$ for AGLASSO: $\lambda_{1}=0.02,\lambda_{2}=0.01,\lambda_{3}=0.03$.]{
  \includegraphics[width=0.3\textwidth]{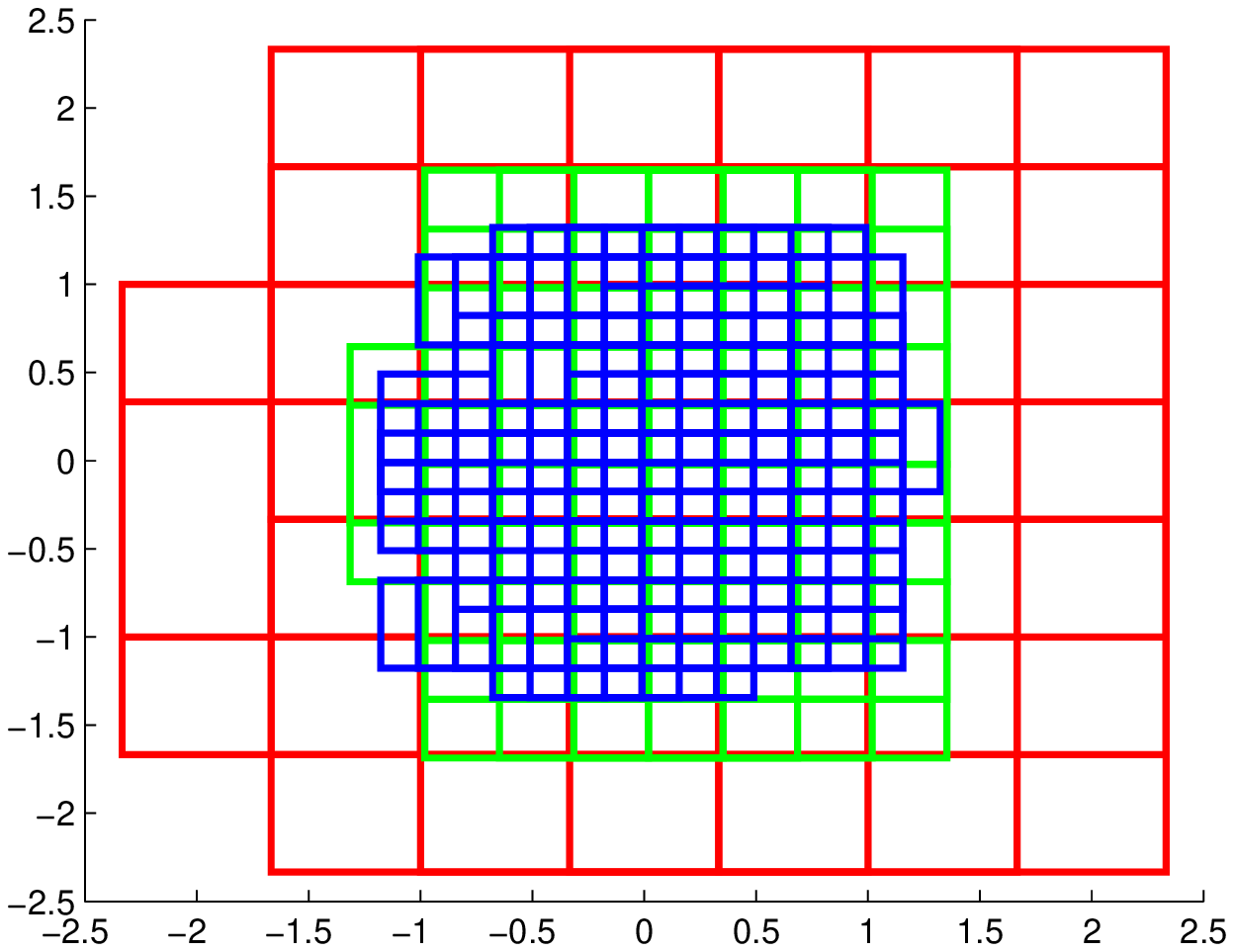}}
  \subfigure[$J=3$ for AGLASSO: $\lambda_{1}=0.03, \lambda_{2}=0.02,\lambda_{3}=0.05$.]{
  \includegraphics[width=0.3\textwidth]{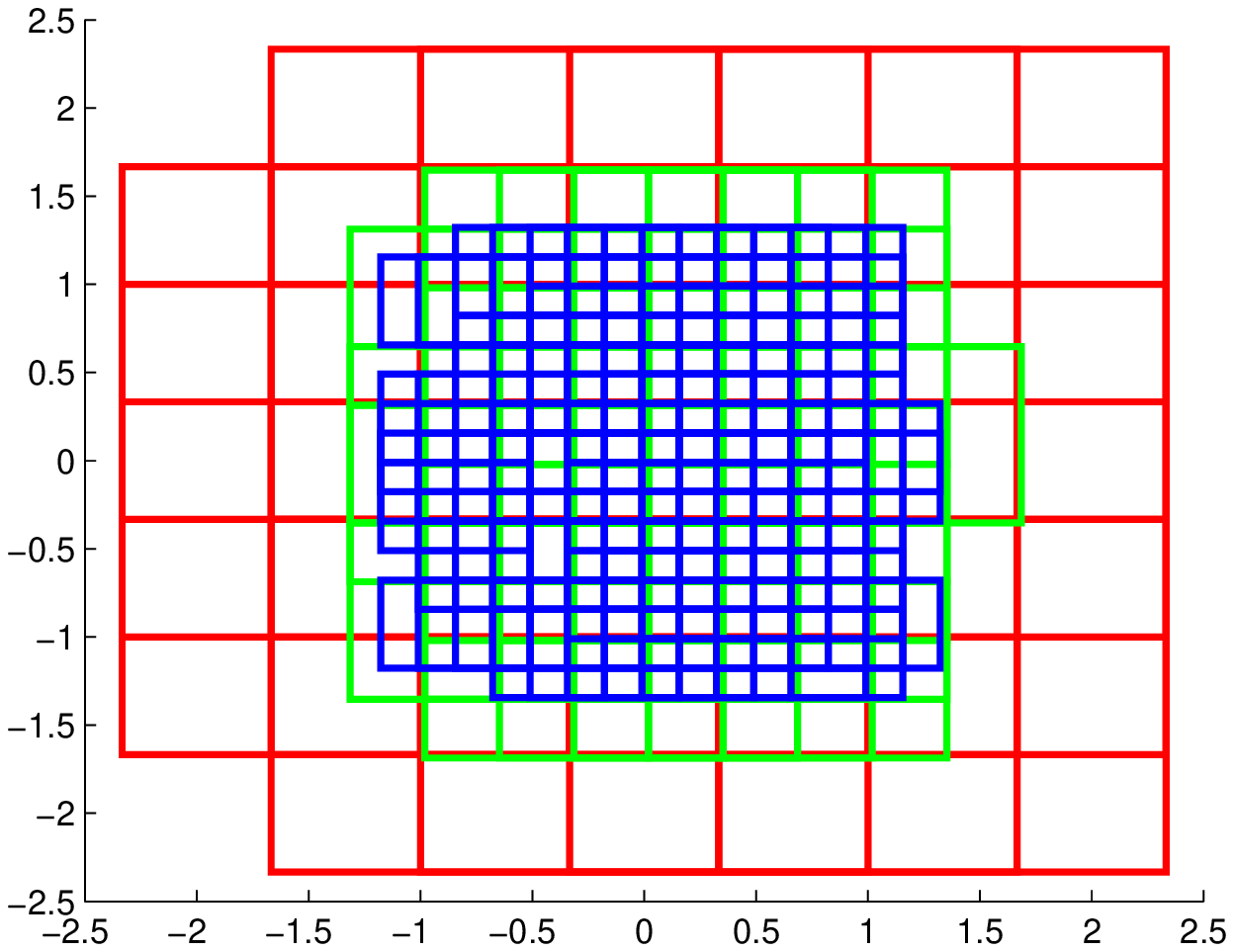}}
  \caption{The distribution of the support for $f_{3}$ with $J=3$.}
  \label{33}
\end{figure}

\begin{table}[htbp]
\center
\caption{The $l_{0}$ norm and the approximation errors for different parameters shown in Fig.\ref{33} for $f_{3}$ with $J=3$.}\label{Tab:33}
\begin{tabular}{cccccc}\hline
                   & $(l_{0}(\textbf{X}_{1}),l_{0}(\textbf{X}_{2}),l_{0}(\textbf{X}_{3}))$  & Error & RMS & Iterations  & Time(sec) \\ \hline
\cite{Lee}         &   (24, 63, 193)      &   3.0285e-3    &  5.8147e-3   & 3      & 0.0617   \\
Fig.\ref{33}(a)    &   (15, 50, 116)      &   3.0273e-3    &  6.3044e-3   & 298    & 0.5042   \\
Fig.\ref{33}(b)    &   (13, 33, 60)       &   4.1409e-3    &  7.6258e-3   & 739    & 1.0354    \\
Fig.\ref{33}(c)    &   (6, 50, 31)        &   4.8321e-3    &  7.6578e-3   & 978    & 1.3277    \\
Fig.\ref{33}(d)    &   (7, 34, 20)        &   5.1157e-3    &  9.2785e-3   & 683    & 1.0317    \\
Fig.\ref{33}(e)    &   (13, 19, 18)       &   4.6342e-2    &  5.1462e-2   & 1      & 0.5745    \\
Fig.\ref{33}(f)    &   (21, 39, 74)       &   5.7749e-2    &  6.0324e-2   & 1      & 0.4182    \\
Fig.\ref{33}(g)    &   (22, 42, 123)      &   6.7765e-2    &  7.0079e-2   & 1      & 0.6242    \\
Fig.\ref{33}(h)    &   (24, 45, 138)      &   1.1091e-1    &  1.0072e-1   & 1      & 0.1562    \\
\hline
\end{tabular}
\end{table}

\begin{figure}[htbp]
\renewcommand{\figurename}{Fig.}
  \centering
  \subfigure[$J=4$ for MLASSO: $\lambda_{1}=\lambda_{2}=\lambda_{3}=\lambda_{4}=0.001$.]{
  \includegraphics[width=0.3\textwidth]{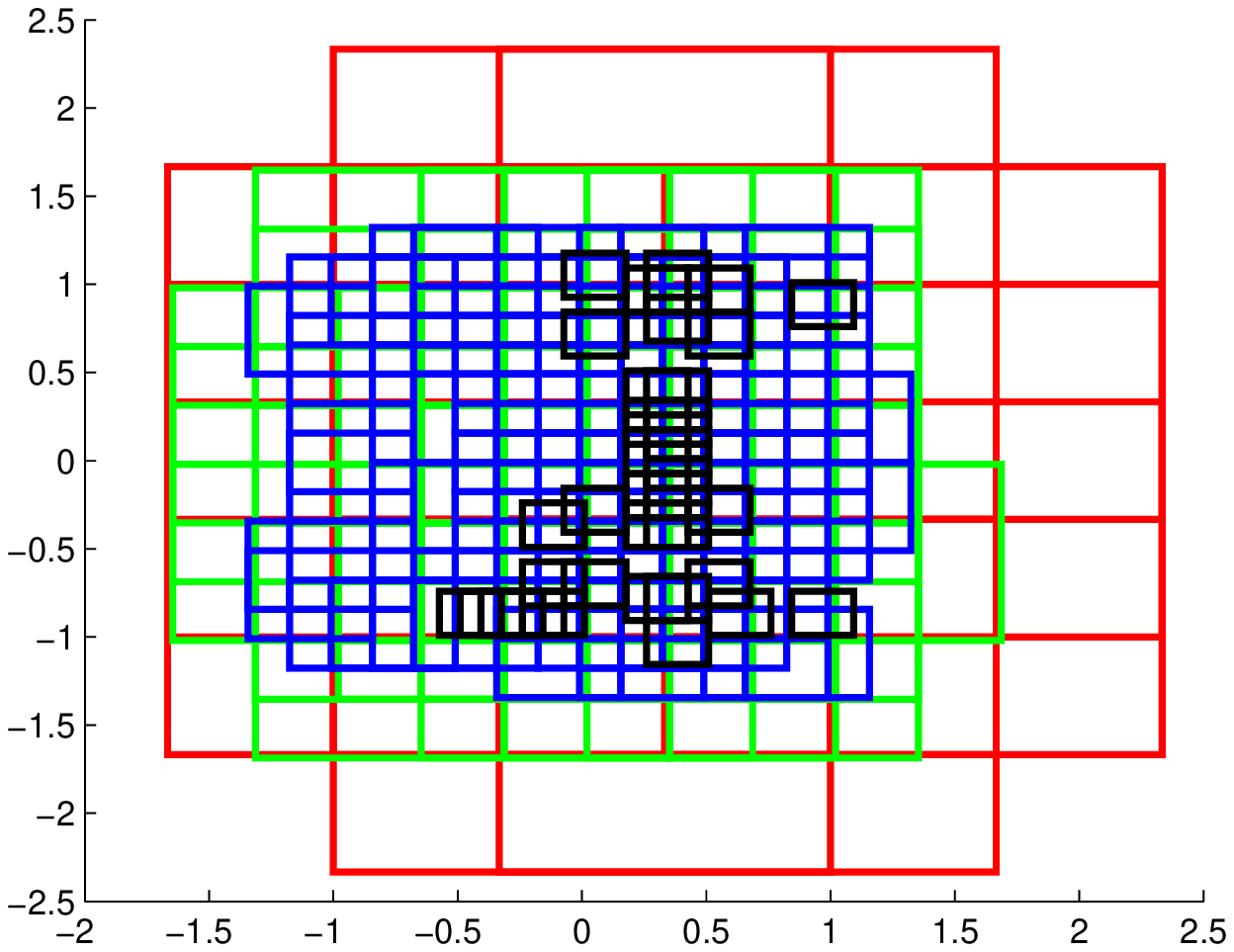}}
  \subfigure[$J=4$ for MLASSO: $\lambda_{1}=\lambda_{2}=\lambda_{3}=\lambda_{4}=0.01$.]{
  \includegraphics[width=0.3\textwidth]{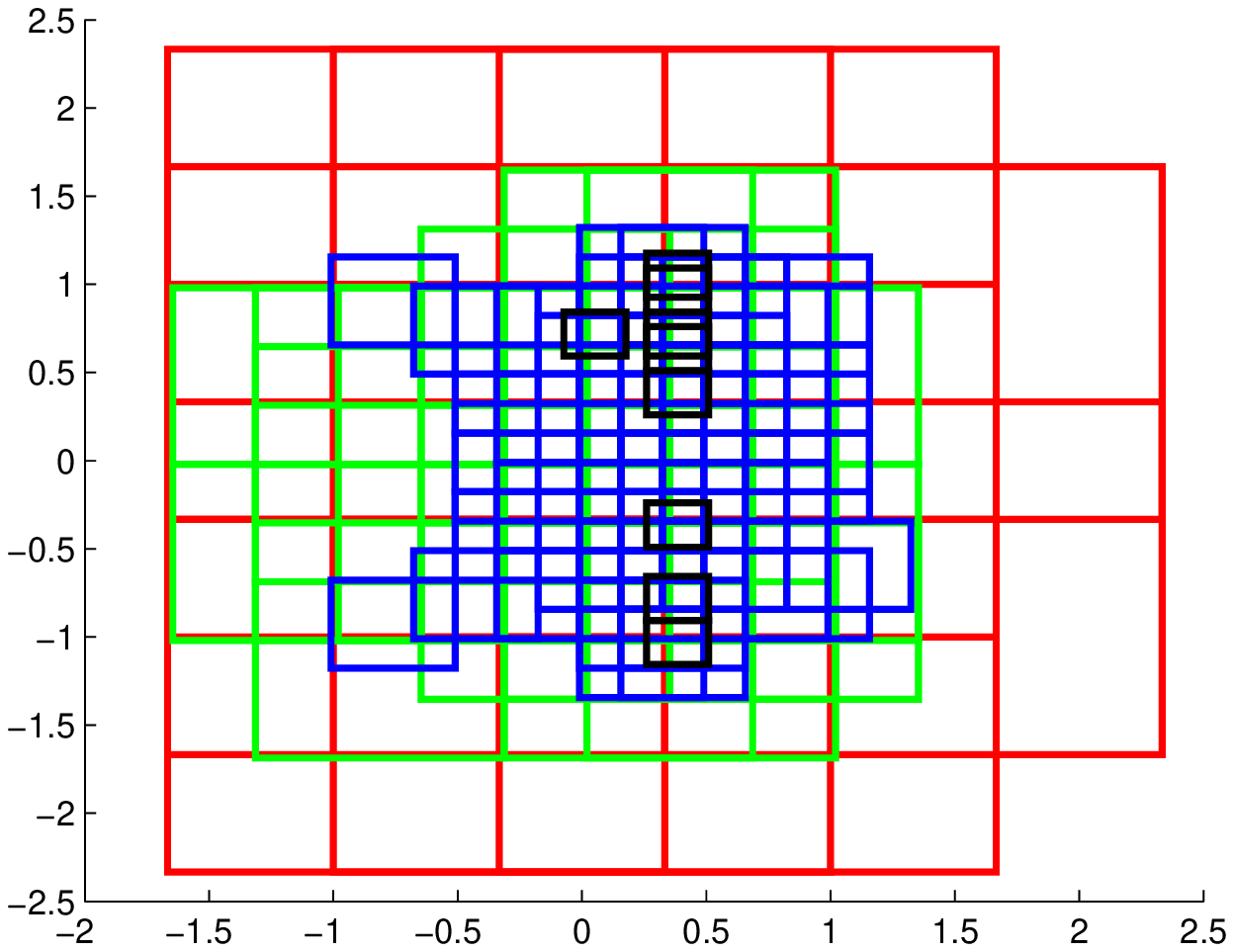}}
  \subfigure[$J=4$ for MLASSO: $\lambda_{1}=0.03,\lambda_{2}=0.02,\lambda_{3}=0.01,\lambda_{4}=0.04$.]{
  \includegraphics[width=0.3\textwidth]{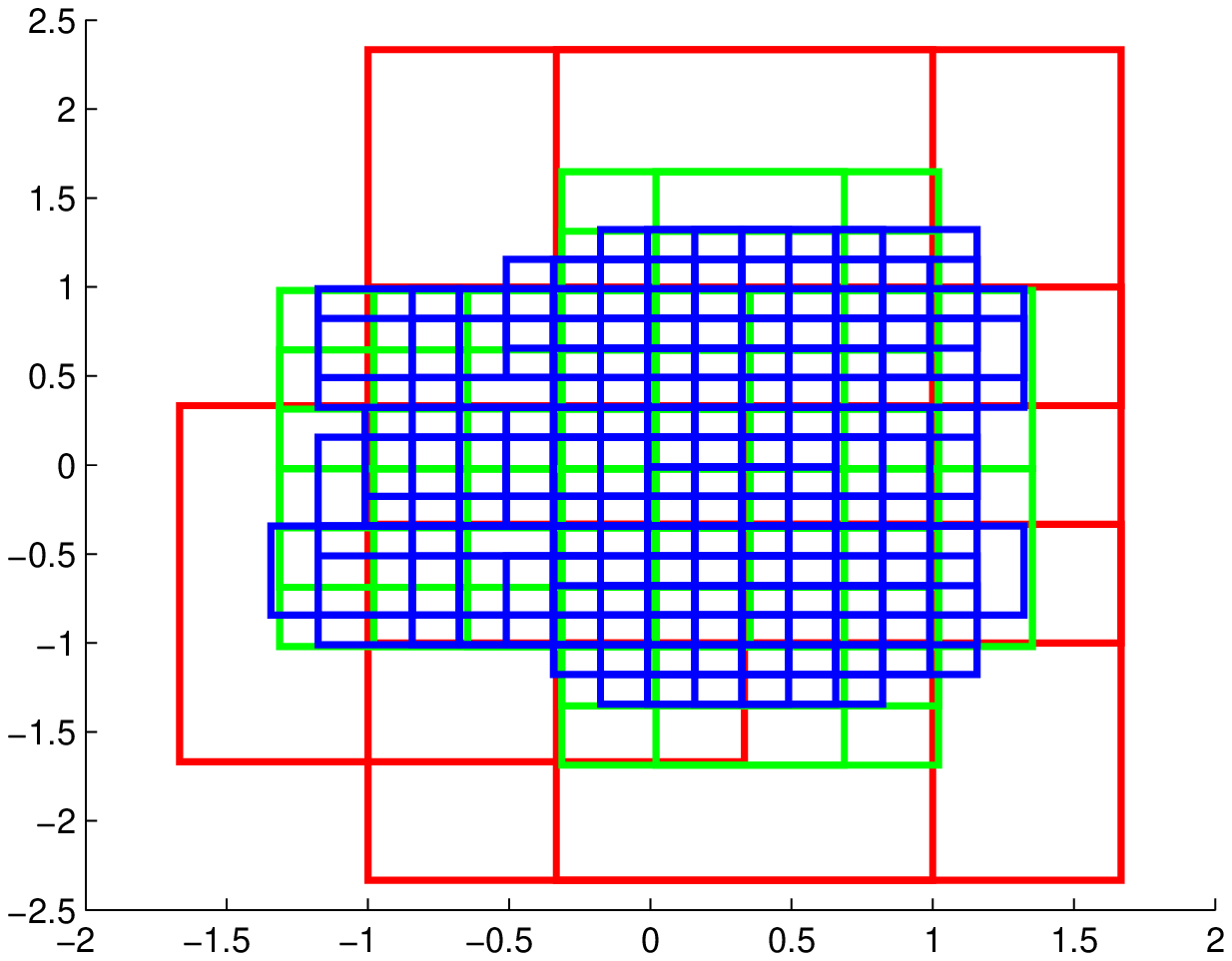}}
  \subfigure[$J=4$ for MLASSO: $\lambda_{1}=0.05,\lambda_{2}=0.01,\lambda_{3}=0.02,\lambda_{4}=0.05$.]{
  \includegraphics[width=0.3\textwidth]{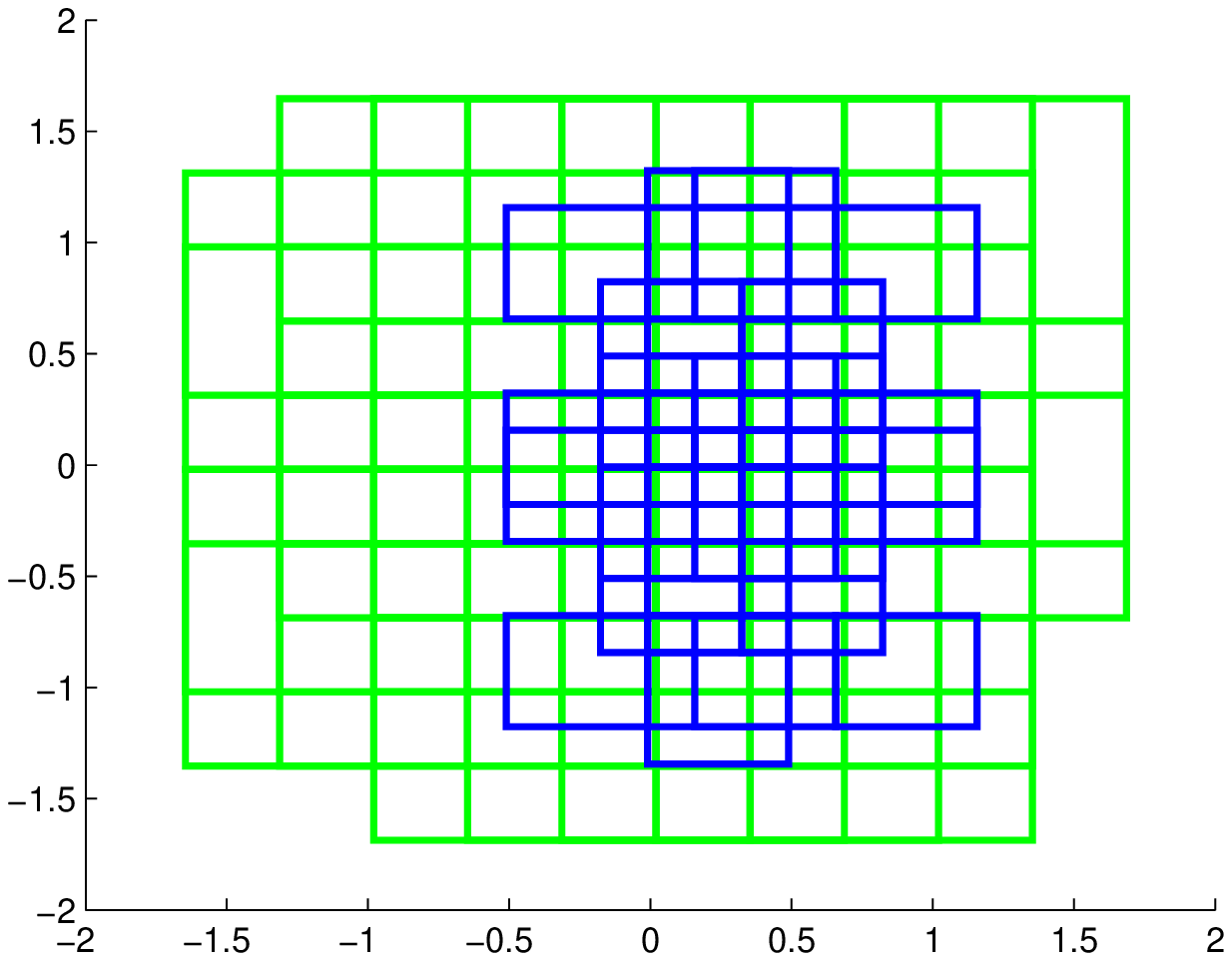}}
  \subfigure[$J=4$ for AGLASSO: $\lambda_{1}=\lambda_{2}=\lambda_{3}=\lambda_{4}=0.001$.]{
  \includegraphics[width=0.3\textwidth]{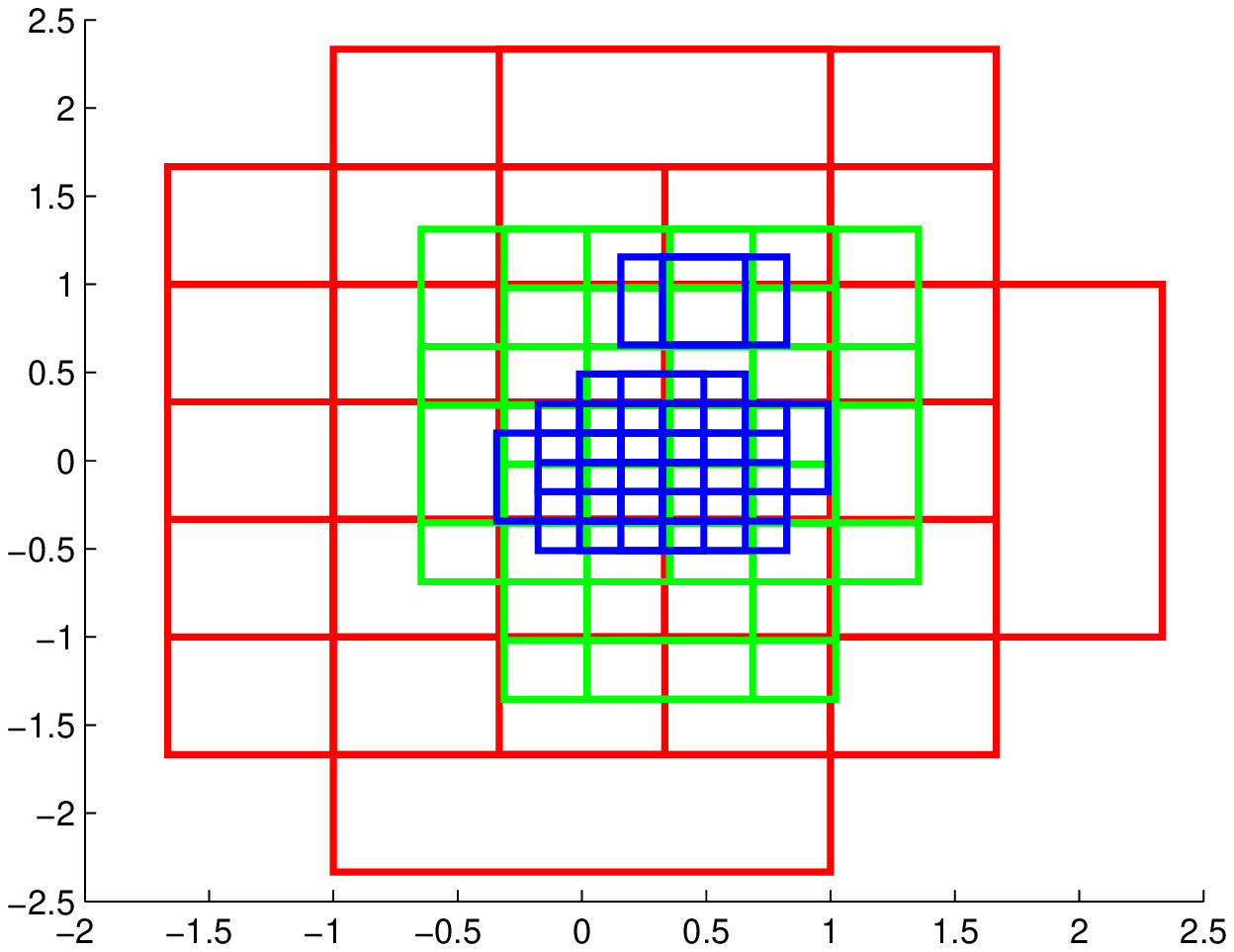}}
  \subfigure[$J=4$ for AGLASSO: $\lambda_{1}=\lambda_{2}=\lambda_{3}=\lambda_{4}=0.01$.]{
  \includegraphics[width=0.3\textwidth]{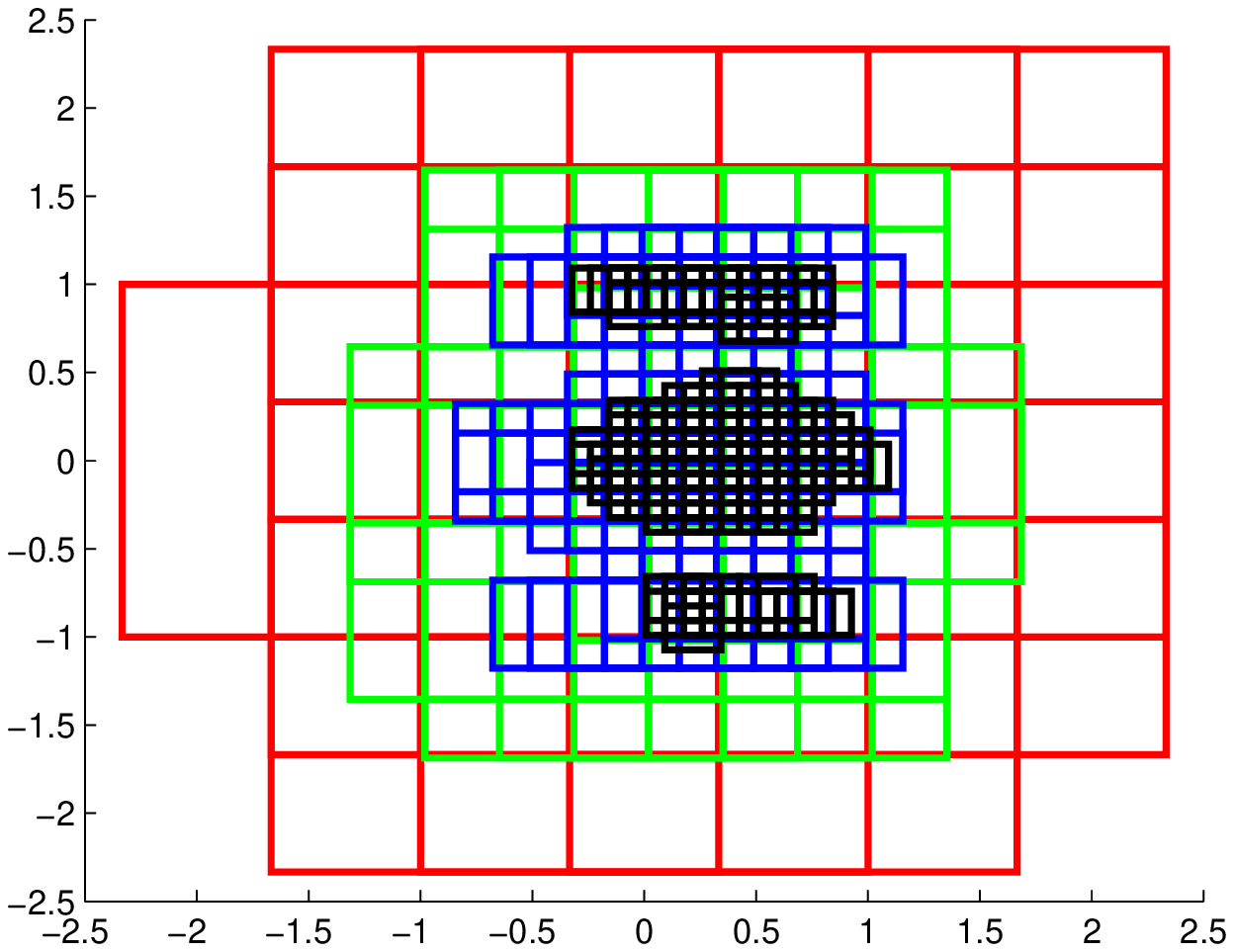}}
  \subfigure[$J=4$ for AGLASSO: $\lambda_{1}=0.03,\lambda_{2}=0.02,\lambda_{3}=0.01,\lambda_{4}=0.04$.]{
  \includegraphics[width=0.3\textwidth]{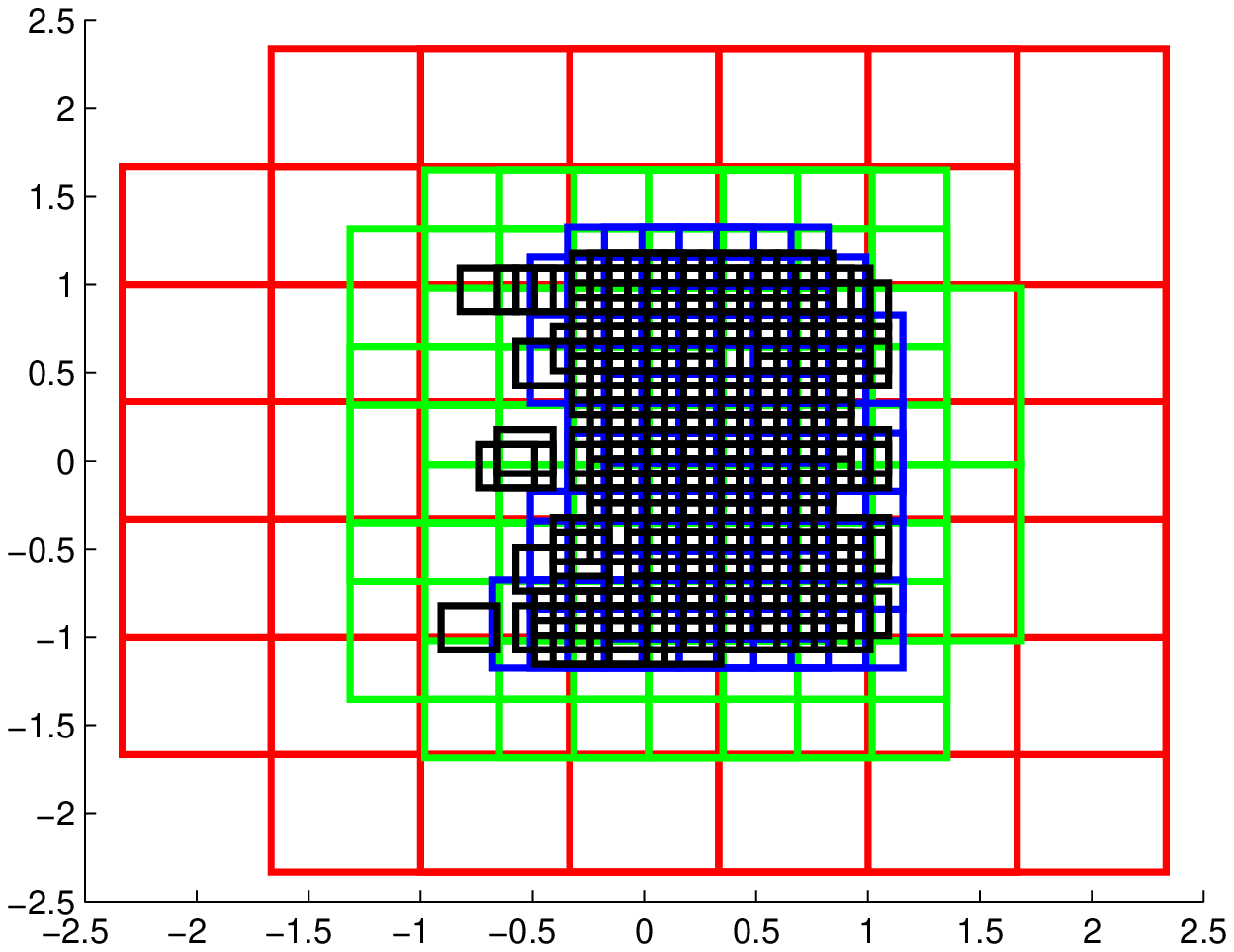}}
  \subfigure[$J=4$ for AGLASSO: $\lambda_{1}=0.05,\lambda_{2}=0.01,\lambda_{3}=0.02,\lambda_{4}=0.05$.]{
  \includegraphics[width=0.3\textwidth]{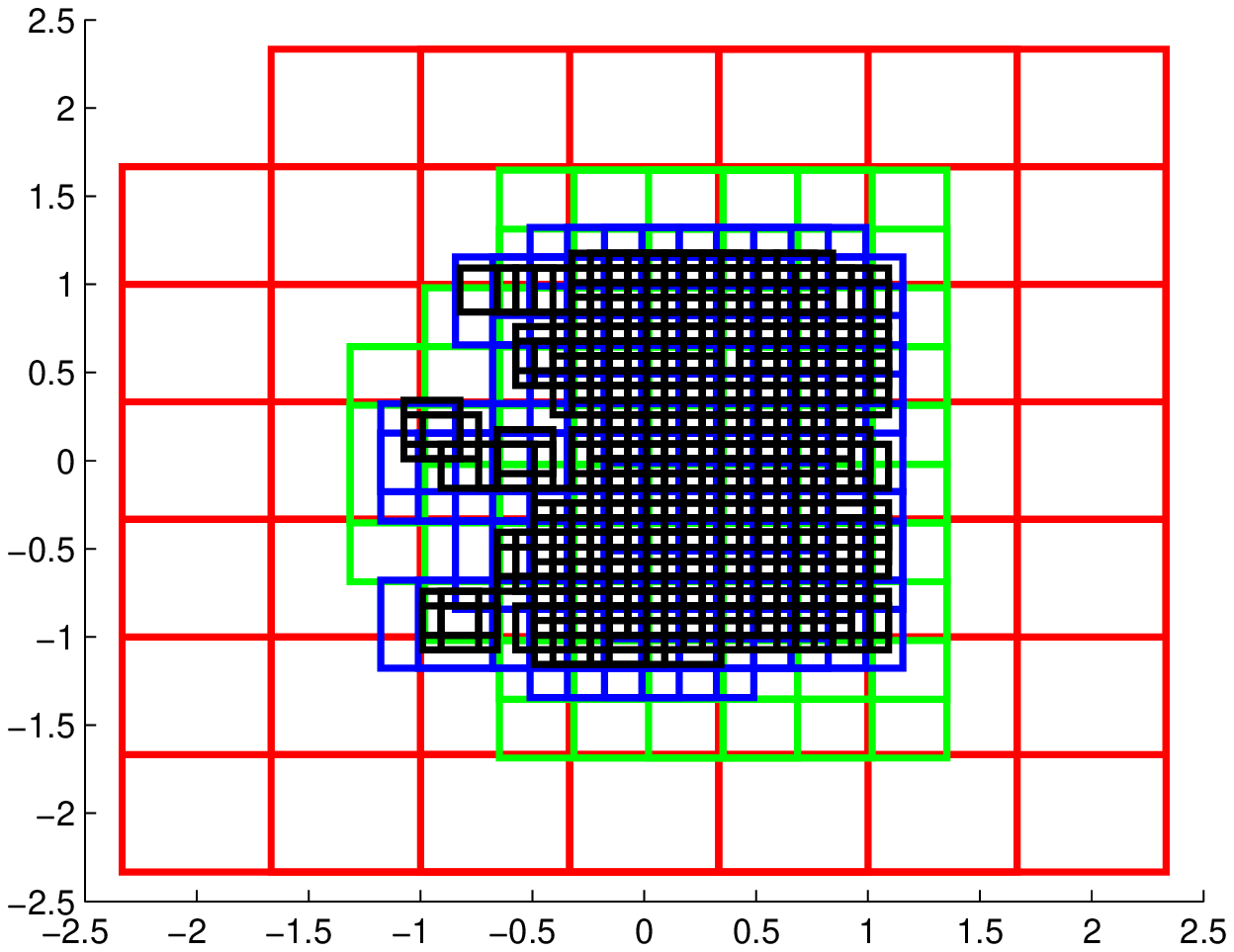}}
  \caption{The distribution of the support for $f_{3}$ with $J=4$.}
  \label{34}
\end{figure}

\begin{table}[htbp]
\center
\caption{The $l_{0}$ norm and the approximation errors for different parameters shown in Fig.\ref{34} for $f_{3}$ with $J=4$.}\label{Tab:34}
\begin{tabular}{cccccc}\hline
                   & $(l_{0}(\textbf{X}_{1}),l_{0}(\textbf{X}_{2}),l_{0}(\textbf{X}_{3}),l_{0}(\textbf{X}_{4}))$  & Error & RMS & Iterations  & Time(sec) \\ \hline
\cite{Lee}         &   (25, 63, 192, 674)     &   1.6373e-3    &  3.8077e-3   & 4      & 0.7348   \\
Fig.\ref{34}(a)    &   (16, 47, 105, 38)      &   1.7665e-3    &  5.6282e-3   & 573    & 8.3451   \\
Fig.\ref{34}(b)    &   (13, 31, 58, 10)        &   2.3442e-3    &  5.6422e-3   & 978    & 13.6259   \\
Fig.\ref{34}(c)    &   (7, 25, 96, 0)         &   3.6609e-3    &  6.4672-3   & 854    & 11.7498    \\
Fig.\ref{34}(d)    &   (0, 51, 36, 0)         &   4.3713e-3    &  8.7693e-3   & 1588   & 19.8770    \\
Fig.\ref{34}(e)    &   (14, 18, 19, 0)        &   4.6868e-2    &  5.3092e-2   & 1      & 5.2829    \\
Fig.\ref{34}(f)    &   (20, 38, 73, 119)      &   5.1103e-2    &  5.5647e-2   & 1      & 5.1456    \\
Fig.\ref{34}(g)    &   (22, 46, 68, 279)      &   6.9144e-2    &  7.0093e-2   & 1      & 4.3378    \\
Fig.\ref{34}(h)    &   (23, 34, 101, 328)     &   1.0138e-1    &  1.1552e-1   & 1      & 3.6595    \\
\hline
\end{tabular}
\end{table}

\begin{figure}[htbp]
\renewcommand{\figurename}{Fig.}
  \centering
  \subfigure[$J=3$ for MLASSO: $\lambda_{1}=\lambda_{2}=\lambda_{3}=0.001$.]{
  \includegraphics[width=0.3\textwidth]{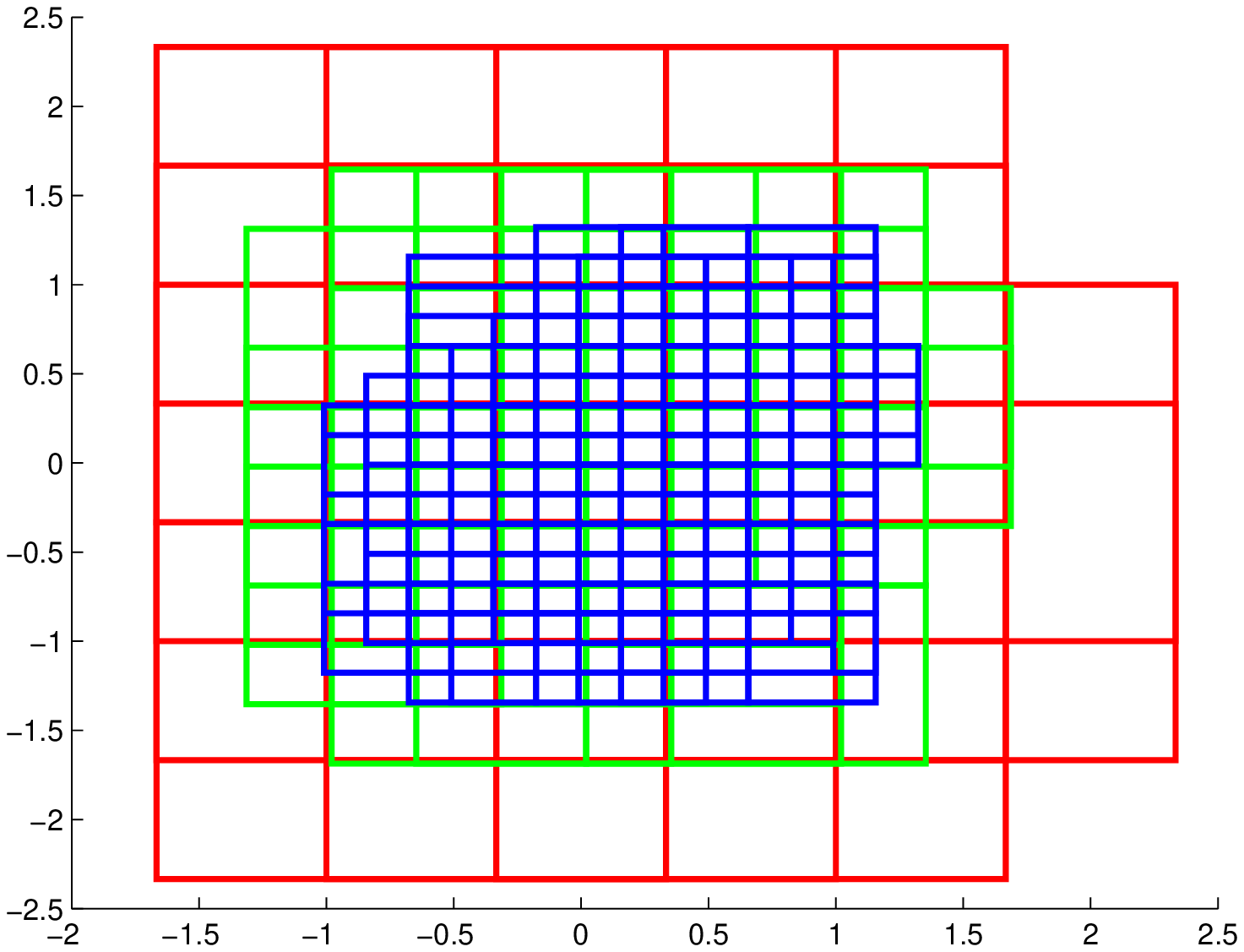}}
  \subfigure[$J=3$ for MLASSO: $\lambda_{1}=\lambda_{2}=\lambda_{3}=0.01$.]{
  \includegraphics[width=0.3\textwidth]{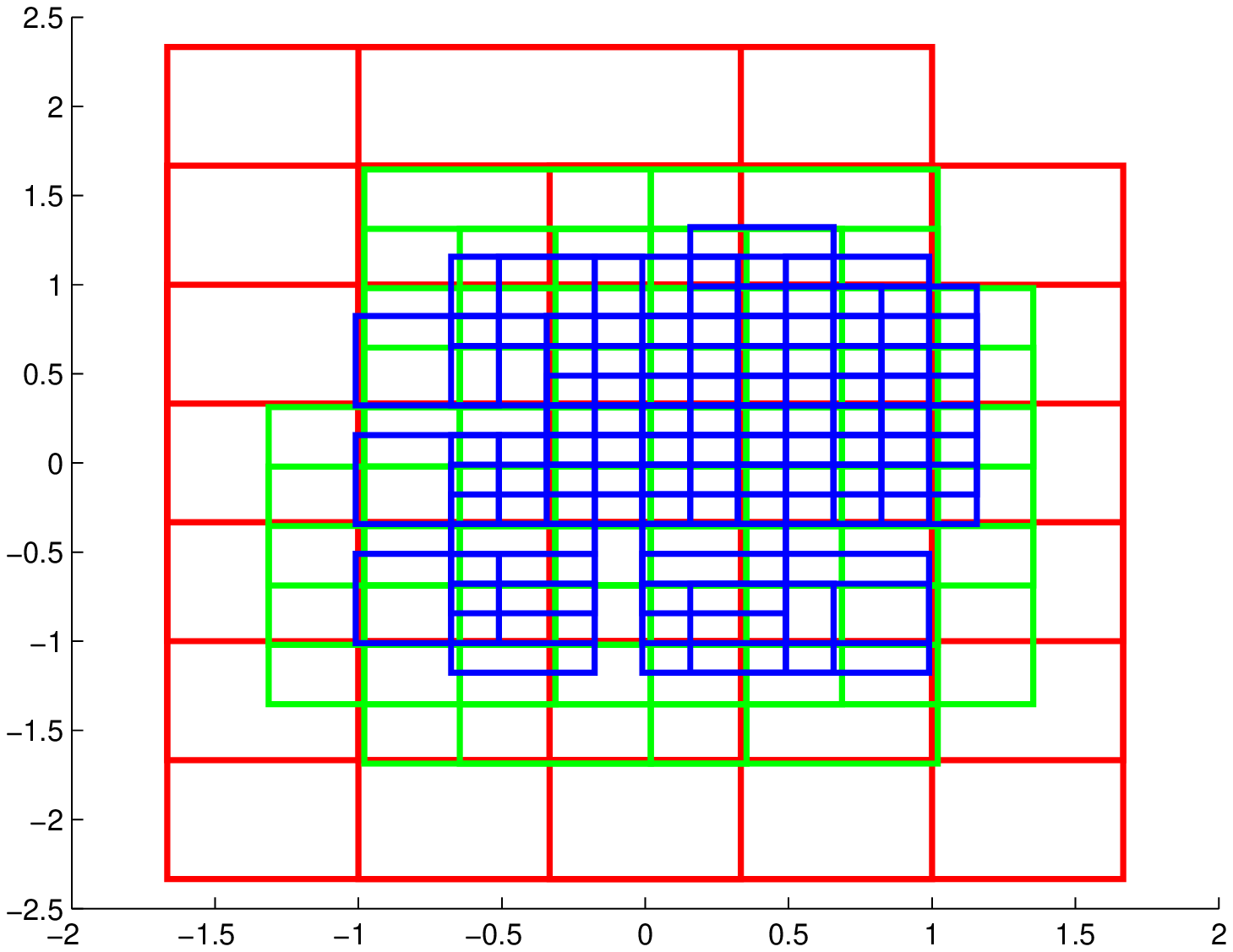}}
  \subfigure[$J=3$ for MLASSO: $\lambda_{1}=0.03,\lambda_{2}=0.01,\lambda_{3}=0.02$.]{
  \includegraphics[width=0.3\textwidth]{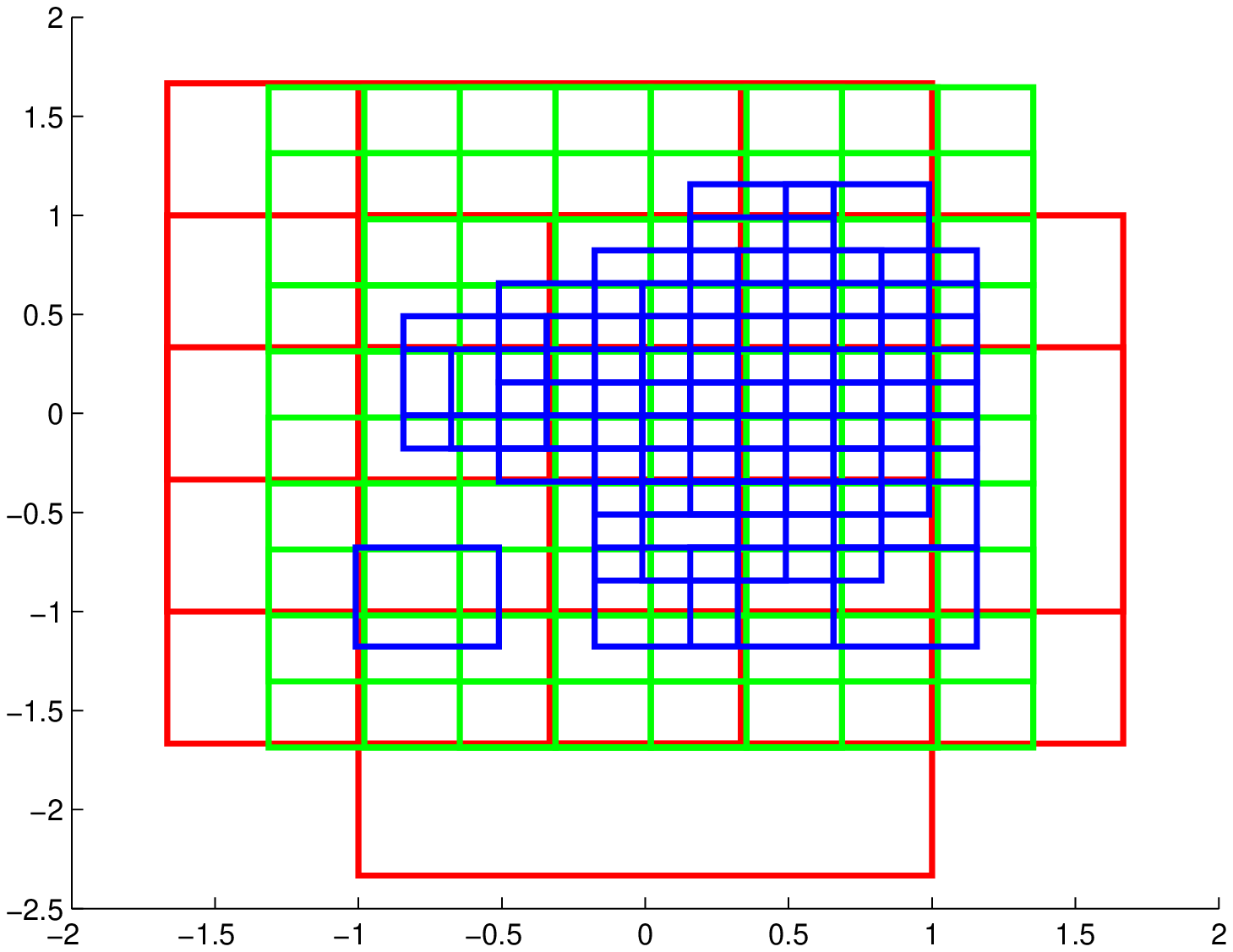}}
  \subfigure[$J=3$ for MLASSO: $\lambda_{1}=0.03, \lambda_{2}=0.02,\lambda_{3}=0.05$.]{
  \includegraphics[width=0.3\textwidth]{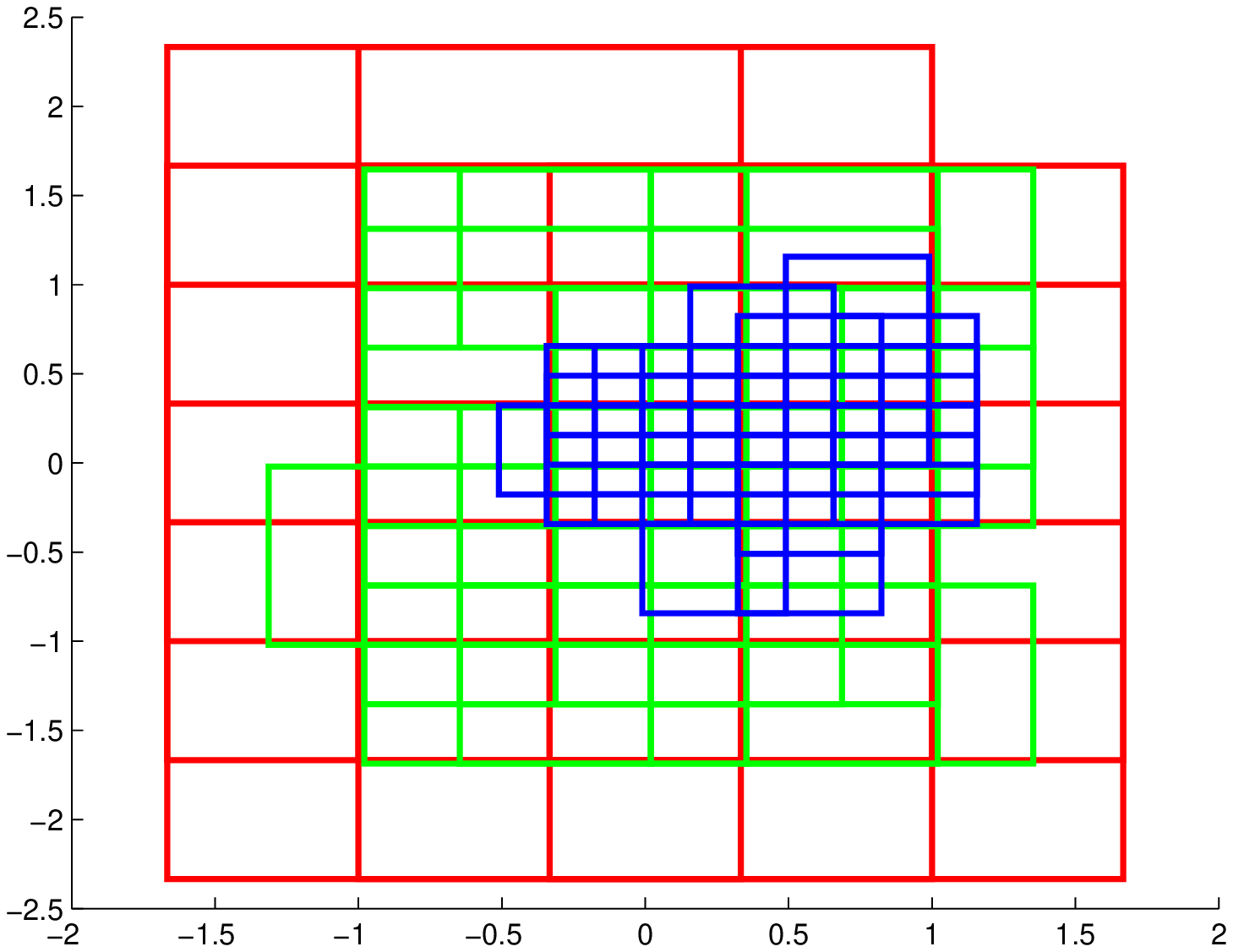}}
  \subfigure[$J=3$ for AGLASSO: $\lambda_{1}=\lambda_{2}=\lambda_{3}=0.001$.]{
  \includegraphics[width=0.3\textwidth]{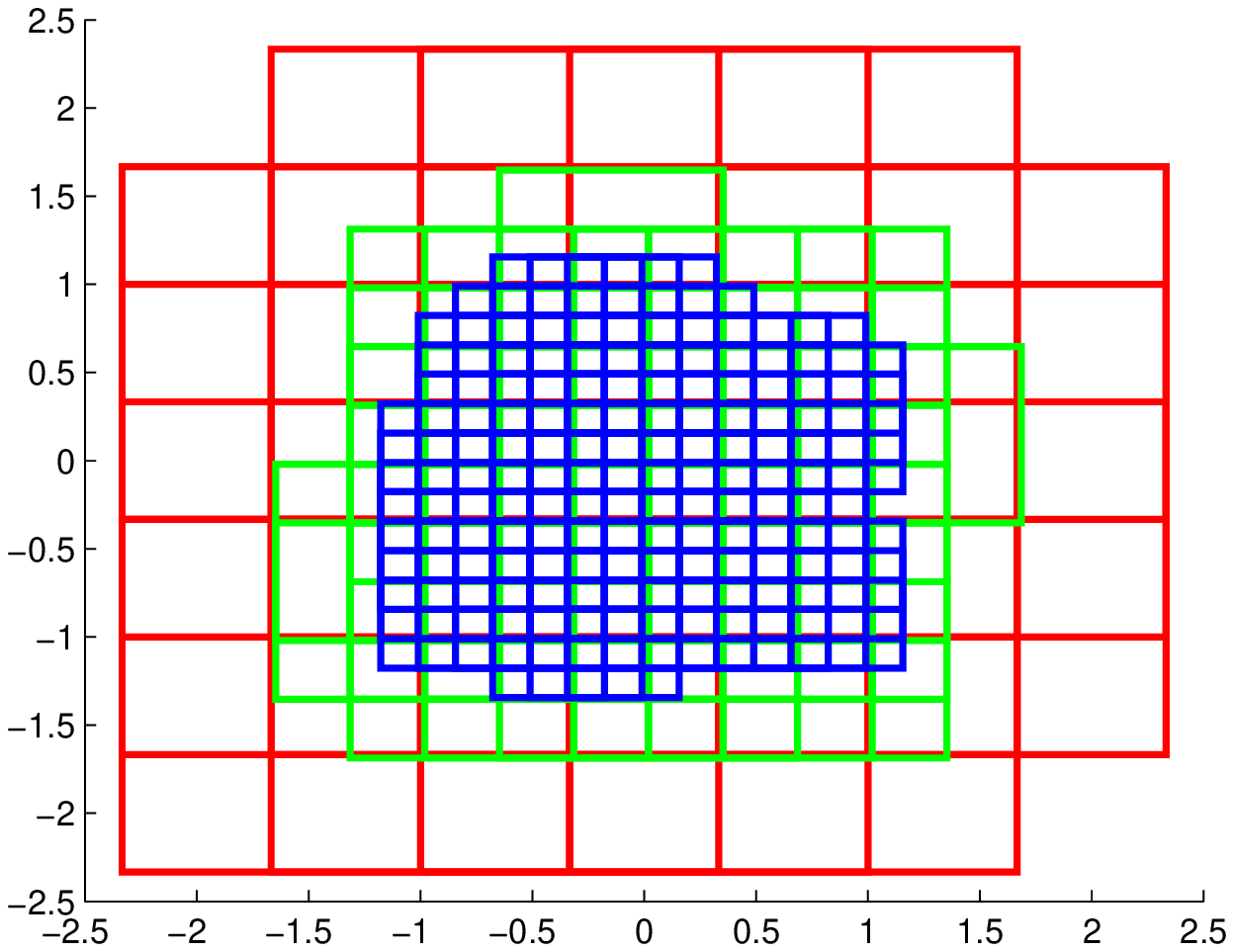}}
  \subfigure[$J=3$ for AGLASSO: $\lambda_{1}=\lambda_{2}=\lambda_{3}=0.01$.]{
  \includegraphics[width=0.3\textwidth]{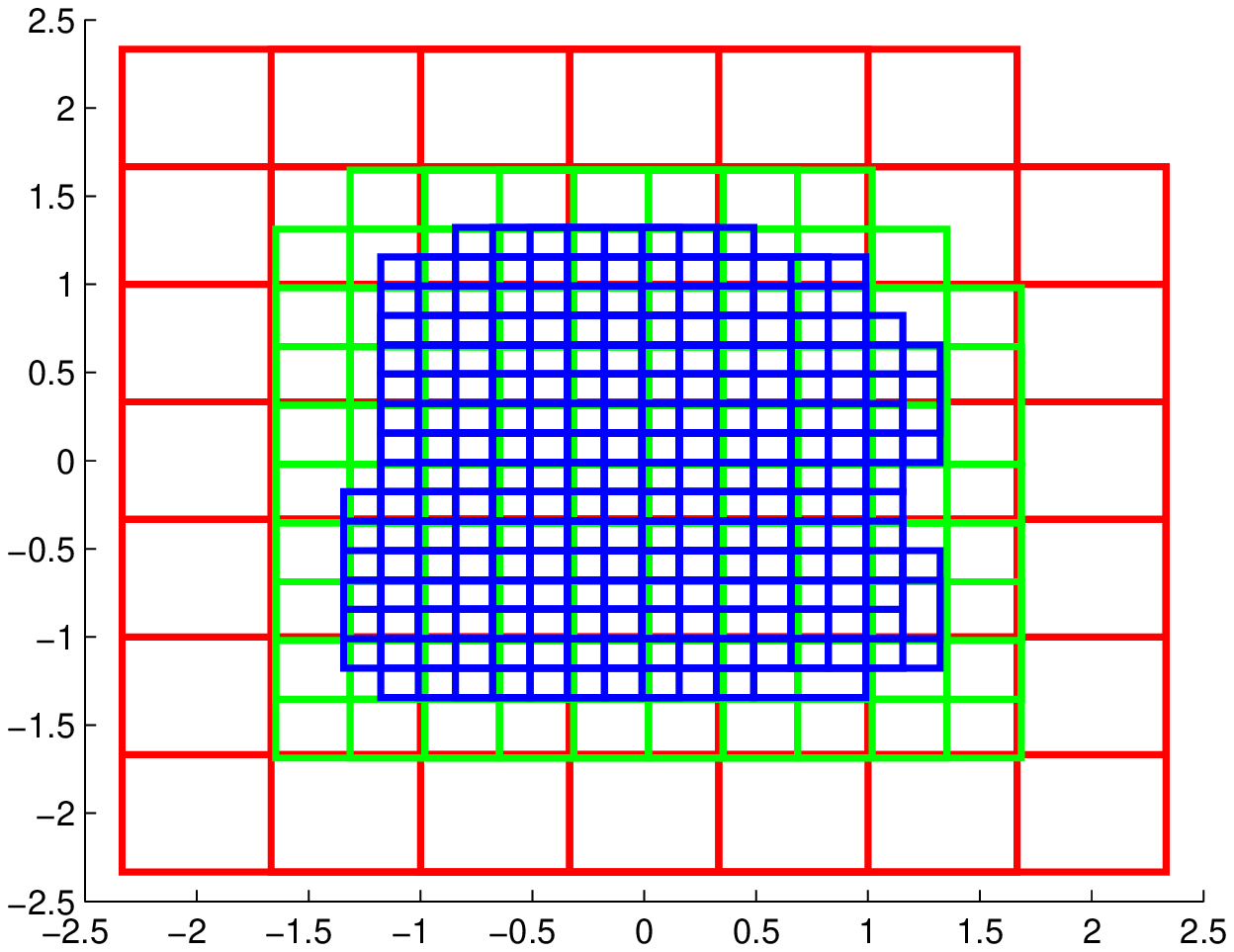}}
  \subfigure[$J=3$ for AGLASSO: $\lambda_{1}=0.03,\lambda_{2}=0.01,\lambda_{3}=0.02$.]{
  \includegraphics[width=0.3\textwidth]{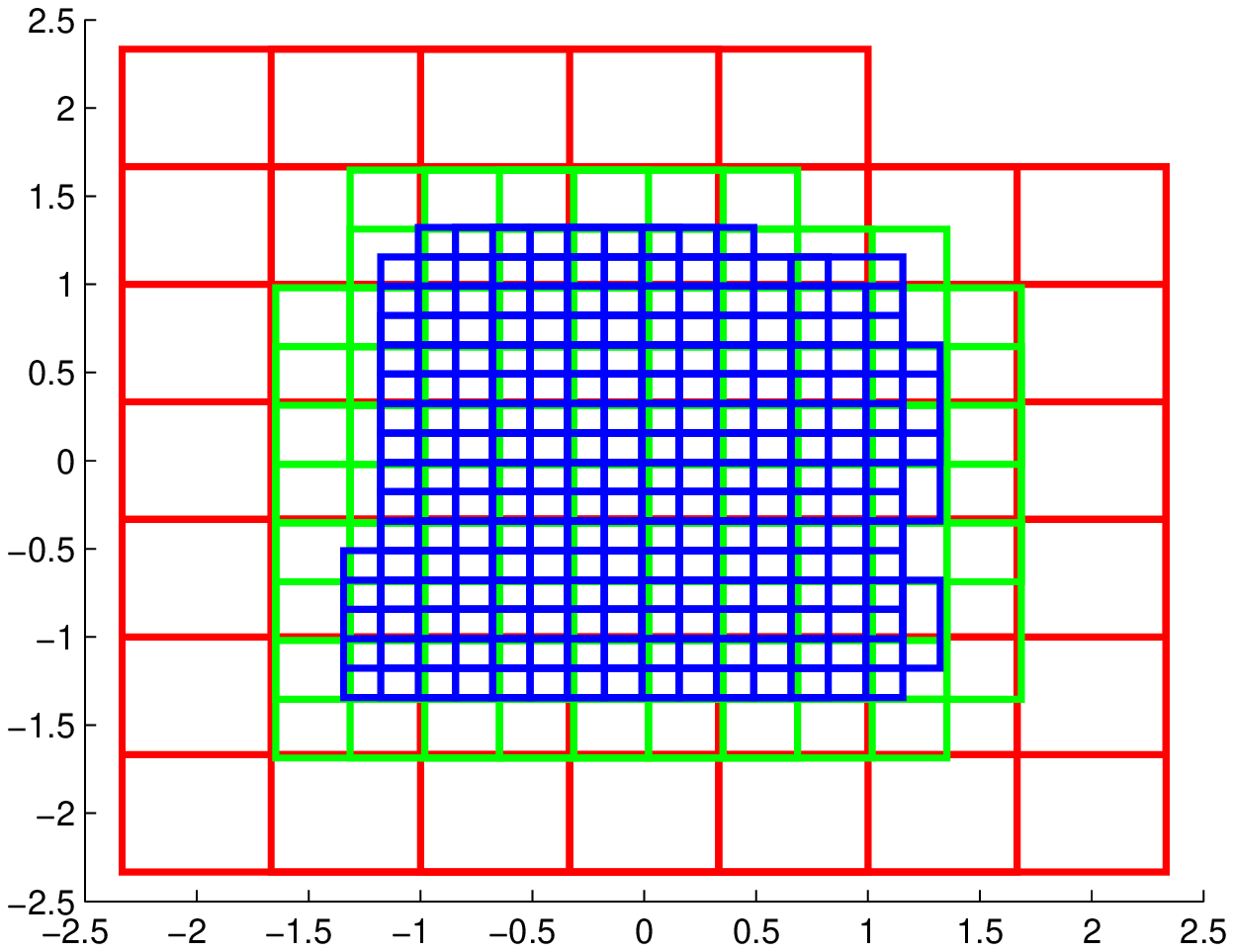}}
  \subfigure[$J=3$ for AGLASSO: $\lambda_{1}=0.03, \lambda_{2}=0.02,\lambda_{3}=0.05$.]{
  \includegraphics[width=0.3\textwidth]{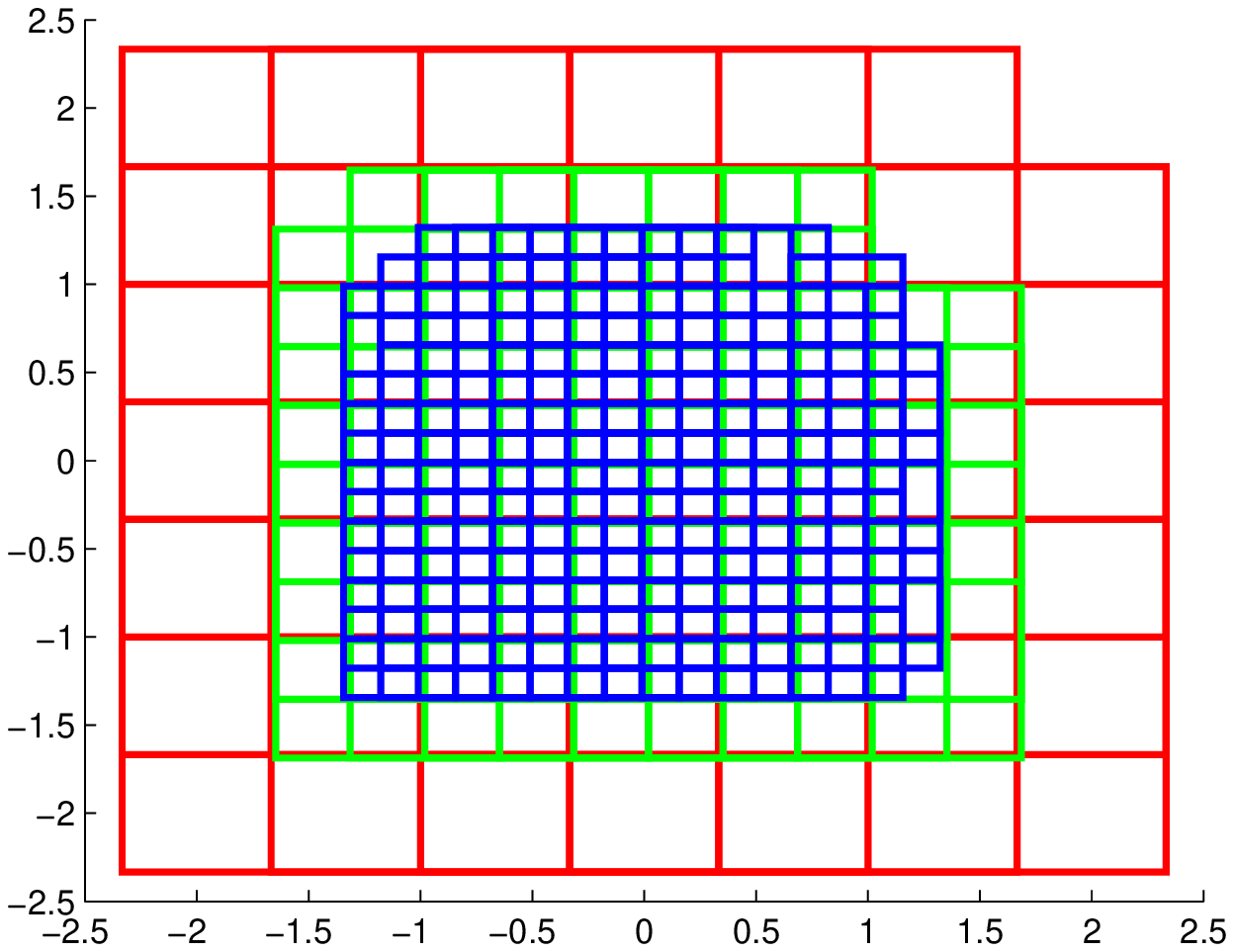}}
  \caption{The distribution of the support for $f_{4}$ with $J=3$.}
  \label{43}
\end{figure}

\begin{table}[htbp]
\center
\caption{The $l_{0}$ norm and the approximation errors for different parameters shown in Fig.\ref{43} for $f_{4}$ with $J=3$.}\label{Tab:43}
\begin{tabular}{cccccc}\hline
                   & $(l_{0}(\textbf{X}_{1}),l_{0}(\textbf{X}_{2}),l_{0}(\textbf{X}_{3}))$  & Error & RMS & Iterations  & Time(sec) \\ \hline
\cite{Lee}         &   (25, 62, 193)      &   3.2254e-3    &  5.0711e-3   & 3      & 0.0834   \\
Fig.\ref{43}(a)    &   (17, 41, 110)      &   3.7421e-3    &  6.3902e-3   & 19     & 1.1642   \\
Fig.\ref{43}(b)    &   (15, 34, 57)       &   4.0324e-3    &  6.2097e-3   & 938    & 1.7439    \\
Fig.\ref{43}(c)    &   (9, 38, 46)        &   4.5021e-3    &  8.1093e-3   & 1463   & 2.4762    \\
Fig.\ref{43}(d)    &   (15, 29, 33)       &   6.2043e-3    &  8.0072e-3   & 2230   & 2.7803    \\
Fig.\ref{43}(e)    &   (21, 45, 126)      &   2.0414e-1    &  6.0056e-1   & 1      & 0.7439    \\
Fig.\ref{43}(f)    &   (24, 60, 161)      &   2.1369e-1    &  3.1057e-1   & 1      & 0.5648    \\
Fig.\ref{43}(g)    &   (23, 57, 166)      &   3.2451e-1    &  2.2731e-1   & 1      & 0.4846    \\
Fig.\ref{43}(h)    &   (24, 58, 179)      &   5.3021e-1    &  2.1237e-1   & 1      & 0.4218    \\
\hline
\end{tabular}
\end{table}

\begin{table}[htbp]
\center
\caption{The $l_{0}$ norm and the approximation errors for different parameters shown in Fig.\ref{44} for $f_{4}$ with $J=4$.}\label{Tab:44}
\begin{tabular}{cccccc}\hline
                   & $(l_{0}(\textbf{X}_{1}),l_{0}(\textbf{X}_{2}),l_{0}(\textbf{X}_{3}),l_{0}(\textbf{X}_{4}))$  & Error & RMS & Iterations  & Time(sec) \\ \hline
\cite{Lee}         &   (25, 63, 194, 675)     &   4.1124e-4    &  3.2017e-3   & 4      & 2.1406   \\
Fig.\ref{44}(a)    &   (14, 34, 84, 108)      &   5.0121e-4    &  5.0512e-3   & 687    & 12.1132   \\
Fig.\ref{44}(b)    &   (15, 33, 48, 31)       &   1.18974e-3    &  4.6645e-3   & 1152   & 18.4120   \\
Fig.\ref{44}(c)    &   (13, 27, 52, 8)        &   4.1321e-3    &  7.1452e-3   & 1333   & 18.2461    \\
Fig.\ref{44}(d)    &   (12, 22, 70, 6)        &   4.2134e-3    &  7.6731e-3   & 2231   & 28.1121    \\
Fig.\ref{44}(e)    &   (22, 47, 126, 145)     &   1.5782e-1    &  1.7287e-1   & 1      & 5.1631    \\
Fig.\ref{44}(f)    &   (24, 58, 163, 491)     &   2.2056e-1    &  2.2142e-1   & 1      & 5.4315    \\
Fig.\ref{44}(g)    &   (24, 58, 164, 557)     &   3.0349e-1    &  3.4397e-1   & 1      & 4.3341    \\
Fig.\ref{44}(h)    &   (25, 60, 160, 524)     &   5.2570e-1    &  5.3329e-1   & 1      & 3.9764    \\
\hline
\end{tabular}
\end{table}

\begin{figure}[htbp]
\renewcommand{\figurename}{Fig.}
  \centering
  \subfigure[$J=4$ for MLASSO: $\lambda_{1}=\lambda_{2}=\lambda_{3}=\lambda_{4}=0.001$.]{
  \includegraphics[width=0.3\textwidth]{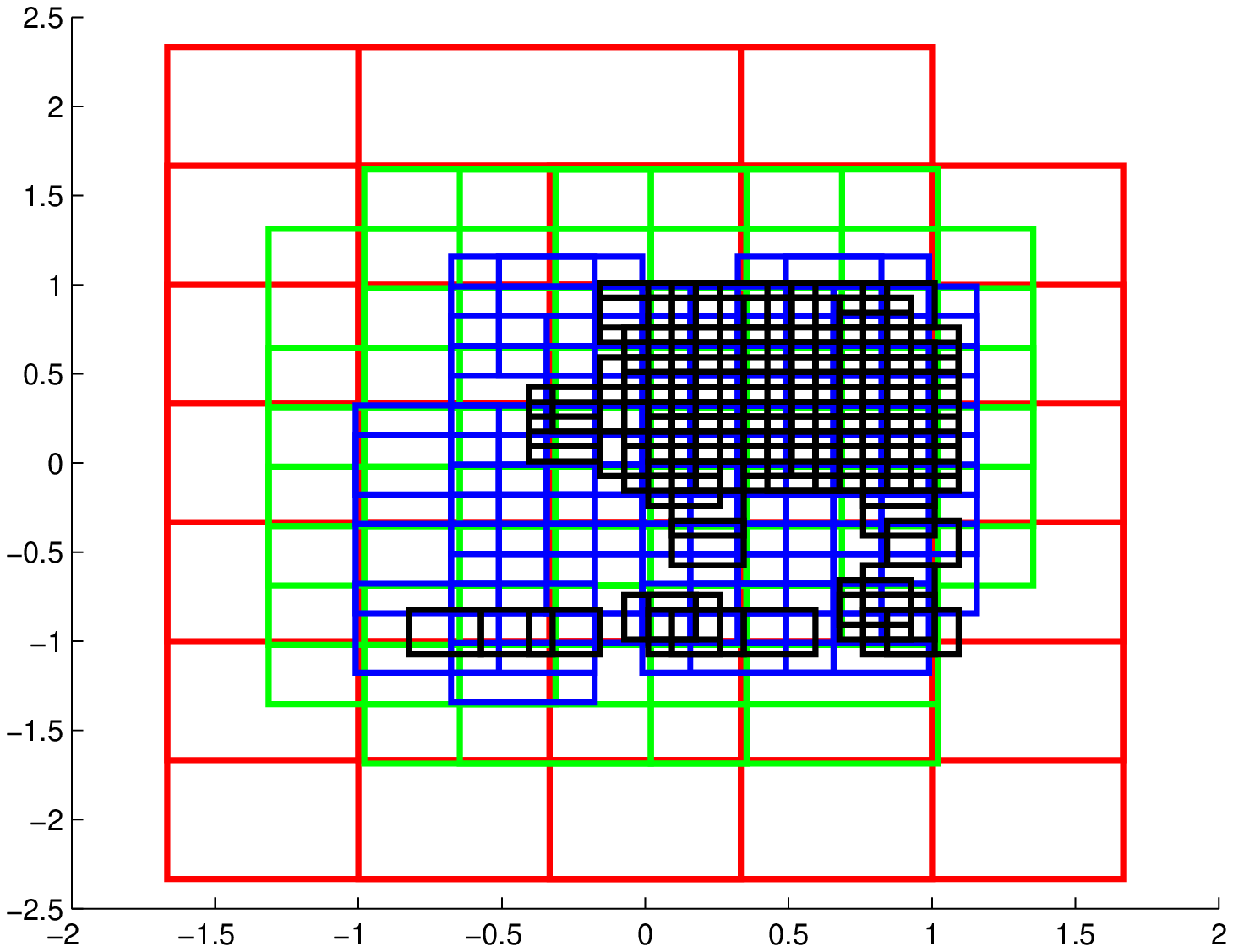}}
  \subfigure[$J=4$ for MLASSO: $\lambda_{1}=\lambda_{2}=\lambda_{3}=\lambda_{4}=0.01$.]{
  \includegraphics[width=0.3\textwidth]{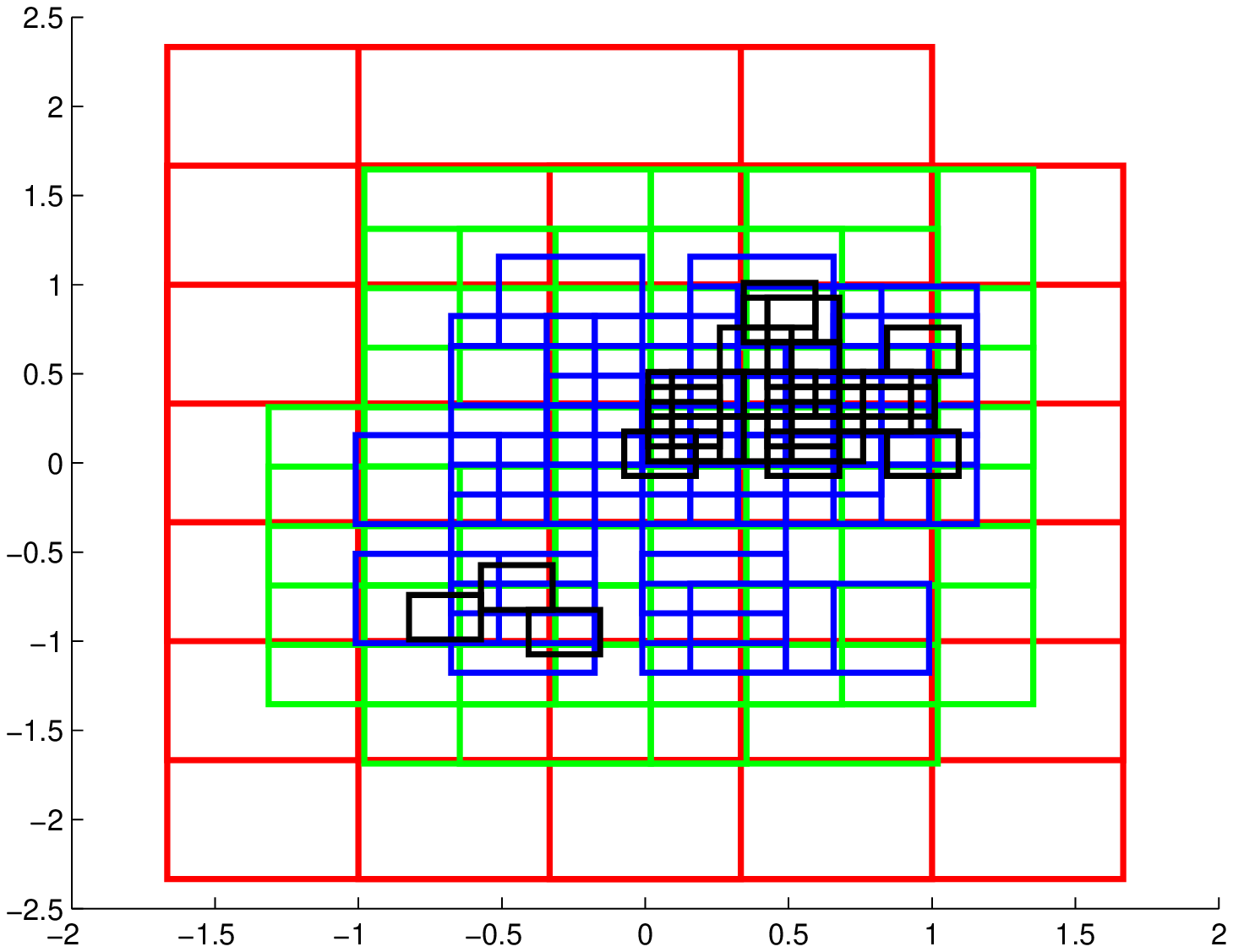}}
  \subfigure[$J=4$ for MLASSO: $\lambda_{1}=0.03,\lambda_{2}=0.02,\lambda_{3}=0.02,\lambda_{4}=0.04$.]{
  \includegraphics[width=0.3\textwidth]{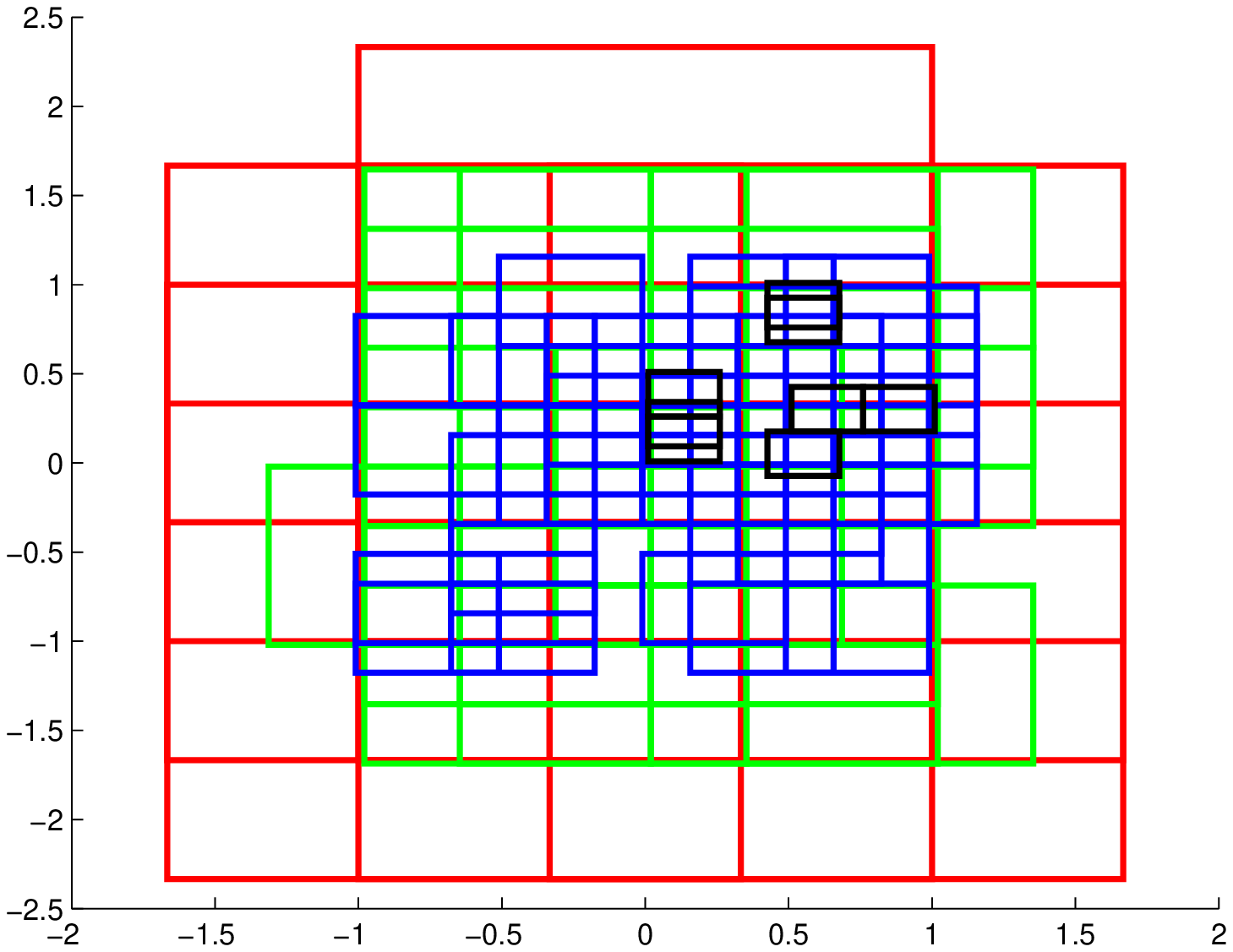}}
  \subfigure[$J=4$ for MLASSO: $\lambda_{1}=0.05,\lambda_{2}=0.03,\lambda_{3}=0.02,\lambda_{4}=0.04$.]{
  \includegraphics[width=0.3\textwidth]{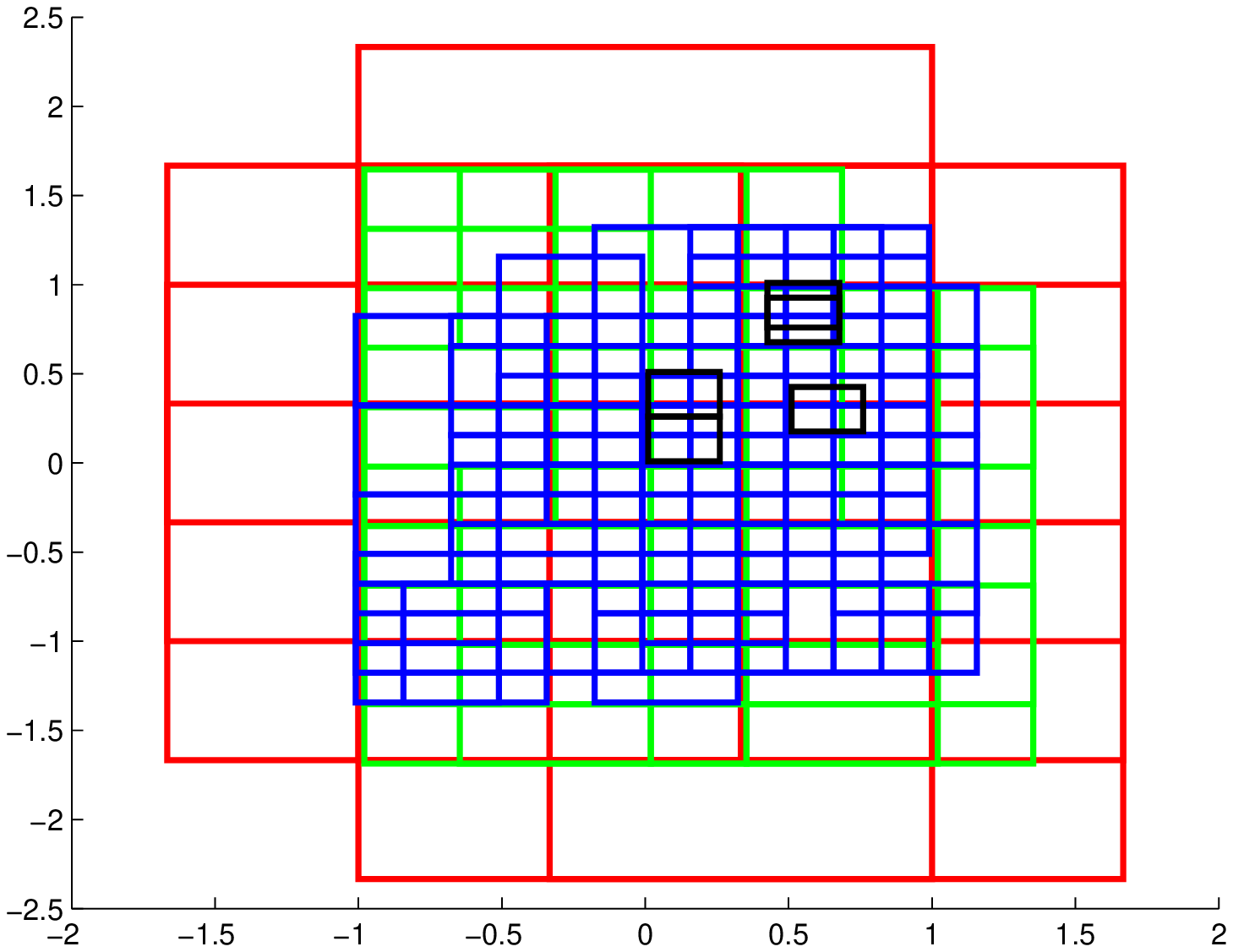}}
  \subfigure[$J=4$ for AGLASSO: $\lambda_{1}=\lambda_{2}=\lambda_{3}=\lambda_{4}=0.001$.]{
  \includegraphics[width=0.3\textwidth]{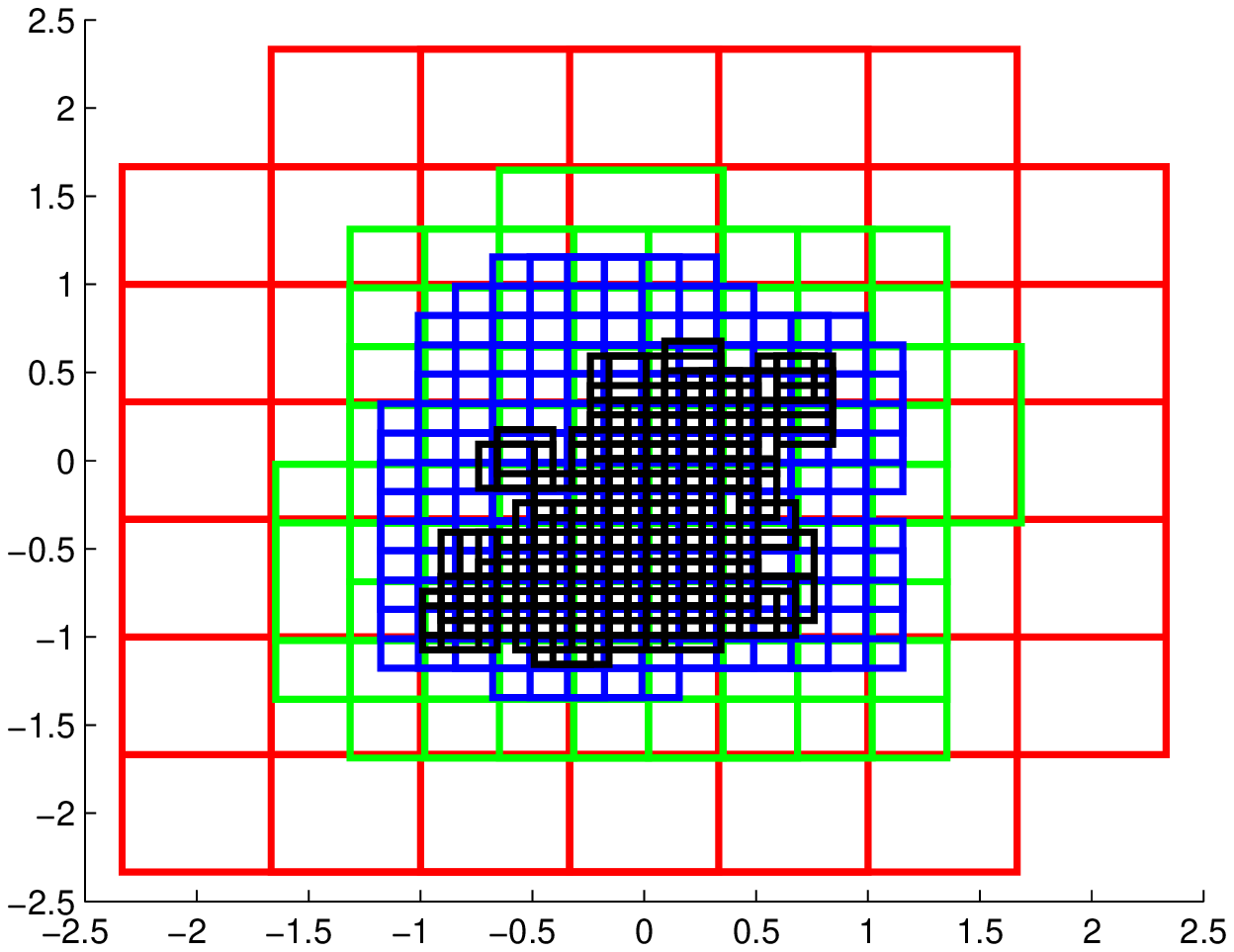}}
  \subfigure[$J=4$ for AGLASSO: $\lambda_{1}=\lambda_{2}=\lambda_{3}=\lambda_{4}=0.01$.]{
  \includegraphics[width=0.3\textwidth]{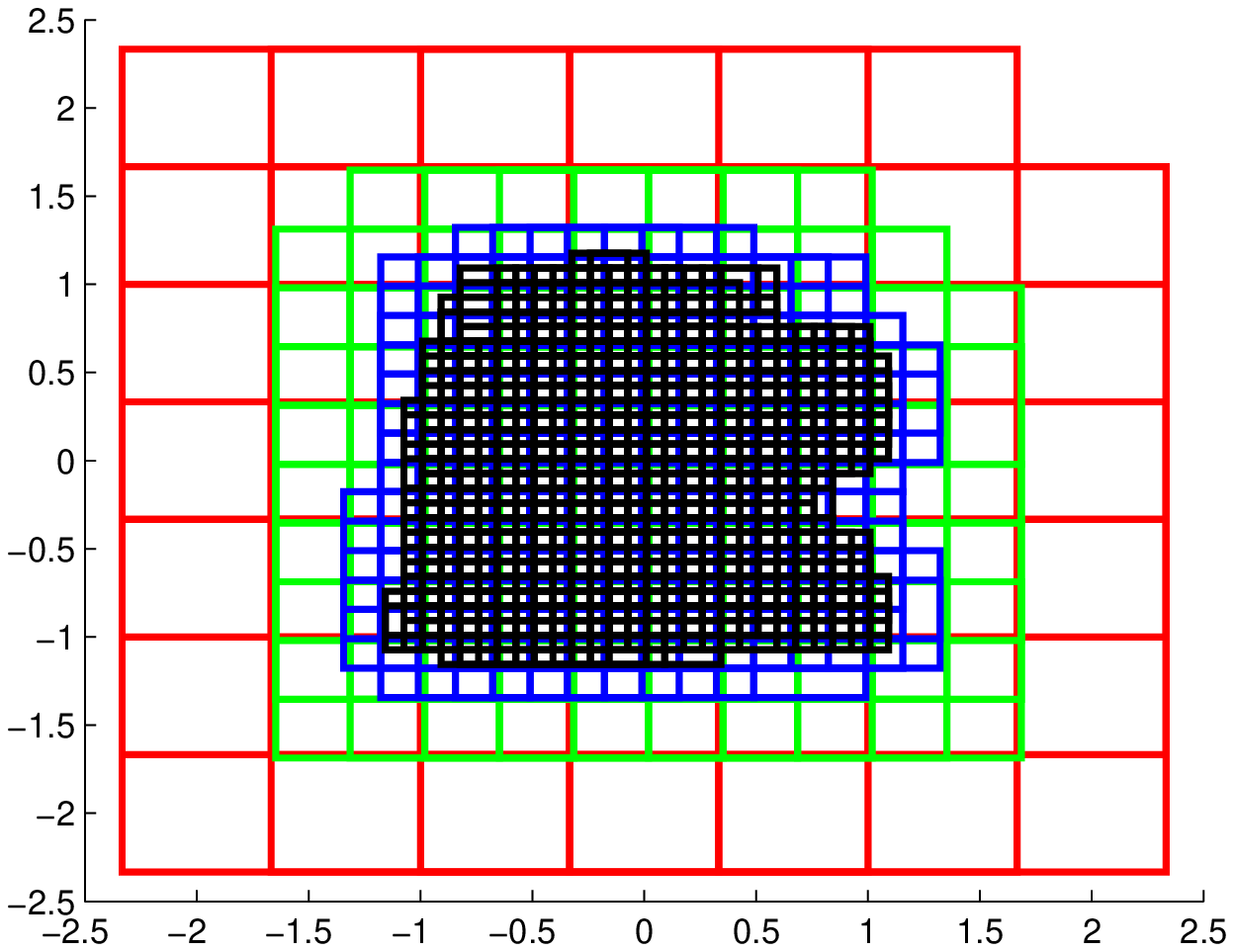}}
  \subfigure[$J=4$ for AGLASSO: $\lambda_{1}=0.03,\lambda_{2}=0.02,\lambda_{3}=0.02,\lambda_{4}=0.04$.]{
  \includegraphics[width=0.3\textwidth]{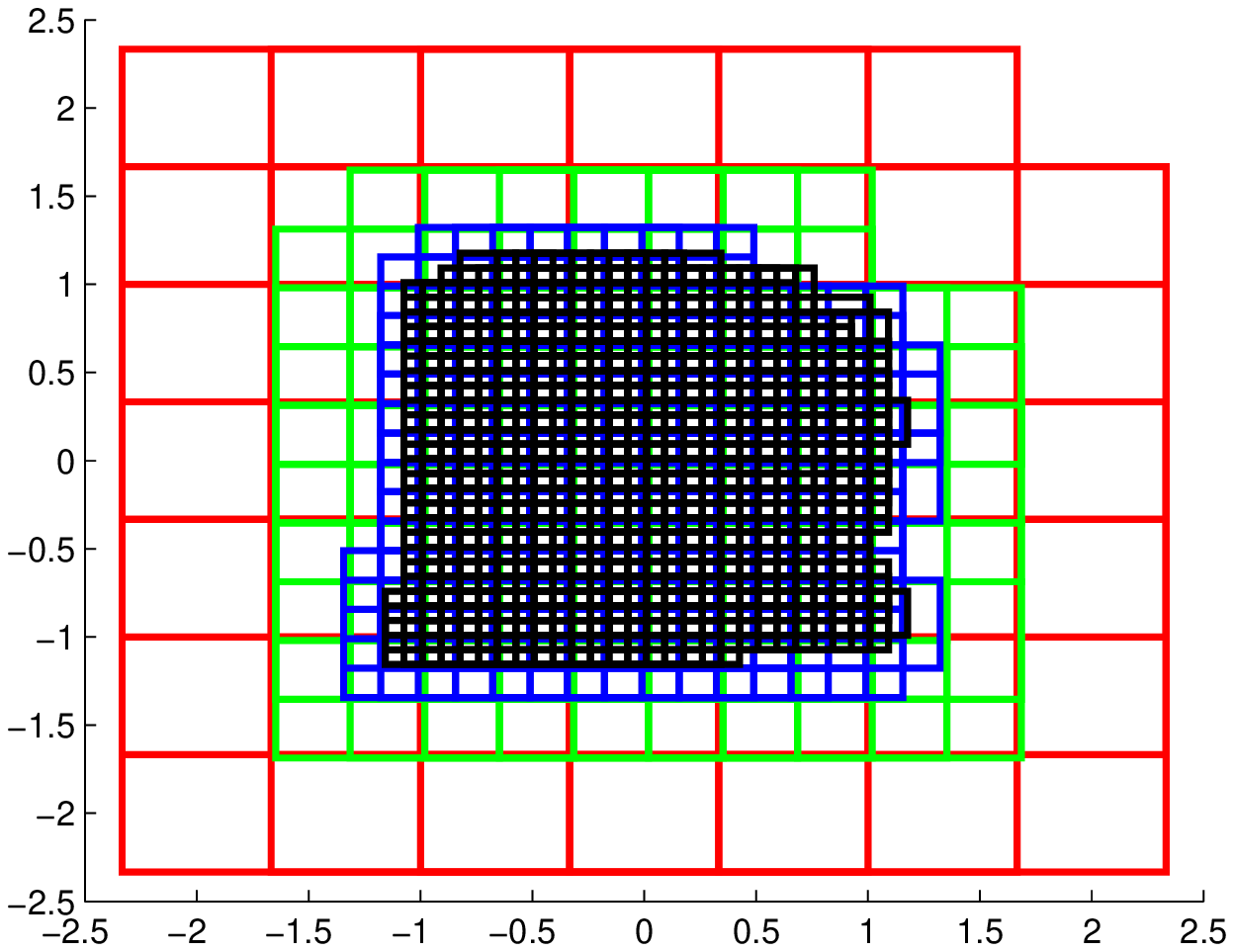}}
  \subfigure[$J=4$ for AGLASSO: $\lambda_{1}=0.05,\lambda_{2}=0.03,\lambda_{3}=0.02,\lambda_{4}=0.04$.]{
  \includegraphics[width=0.3\textwidth]{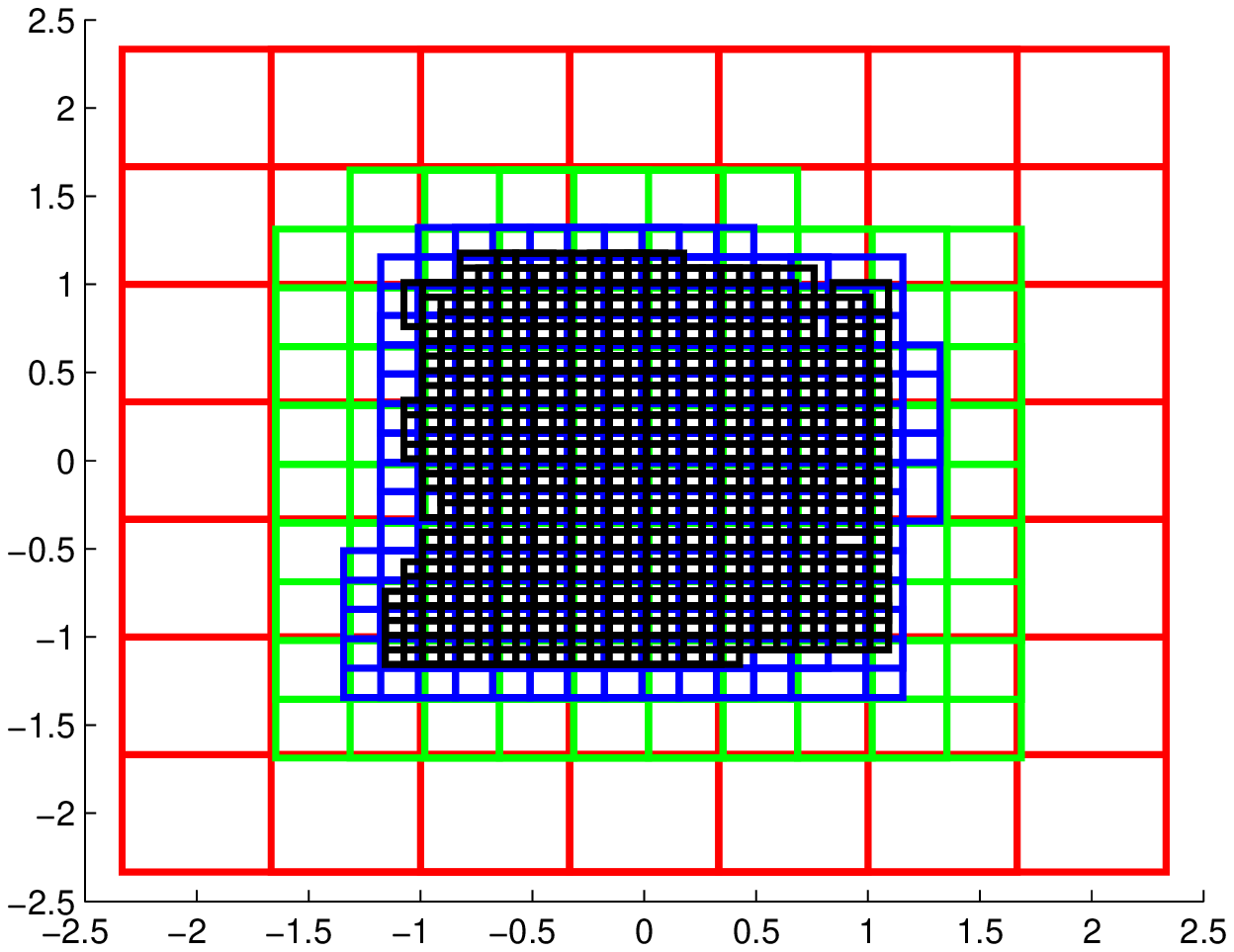}}
  \caption{The distribution of the support for $f_{4}$ with $J=4$. }
  \label{44}
\end{figure}

For the numerical implementation in this paper, we employ the 2D tensor product quadratic B-spline
as the function $\varphi$. We show that such a simple system with the Algorithm 1 can be used to
effectively reconstruct surface of sparse representation from a scattered data set. The choices of
the threshold $\sigma, \epsilon$ can be chosen according to the accuracy and sparsity. We consider
the four functions with $\sigma=10^{-3},\epsilon=10^{-4},N=900$ and the 900 scattered points are chosen randomly.
Moreover, we compare Algorithm 1 with the basic MBA algorithm presented in \cite{Lee} under the termination condition $|\triangle g_{j}|\leq O(10^{-3})$ and the AGLASSO model, both applying the same thresholding $\sigma=10^{-3}$ as Algorithm 1 to their results.

We experiment the above three methods for $J=3$ and $J=4$ starting from level 1 with $5\times 5$ biquadratic B-spline functions respectively.
Fig.\ref{11} shows the scattered data points and the corresponding approximations of the four functions with Algorithm 1,
the method in \cite{Lee} and the AGLASSO model respectively. Fig.2-9 illustrate the distribution of the support of the B-spline
functions with nonzero coefficients for $J=3$ and $J=4$ respectively. The red, green, blue and black rectangles denote
the support of the B-splines corresponding to $\textbf{X}_{1}, \textbf{X}_{2},\textbf{X}_{3}$ and $\textbf{X}_{4}$ respectively. Moreover, Tables 1-8 give the approximation accuracy, the iterations, the running time and the $l_{0}$ norm of the solution with different parameters for $f_{1}-f_{4}$, where the length of $\textbf{X}_{1}, \textbf{X}_{2},\textbf{X}_{3}$ and $\textbf{X}_{4}$ is 25, 64, 196 and 625 respectively.
In addition, all our calculations are done in Matlab on a laptop with Inter Core i7 (2.90GHZ) CPU and 8.0G RAM.

\begin{discussion}
The numerical results demonstrate that Algorithm 1 and the method in \cite{Lee} have almost
the same approximation errors, while Algorithm 1 obtains the sparse solution. Compared with the AGLASSO model,
Algorithm 1 provides the sparser solutions with less error, though more iterations and more time.

Through the first two steps of Algorithm 1, we can obtain an approximation solution of the MLASSO model with
no sparsity and the solution has some big components and some small ones which reflect different importance and contribution.
Then by step three, we throw out small components of the computed solution which means we keep only the important ones
with great contribution to the solution. Therefore, by choosing appropriate regularization parameters, the final solution
can indicate the important parts we are interested in and identify the important features within the selected levels simultaneously
of the exact surface since they are all determined by the big components. Experiment results verify this conclusion.

Taken together, Algorithm 1 can reconstruct the test functions reasonably with a sparse representation within a few levels by choosing some appropriate regularization parameters.

\end{discussion}

\section{Conclusion}
This paper presents an approach for scattered data fitting using the PSI space and the $l_{1}$ regularization.
It is concluded into the MLASSO model which allows us to balance the accuracy and the sparsity of the fitting surface.
The model can be solved using Algorithm 1 with the ADMM algorithm.
Numerical examples demonstrate that compared to the basic MBA algorithm in \cite{Lee} and the AGLASSO model,
the MLASSO model provides an efficient, sparse, flexible and reasonable solution.
Moreover, the distribution of the basis functions of the sparse solution can identify the regions
of the underlying surface where large fluctuations occur.

\Acknowledgements{This work was supported by the National Natural Science Foundation of China
(Nos.11526098, 11001037, 11290143, 11471066), the Research Foundation for Advanced Talents of
Jiangsu University (No.14JDG034), the Natural Science Foundation of Jiangsu Province (No.BK20160487)
and the Fundamental Research Funds for the Central Universities (DUT15LK44).}


\end{document}